\documentclass[a4paper, twoside]{article}

\usepackage[
  inner=2.5cm,
  textwidth=17.5cm,tmargin=2cm, bmargin=1cm,  includefoot,
]{geometry}
\usepackage{xfrac}
\usepackage{amsfonts}
\usepackage{amsmath}
\usepackage{amssymb}
\usepackage{graphicx}
\usepackage{ntheorem}

\usepackage[hidelinks]{hyperref}
\let\inf\relax \DeclareMathOperator*\inf{\vphantom{p}inf}
\let\liminf\relax \DeclareMathOperator*\liminf{\vphantom{p}liminf}
\let\lim\relax \DeclareMathOperator*\lim{\vphantom{p}lim}
\let\max\relax \DeclareMathOperator*\max{\vphantom{p}max}
\newtheorem{theorem}{Theorem}[section]

\newtheorem{corollary}[theorem]{Corollary}

\newtheorem{lemma}[theorem]{Lemma}

\theorembodyfont{\upshape}\newtheorem{example}[theorem]{Example}
\newtheorem{definition}[theorem]{Definition}
\newtheorem{exercise}[theorem]{Exercise}

\newcommand{\be}{\begin{equation}}
\newcommand{\ee}{\end{equation}}
\addcontentsline{toc}{section}{Abstract}

\usepackage[hidelinks]{hyperref}

\newcommand\independent{\protect\mathpalette{\protect\independenT}{\perp}}
\def\independenT#1#2{\mathrel{\rlap{$#1#2$}\mkern4mu{#1#2}}}

\begin{document}
\title{Weak Convergence of Probability Measures }
\date{}

\author{{\sc Serik Sagitov}, {\small Chalmers University of Technology and Gothenburg University}}
%\date{}
\maketitle

\begin{abstract}
This text contains my lecture notes for the graduate course ``Weak Convergence" given in September-October 2013 and then in March-May 2015. 
The course is based on the book Convergence of Probability Measures by Patrick Billingsley, partially covering Chapters 1-3, 5-9, 12-14, 16, as well as appendices. 
In this text the formula label $(*)$ operates locally.  
%The visible theorem labels often show the theorem numbers in the book, labels involving PM refer to the other book by Billingsley - "Probability and Measure". 

I am grateful to Timo Hirscher whose numerous valuable suggestions helped me to improve earlier versions of these notes.  Last updated: \today.
\end{abstract}

%\begin{quote}
%\scriptsize{
% I know of scarcely anything so apt to impress the imagination as the wonderful form of cosmic order expressed by the ``Law of Frequency of Error". The law would have been personified by the Greeks and deified, if they had known of it. It reigns with serenity and in complete self-effacement, amidst the wildest confusion. The huger the mob, and the greater the apparent anarchy, the more perfect is its sway. It is the supreme law of Unreason. Whenever a large sample of chaotic elements are taken in hand and marshaled in the order of their magnitude, an unsuspected and most beautiful form of regularity proves to have been latent all along.
% \hfill Galton F. (1889) Natural Inheritance.}
%\end{quote}

\tableofcontents

\newpage
\section*{Introduction}
\begin{figure}
\centering
\includegraphics[width=12cm,height=4cm]{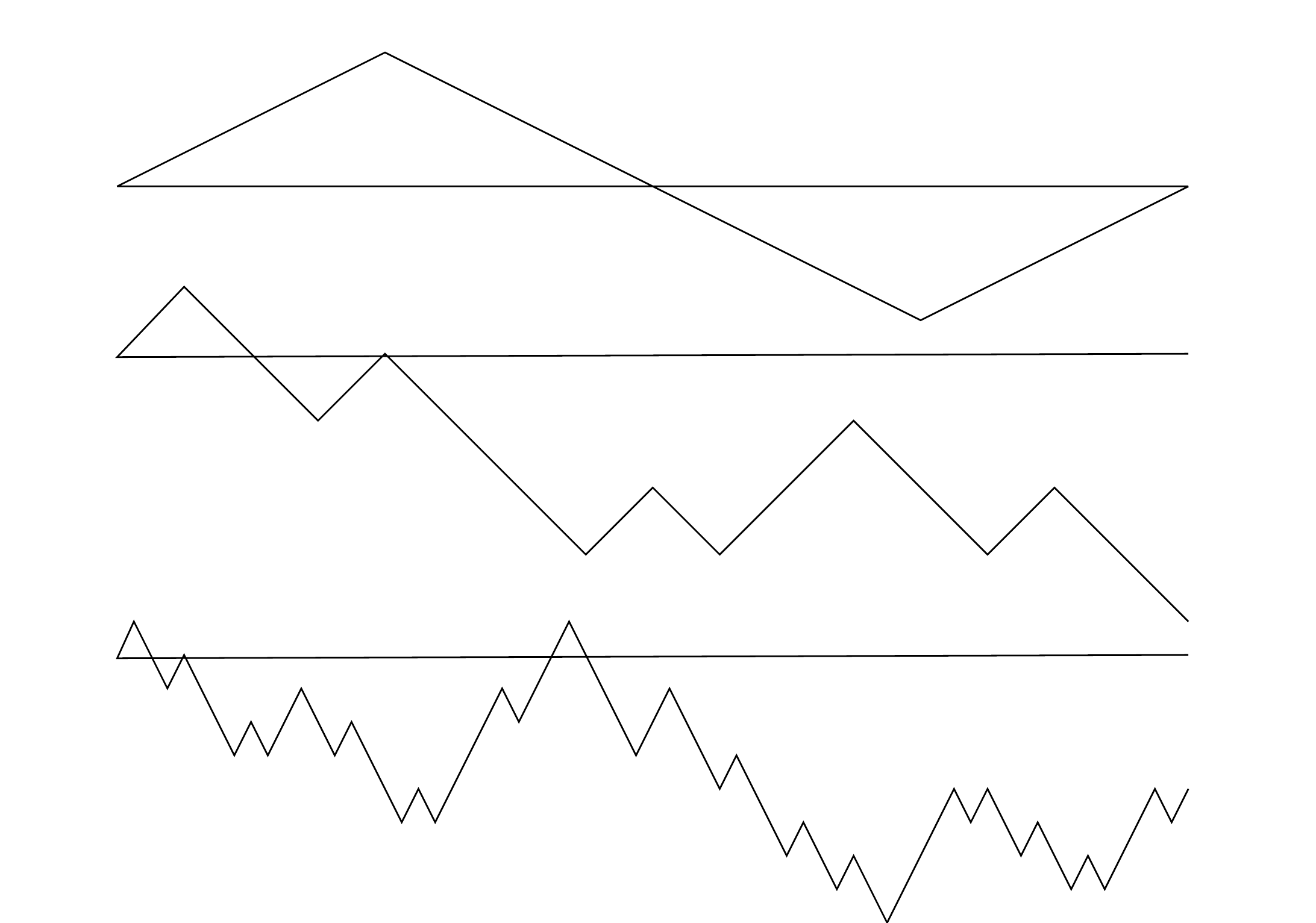}
\caption{Scaled symmetric simple random walk $X^n_t(\omega)$ for a fixed $\omega\in\Omega$ and $n=4,16,64$.}
\label{frw}
\end{figure}

Throughout these lecture notes we use the following notation
\[\Phi(z)={1\over\sqrt{2\pi}}\int_{-\infty}^z e^{-{u^2/2}}du.\]
Consider a symmetric simple random walk  $S_n=\xi_1+\ldots+\xi_n$ with $\mathbb P(\xi_i=1)=\mathbb P(\xi_i=-1)=1/2$. The random sequence $S_n$ has no limit in the usual sense. However, by de Moivre's theorem (1733),  
\[\mathbb P(S_n\le z\sqrt{n})\to\Phi(z) \mbox{ as }n\to\infty  \mbox{ for any }z\in\boldsymbol R.\]
This is an example of convergence in distribution ${S_n\over \sqrt{n}}\Rightarrow Z$ to a normally distributed random variable. 
Define a sequence of stochastic processes $X^n=(X^n_t)_{t\in[0,1]}$ by linear interpolation between its values $X^n_{i/n}(\omega)={S_i(\omega)\over\sigma\sqrt n}$ at the points $t=i/n$, see Figure \ref{frw}. 
The much more powerful functional CLT claims convergence in distribution  towards the Wiener process 
$X^n\Rightarrow W$.

%What we need is to expand the notion of convergence in distribution to random elements of metric spaces like $\boldsymbol C[0,1]$ and $\boldsymbol D[0,1]$.

This course deals with weak convergence of probability measures on Polish spaces $({\boldsymbol S},\mathcal S)$. For us, the principal examples of Polish spaces (complete separable metric spaces) are
 
 the space $\boldsymbol C=\boldsymbol C[0,1]$ of continuous trajectories $x:[0,1]\to\boldsymbol R$ (Section \ref{secC}),
 
 the space $\boldsymbol D=\boldsymbol D[0,1]$ of cadlag trajectories $x:[0,1]\to\boldsymbol R$ (Section \ref{secD}),

 the space $\boldsymbol D[0,\infty)$ of cadlag trajectories $x:[0,\infty)\to\boldsymbol R$ (Section \ref{secDi}).\\

To prove the functional CLT $X^n\Rightarrow W$,  we have to check that $\mathbb E f(X^n)\to\mathbb Ef(W)$ for all bounded continuous functions $f:{\boldsymbol C}[0,1]\to \boldsymbol R$, which is not practical to do straightforwardly. Instead, one starts with the finite-dimensional distributions 
$$(X^n_{t_1},\ldots,X^n_{t_k})\Rightarrow (W_{t_1},\ldots,W_{t_k}).$$
To prove the weak convergence of the finite-dimensional distributions, it is enough to check the convergence of moment generating functions, thus allowing us to focus on a special class of continuous functions $f_{\lambda_1,\ldots,\lambda_k}: \boldsymbol R^k\to \boldsymbol R$, where $\lambda_i\ge0$ and 
\[f_{\lambda_1,\ldots,\lambda_k}(x_1,\ldots,x_k)=\exp(\lambda_1x_1+\ldots+\lambda_kx_k).\] 

For the weak convergence in the infinite-dimensional space $\boldsymbol C[0,1]$, the usual additional step is to verify tightness of the distributions of the family of processes $(X^n)$. Loosely speaking, tightness means that no probability mass escapes to infinity.
By Prokhorov theorem  (Section \ref{secP}), tightness implies relative compactness, which means that each subsequence of $X^n$ contains a further subsequence converging weakly. Since all possible limits have the finite-dimensional distributions of $W$, we conclude that all subsequences converge to the same limit $W$, and by this we establish the convergence  $X^n\Rightarrow W$. 

This approach makes it crucial to find tightness criteria in $\boldsymbol C[0,1]$, $\boldsymbol D[0,1]$, and then in $\boldsymbol D[0,\infty)$.

\section{The Portmanteau and mapping theorems}

\subsection{Metric spaces}

Consider a metric space ${\boldsymbol S}$  with metric  $\rho(x,y)$. For subsets $A\subset {\boldsymbol S}$, denote the closure by $A^-$, the interior by $A^\circ$, and the boundary by $\partial A=A^--A^\circ$. We write 
\[\rho(x,A)=\inf\{\rho(x,y):y\in A\},\qquad A^\epsilon=\{x:\rho(x,A)<\epsilon\}.\]

\begin{definition}
Open balls $B(x,r)=\{y\in {\boldsymbol S}: \rho(x,y)<r\}$ form a {\it  base} for ${\boldsymbol S}$: each {\it open set} in ${\boldsymbol S}$ is a union of open balls. Complements to the open sets are called  {\it closed sets}. The   {\it Borel  $\sigma$-algebra} $\mathcal S$ is formed from the open and closed sets  in ${\boldsymbol S}$ using the operations of countable intersection, countable union, and set difference.
\end{definition}

\begin{definition}
A collection $\mathcal A$  of ${\boldsymbol S}$-subsets is called a  $\pi$-system if it is closed under intersection, that is if $A,B\in\mathcal A$, then $A\cap B\in\mathcal A$. 
We say that $\mathcal L$ is a $\lambda$-system if: (i) ${\boldsymbol S}\in\mathcal L$, (ii) $A\in\mathcal L$ implies $A^c \in\mathcal L$, (iii) for any sequence of disjoint sets  $A_n\in\mathcal L$,  $\cup_nA_n\in \mathcal L$.
\end{definition}

\begin{theorem}\label{Dyn} Dynkin's $\pi$-$\lambda$  lemma.
If $\mathcal A$ is a  $\pi$-system such that  $\mathcal A\subset\mathcal L$, where $\mathcal L$ is a $\lambda$-system, then 
 $\sigma(\mathcal A)\subset\mathcal L$, where $\sigma(\mathcal A)$ is the  $\sigma$-algebra generated by $\mathcal A$. 
\end{theorem}
\begin{definition}
A metric space ${\boldsymbol S}$ is called
{\it separable} if it contains a countable dense subset. It is called 
{\it complete} if every Cauchy (fundamental) sequence has a limit lying in ${\boldsymbol S}$.  A complete separable metric space is called a {\it Polish} space. 
\end{definition}

Separability is a topological property, while completeness is a property of the metric and not of the topology. 
\begin{definition}
 An {\it open cover} of $A\subset {\boldsymbol S}$ is a class of open sets whose union contains $A$. 
\end{definition}
\begin{theorem}\label{M3} 
 These three conditions are equivalent:
 
 (i) ${\boldsymbol S}$ is separable,
 
 (ii) ${\boldsymbol S}$ has a countable base (a class of open sets such that each open set is a union of sets in the class),
 
 (iii) Each open cover of each subset of ${\boldsymbol S}$ has a countable subcover. 
\end{theorem}
\begin{theorem}\label{M3'} 
 Suppose that the subset $M$ of ${\boldsymbol S}$ is separable.
 
 (i) There is a countable class $\mathcal A$ of open sets with the property that, if $x\in G\cap M$ and $G$ is open, then $x\in A\subset A^-\subset G$ for some $A\in\mathcal A$.
 
 (ii) Lindel\"{o}f property. Each open cover of $M$ has a countable subcover. 
\end{theorem}

\begin{definition}
 A set $K$ is called {\it compact} if each open cover of $K$ has a finite subcover.
 A set $A\subset {\boldsymbol S}$  is called {\it relatively compact} if each sequence in $A$ has a convergent subsequence the limit of which may not lie in $A$. 
\end{definition}

\begin{theorem}\label{M5}  
Let $A$ be a subset of a metric space ${\boldsymbol S}$. The following three conditions are equivalent:

(i) $A^-$ is compact,

(ii) $A$ is relatively compact,

(iii) $A^-$ is complete and $A$ is totally bounded (that is for any $\epsilon>0$, $A$ has a finite $\epsilon$-net the points of which are not required to lie in $A$).
 
\end{theorem}

\begin{theorem}\label{M10}  
Consider two metric spaces $({\boldsymbol S},\rho)$ and $({\boldsymbol S}',\rho')$ and maps $h,h_n:{\boldsymbol S}\to{\boldsymbol S}'$. If $h$ is continuous, then it is measurable $\mathcal S/\mathcal S'$. If each $h_n$ is measurable $\mathcal S/\mathcal S'$, and if $h_nx\to hx$ for every $x\in {\boldsymbol S}$, then $h$ is also measurable $\mathcal S/\mathcal S'$.
 
\end{theorem}

\subsection{Convergence in distribution and weak convergence}

\begin{definition}\label{p7}
 Let  $P_n, P$ be probability measures on $({\boldsymbol S},\mathcal S)$. We say $P_n\Rightarrow  P$ {\it weakly converges} as $n\to \infty$ if for any bounded continuous function $f:{\boldsymbol S}\to \boldsymbol R$
\[\int_{{\boldsymbol S}}f(x)P_n(dx)\to\int_{{\boldsymbol S}}f(x)P(dx), \quad n\to\infty. \]
\end{definition}
\begin{definition}
 Let  $X$ be a $({\boldsymbol S},\mathcal S)$-valued random element defined on the probability space $(\Omega,\mathcal F,\mathbb P)$. 
We say that a probability measure $P$ on ${\boldsymbol S}$ is the probability distribution of $X$ if $P(A)=\mathbb P(X\in A)$ for all $A\in\mathcal S$.
\end{definition}
\begin{definition}\label{p25}
 Let  $X_n,X$ be $({\boldsymbol S},\mathcal S)$-valued random elements defined on the probability spaces $(\Omega_n,\mathcal F_n,\mathbb P_n)$, $(\Omega,\mathcal F,\mathbb P)$. 
We say $X_n$ converge in distribution to $X$ as $n\to \infty$ and write $X_n\Rightarrow X$, if for any bounded continuous function $f:{\boldsymbol S}\to \boldsymbol R$,
\[\mathbb E_n(f(X_n))\to\mathbb E(f(X)), \quad n\to\infty. \]
This is equivalent to the weak convergence $P_n\Rightarrow  P$ of the respective probability distributions.
\end{definition}

\begin{example}
 The function $f(x)=1_{\{x\in A\}}$ is bounded but not continuous, therefore if $P_n\Rightarrow P$, then $P_n(A)\to P(A)$ does not always hold.
For ${\boldsymbol S}={{\boldsymbol R}}$, the function $f(x)=x$ is continuous but not  bounded, therefore if $X_n\Rightarrow X$, then $\mathbb E_n(X_n)\to\mathbb E(X)$ does not always hold.
\end{example}

\begin{definition}
 Call $A\in\mathcal S$ a {\it $P$-continuity set} if $P(\partial A)=0$.
\end{definition}

\begin{theorem}\label{2.1} Portmanteau's theorem. The following five statements are equivalent.

(i) $P_n\Rightarrow P$.

(ii) $\int f(x)P_n(dx)\to\int f(x)P(dx)$ for all bounded uniformly continuous $f:{\boldsymbol S}\to \boldsymbol R$.

(iii) $\limsup_{n\to\infty}P_n(F)\le P(F)$ for all closed $F\in\mathcal S$.

(iv) $\liminf_{n\to\infty}P_n(G)\ge P(G)$ for all open $G\in\mathcal S$.

(v) $P_n(A)\to P(A)$ for all $P$-continuity sets  $A$.
\end{theorem}
Proof. (i) $\to$ (ii) is trivial. 

\noindent (ii) $\to$ (iii). For a closed $F\in\mathcal S$ put 
$$g(x)=(1-\epsilon^{-1}\rho(x,F))\vee0.$$ 
This function is bounded and uniformly continuous since $|g(x)-g(y)|\le\epsilon^{-1}\rho(x,y)$. Using 
\[
1_{\{x\in F\}}\le g(x)\le 1_{\{x\in F^\epsilon\}},%\mbox{ where }F^\epsilon=\{x: \rho(x,F)<\epsilon\},
\] 
we derive (iii) from (ii):
$$\limsup_{n\to\infty}P_n(F)\le\limsup_{n\to\infty}\int g(x)P_n(dx)=\int g(x)P(dx)\le P(F^\epsilon)\to P(F),\quad \epsilon\to0.$$ 

\noindent (iii) $\to$ (iv) follows by complementation.

\noindent (iii) + (iv) $\to$ (v). If $P(\partial A)=0$, then then the leftmost and rightmost probabilities coincide:
\begin{align*}
  P(A^-)&\ge \limsup_{n\to\infty} P_n(A^-)\ge \limsup_{n\to\infty} P_n(A)\\
  &\ge \liminf_{n\to\infty}P_n(A)\ge \liminf_{n\to\infty}P_n(A^\circ)\ge P(A^\circ).
\end{align*}
\noindent (v) $\to$ (i). By linearity we may assume that the bounded continuous function $f$ satisfies $0\le f\le1$. Then putting $A_t=\{x:f(x)>t\}$ we get
\[\int_{\boldsymbol S} f(x)P_n(dx)=\int_0^1 P_n(A_t)dt \to\int_0^1 P(A_t)dt =\int_{\boldsymbol S} f(x)P(dx). \]
Here the convergence follows from (v) since $f$ is continuous, implying that 
$\partial A_t=\{x:f(x)=t\}$,
and since $\{x:f(x)=t\}$ are $P$-continuity sets except for countably many $t$. We also used the bounded convergence theorem.
\begin{example}
 Let $F(x)=\mathbb P(X\le x)$. Then $X_n=X+n^{-1}$ has distribution $F_n(x)=F(x-n^{-1})$. As $n\to\infty$, $F_n(x)\to F(x-)$, so convergence only occurs at continuity points.
\end{example}
\begin{corollary}\label{1.2}
A single sequence of probability measures can not weakly converge to each of two different limits. 
\end{corollary}
Proof. It suffices to prove that if $\int_{{\boldsymbol S}}f(x)P(dx)=\int_{{\boldsymbol S}}f(x)Q(dx)$ for all bounded, uniformly continuous functions $f:{\boldsymbol S}\to \boldsymbol R$, then $P=Q$. Using the 
bounded, uniformly continuous functions $g(x)=(1-\epsilon^{-1}\rho(x,F))\vee0$ we get
$$P(F)\le\int_{{\boldsymbol S}}g(x)P(dx)=\int_{{\boldsymbol S}}g(x)Q(dx)\le Q(F^\epsilon).$$
Letting $\epsilon\to0$ it gives for any closed set $F$, that $P(F)\le Q(F)$ and by symmetry we conclude that $P(F)=Q(F)$. It follows that $P(G)=Q(G)$ for all open sets $G$.

It remains to use regularity of any probability measure $P$: if $A\in\mathcal S$ and $\epsilon>0$, then there exist a closed set $F_\epsilon$ and an open set $G_\epsilon$ such that $F_\epsilon\subset A\subset G_\epsilon$ and $P(G_\epsilon-F_\epsilon)<\epsilon$. To this end we denote by $\mathcal G_P$ the class of $\mathcal S$-sets with the just stated property. If $A$ is closed, we can take $F=A$ and $G=F^\delta$, where $\delta$ is small enough. Thus all closed sets belong to 
$\mathcal G_P$, and we need to show that $\mathcal G_P$ forms a $\sigma$-algebra. Given $A_n\in\mathcal G_P$, choose closed sets $F_n$ and open sets $G_n$ such that $F_n\subset A_n\subset G_n$ and $P(G_n-F_n)<2^{-n-1}\epsilon$. If $G=\cup_{n}G_n$ and $F=\cup_{n\le n_0}F_n$ with $n_0$ chosen so that $P(\cup_{n}F_n-F)<\epsilon/2$, then $F\subset \cup_{n}A_n\subset G$ and $P(G-F)<\epsilon$. Thus $\mathcal G_P$ is closed under the formation of countable unions. Since it is closed under complementation, $\mathcal G_P$ is a $\sigma$-algebra.

\begin{theorem}\label{2.7}
 Mapping theorem. Let $X_n$ and $X$ be random elements of a metric space ${\boldsymbol S}$. Let $h:{\boldsymbol S}\to {\boldsymbol S}'$ be a $\mathcal S/\mathcal S'$-measurable mapping and $D_h$ be the set of its discontinuity points. If $X_n\Rightarrow X$ and $\mathbb P(X\in D_h)=0$, then $h(X^n)\Rightarrow h(X)$.
 
 In other terms, if $P_n\Rightarrow P$ and $P(D_h)=0$, then $P_nh^{-1}\Rightarrow Ph^{-1}$.
\end{theorem}
Proof. We show first that $D_h$ is a Borel subset of ${\boldsymbol S}$. For any pair $(\epsilon,\delta)$ of positive rationals, the set
\[A_{\epsilon\delta}=\{x\in {\boldsymbol S}: \mbox{ there exist $y,z\in {\boldsymbol S}$ such that }\rho(x,y)<\delta,\rho(x,z)<\delta,\rho'(hy,hz)\ge \epsilon\}\]
is open. Therefore, $D_h=\cup_\epsilon\cap_\delta A_{\epsilon\delta}\in\mathcal S$. Now, for each $F\in\mathcal S'$,
\begin{align*}
\limsup_{n\to\infty} P_n(h^{-1}F)&\le\limsup_{n\to\infty} P_n( (h^{-1}F)^-)\le P((h^{-1}F)^-)\\
&\le P( h^{-1}(F^-)\cup D_h)= P( h^{-1}(F^-)).
\end{align*}
To see that $(h^{-1}F)^-\subset h^{-1}(F^-)\cup D_h$ take an element $x\in(h^{-1}F)^-$. There is a sequence $x_n\to x$ such that $h(x_n)\in F$, and therefore, either $h(x_n)\to h(x)$ or $x\in D_h$. By the Portmanteau theorem, the last chain of inequalities implies $P_nh^{-1}\Rightarrow Ph^{-1}$.

\begin{example}
Let $P_n\Rightarrow P$. If $A$ is a $P$-continuity set and $h(x)=1_{\{x\in A\}}$, then by the mapping theorem, $P_nh^{-1}\Rightarrow Ph^{-1}$.
 \end{example}
 
 \subsection{Convergence in probability and in total variation. Local limit theorems}

\begin{definition}
 Suppose $X_n$ and $X$ are random elements of ${\boldsymbol S}$ defined on the same probability space. If $\mathbb P(\rho(X_n,X)<\epsilon)\to1$ for each positive $\epsilon$, we say $X_n$ converge to $X$ {\it in probability} and write $X_n\stackrel{\rm P}{\to} X$. 
\end{definition}
\begin{exercise}\label{inP}
Convergence in probability $X^n\stackrel{\rm P}{\to} X$ is equivalent to the weak convergence $\rho(X^n,X)\Rightarrow0$. Moreover, $(X^n_1,\ldots,X^n_k)\stackrel{\rm P}{\to} (X_1,\ldots,X_k)$ if and only if $X^n_i\stackrel{\rm P}{\to} X_i$ for all $i=1,\ldots,k$.
\end{exercise}
\begin{theorem}\label{3.2}
 Suppose $(X_n,X_{u,n})$ are random elements of ${\boldsymbol S}\times {\boldsymbol S}$. If $X_{u,n}\Rightarrow Z_u$ as $n\to\infty$ for any fixed $u$, and $Z_u\Rightarrow X$ as $u\to\infty$, and
 \[\lim_{u\to\infty}\limsup_{n\to\infty}\mathbb P(\rho(X_{u,n},X_n)\ge\epsilon)=0,\mbox{ for each positive }\epsilon,\]
 then $X_n\Rightarrow X$.
\end{theorem}
Proof. Let $F\in\mathcal S$ be closed and define $F_\epsilon$ as the set $\{x:\rho(x,F)\le\epsilon\}$. Then
\begin{align*}
 \mathbb P(X_n\in F)&=\mathbb P(X_n\in F,X_{u,n}\notin F_\epsilon)+\mathbb P(X_n\in F,X_{u,n}\in F_\epsilon)\\
 &\le\mathbb P(\rho(X_{u,n},X_n)\ge\epsilon)+\mathbb P(X_{u,n}\in F_\epsilon).
\end{align*}
Since $F_\epsilon$ is also closed and $F_\epsilon\downarrow F$ as $\epsilon\downarrow 0$, we get
\begin{align*}
 \limsup_{n\to\infty}\mathbb P(X_n\in F)&\le\limsup_{\epsilon\to0}\limsup_{u\to\infty}\limsup_{n\to\infty}\mathbb P(X_{u,n}\in F_\epsilon)\\
 &\le\limsup_{\epsilon\to0}\mathbb P(X\in F_\epsilon)=\mathbb P(X\in F).
\end{align*}
\begin{corollary}\label{3.1}
 Suppose $(X_{n},Y_n)$ are random elements of ${\boldsymbol S}\times {\boldsymbol S}$. If $Y_{n}\Rightarrow X$ as $n\to\infty$ and $\rho(X_n,Y_n)\Rightarrow0$, then $X_n\Rightarrow X$. Taking $Y_n\equiv X$, we conclude that convergence in probability implies convergence in distribution.
\end{corollary}

\begin{definition}
Convergence in {\it total variation} $P_n\stackrel{\rm TV}{\to}P$ means
\[\sup_{A\in \mathcal S}|P_n(A)-P(A)|\to0.\]
\end{definition}
\begin{theorem}\label{(3.10)} Scheffe's theorem. Suppose $P_n$ and $P$ have densities $f_n$ and $f$ with respect to a measure $\mu$ on $({\boldsymbol S},\mathcal S)$. If $f_n\to f$ almost everywhere  with respect to $\mu$, then $P_n\stackrel{\rm TV}{\to}P$ and therefore $P_n\Rightarrow P$.
\end{theorem}
Proof. For any $A\in \mathcal S$
\begin{align*}
 |P_n(A)-P(A)|&=\Big|\int_A(f_n(x)-f(x))\mu(dx)\Big|\le \int_{\boldsymbol S}|f(x)-f_n(x)|\mu(dx)\\
 &=2\int_{\boldsymbol S}(f(x)-f_n(x))^+\mu(dx),
\end{align*}
where the last equality follows from 
\begin{align*}
 0=\int_{\boldsymbol S} (f(x)-f_n(x))\mu(dx)=\int_{\boldsymbol S}(f(x)-f_n(x))^+\mu(dx)-\int_{\boldsymbol S}(f(x)-f_n(x))^-\mu(dx).
\end{align*}
On the other hand, by the dominated convergence theorem, $\int(f(x)-f_n(x))^+\mu(dx)\to0$.

\begin{example}\label{E3.3} According to Theorem \ref{(3.10)} the local limit theorem implies the integral limit theorem $P_n\Rightarrow P$. The reverse implication is false. Indeed, let $P=\mu$ be Lebesgue measure on ${\boldsymbol S}=[0,1]$ so that $f\equiv1$. Let $P_n$ be the uniform distribution on the set 
\[B_n=\bigcup_{k=0}^{n-1} (kn^{-1},kn^{-1}+n^{-3}) \]
with density $f_n(x)=n^{2}1_{\{x\in B_n\}}$. Since $\mu(B_n)=n^{-2}$, the Borel-Cantelli lemma implies that $\mu(B_n \mbox{ i.o.})=0$. Thus $f_n(x)\to0$ for almost all $x$ and there is no local theorem. On the other hand, $|P_n[0,x]-x|\le n^{-1}$ implying $P_n\Rightarrow P$.
\end{example}

\begin{theorem}\label{3.3} Let ${\boldsymbol S}=\boldsymbol R^k$. Denote by $L_n\subset \boldsymbol R^k$ a lattice with cells having dimensions $(\delta_1(n),\ldots,\delta_k(n))$ so that the cells of the lattice $L_n$ all having the form
\[B_n(x)=\{y:x_1-\delta_1(n)<y_1\le x_1,\ldots,x_k-\delta_k(n)<y_k\le x_k\},\quad x\in L_n\] 
have size $v_n=\delta_1(n)\cdots\delta_k(n)$. Suppose that  $(P_n)$ is a sequence of probability measures on $\boldsymbol R^k$, where $P_n$ is supported by $L_n$ with probability mass function $p_n(x)$. 

Suppose that  $P$ is a probability measure on $\boldsymbol R^k$ having density $f$ with respect to Lebesgue measure.  Assume that all $\delta_i(n)\to0$ as $n\to\infty$. If 
${p_n(x_n)\over v_n} \to f(x)$ whenever $x_n\in L_n$ and $x_n\to x$, then $P_n\Rightarrow P$.
\end{theorem}
Proof. Define a probability density $f_n$ on $\boldsymbol R^k$ by setting $f_n(y)={p_n(x)\over v_n}$ for $y\in B_n(x)$. It follows that $f_n(y)\to f(y)$ for all $y\in\boldsymbol R^k$. Let a random vector $Y_n$ have the density $f_n$ and $X$ have the density $f$. By Theorem \ref{(3.10)}, $Y_n\Rightarrow X$. Define $X_n$ on the same probability space as $Y_n$ by setting $X_n=x$ if $Y_n$ lies in the cell $B_n(x)$. Since $\|X_n-Y_n\|\le\|\delta(n)\|$, we conclude using Corollary \ref{3.1} that $X_n\Rightarrow X$.

\begin{example}\label{E3.4} If $S_n$ is the number of successes in $n$ Bernoulli trials, then according to the local form of the de Moivre-Laplace theorem,
\[\mathbb P(S_n=i) \sqrt{npq}={n\choose i}p^iq^{n-i}\sqrt{npq}\to{1\over\sqrt{2\pi}}e^{-z^2/2}\]
 provided $i$ varies with $n$ in such a way that ${i-np\over\sqrt{npq}}\to z$. Therefore, Theorem \ref{3.3} applies to the lattice 
 \[L_n=\Big\{{i-np\over\sqrt{npq}}, i\in\mathbb Z\Big\}\]
with $v_n={1\over\sqrt{npq}}$ and the probability mass function $p_n({i-np\over\sqrt{npq}})=\mathbb P(S_n=i)$ for $i=0,\ldots,n$. As a result we get the integral form of the de Moivre-Laplace theorem:
\[\mathbb P\Big({S_n-np\over \sqrt{npq}}\le z\Big)\to\Phi(z) \mbox{ as }n\to\infty  \mbox{ for any }z\in\boldsymbol R.\]
\end{example}

\section{Convergence of finite-dimensional distributions}
\subsection{Separating and convergence-determining classes}
\begin{definition}\label{p18}
 Call a subclass  $\mathcal A\subset \mathcal S$ a {\it separating class} if any two probability measures with $P(A)=Q(A)$ for all $A\in \mathcal A$, must be identical: $P(A)=Q(A)$ for all $A\in \mathcal S$.
 
Call a subclass  $\mathcal A\subset \mathcal S$ a {\it convergence-determining class} if, for every $P$ and every sequence $(P_n)$, convergence $P_n(A)\to P(A)$ for all $P$-continuity sets  $A\in\mathcal A$ implies $P_n\Rightarrow P$.
\end{definition}

\begin{lemma}\label{PM.42}
If $\mathcal A\subset\mathcal S$  is a  $\pi$-system and $\sigma(\mathcal A)=\mathcal S$, then $\mathcal A$ is a separating class.
\end{lemma}
Proof. Consider a pair of probability measures such that $P(A)=Q(A)$ for all $A\in \mathcal A$. Let $\mathcal L=\mathcal L_{P,Q}$ be the class of all sets $A\in \mathcal S$ such that $P(A)=Q(A)$.  Clearly, ${\boldsymbol S}\in \mathcal L$. If  $A\in \mathcal L$, then $A^c\in \mathcal L$ since $P(A^c)=1-P(A)=1-Q(A)=Q(A^c)$.
If  $A_n$ are disjoint sets in $\mathcal L$, then $\cup_nA_n\in \mathcal L$ since 
$$P(\cup_nA_n)=\sum_nP(A_n)=\sum_nQ(A_n)=Q(\cup_nA_n).$$
Therefore $\mathcal L$ is a $\lambda$-system, and since $\mathcal A\subset\mathcal L$, Theorem \ref{Dyn} gives $\sigma(\mathcal A)\subset\mathcal L $, and $\mathcal L=\mathcal S$. 

\begin{theorem}\label{2.3} Suppose that $P$ is a probability measure on a separable ${\boldsymbol S}$,  and a subclass $\mathcal A_P\subset \mathcal S$ satisfies

(i) $\mathcal A_P$ is a $\pi$-system,

(ii) for every $x\in {\boldsymbol S}$ and $\epsilon>0$, there is an $A\in\mathcal A_P$ for which $x\in A^\circ\subset A\subset B(x,\epsilon)$.

\noindent If $P_n(A)\to P(A)$ for every $A\in\mathcal A_P$, then $P_n\Rightarrow P$.
 \end{theorem}
Proof. If $A_1,\ldots,A_r$ lie in $\mathcal A_P$, so do their intersections. 
Hence, by the inclusion-exclusion formula and a theorem assumption, 
\begin{align*}
 P_n\Big(\bigcup_{i=1}^r A_i\Big)&=\sum_iP_n(A_i)-\sum_{ij}P_n(A_i\cap A_j)+\sum_{ijk}P_n(A_i\cap A_j\cap A_k)-\ldots\\
&\to\sum_iP(A_i)-\sum_{ij}P(A_i\cap A_j)+\sum_{ijk}P(A_i\cap A_j\cap A_k)-\ldots =P\Big(\bigcup_{i=1}^r A_i\Big).
\end{align*}
If  $G\subset {\boldsymbol S}$ is open,  then for each $x\in G$, $x\in A_x^\circ\subset A_x\subset G$ holds for some $A_x\in\mathcal A_P$. Since ${\boldsymbol S}$ is separable, by Theorem \ref{M3} (iii), there is a countable sub-collection $(A_{x_i}^\circ)$ that covers $G$. Thus  $G=\cup_i A_{x_i}$, where all $A_{x_i}$ are $\mathcal A_P$-sets.

With $A_i=A_{x_i}$ we have $G=\cup_i A_i$. Given $\epsilon$, choose $r$ so that $P\big(\cup_{i=1}^r A_i\big)>P(G)-\epsilon$. Then,
\[P(G)-\epsilon<P\Big(\bigcup_{i=1}^r A_i\Big)=\lim_nP_n\Big(\bigcup_{i=1}^r A_i\Big)\le\liminf_nP_n(G).\]
Now, letting $\epsilon\to0$ we find that for any open set $\liminf_nP_n(G)\ge P(G)$.

\begin{theorem}\label{2.4} Suppose that  ${\boldsymbol S}$ is separable and consider a subclass $\mathcal A\subset \mathcal S$. Let $\mathcal A_{x,\epsilon}$ be the class of $A\in\mathcal A$ satisfying $x\in A^\circ\subset A\subset B(x,\epsilon)$, and let $\partial\mathcal A_{x,\epsilon}$ be the class of their boundaries. If 

(i) $\mathcal A$ is a $\pi$-system,

(ii) for every $x\in {\boldsymbol S}$ and $\epsilon>0$,  $\partial\mathcal A_{x,\epsilon}$ contains uncountably many disjoint sets,

\noindent then  $\mathcal A$ is a convergence-determining class.
 \end{theorem}
Proof. For an arbitrary $P$ let  $\mathcal A_P$ be the class of $P$-continuity sets in $\mathcal A$. We have to show that if $P_n(A)\to P(A)$ holds for every $A\in\mathcal A_P$, then $P_n\Rightarrow P$. Indeed, by (i), since $\partial (A\cap B)\subset\partial (A)\cup\partial (B)$, $\mathcal A_P$  is a $\pi$-system. By (ii), there is an $A_x\in\mathcal A_{x,\epsilon}$ such that $P(\partial A_x)=0$ so that $A_x\in\mathcal A_P$. It remains to apply Theorem \ref{2.3}.

\subsection{Weak convergence in product spaces}

\begin{definition}\label{dpr}
 Let $P$ be a probability measure on ${\boldsymbol S}={\boldsymbol S}'\times {\boldsymbol S}''$ with the product metric 
 \[\rho((x',x''),(y',y''))=\rho'(x',y')\vee\rho''(x'',y'').\]
 Define the marginal distributions by $P'(A')=P(A'\times {\boldsymbol S}'')$ and $P''(A'')=P({\boldsymbol S}'\times A'')$. If the marginals are independent, we write $P=P'\times P''$. 
 We denote by $\mathcal S'\times \mathcal S''$ the product  $\sigma$-algebra  generated by the measurable rectangles $A'\times A''$ for $A'\in \mathcal S'$ and $A''\in \mathcal S''$. 
\end{definition}

\begin{lemma}\label{M10}
 If ${\boldsymbol S}={\boldsymbol S}'\times {\boldsymbol S}''$ is separable, then the three Borel $\sigma$-algebras are related by $\mathcal S=\mathcal S'\times \mathcal S''$.
\end{lemma}
Proof. Consider the projections $\pi':{\boldsymbol S}\to {\boldsymbol S}'$ and $\pi'':{\boldsymbol S}\to {\boldsymbol S}''$ defined by $\pi'(x',x'')=x'$ and $\pi''(x',x'')=x''$, each is continuous. 
For $A'\in \mathcal S'$ and $A''\in \mathcal S''$, we have
$$A'\times A''=(\pi')^{-1}A'\cap(\pi'')^{-1}A''\in \mathcal S,$$
since the two projections are continuous and therefore measurable. Thus $\mathcal S'\times \mathcal S''\subset \mathcal S$.
On the other hand, if ${\boldsymbol S}$ is separable, then each open set in ${\boldsymbol S}$ is a countable union of the balls
\[B((x',x''),r)=B'(x',r)\times B''(x'',r)\]
and hence lies in $\mathcal S'\times \mathcal S''$. Thus $\mathcal S\subset \mathcal S'\times \mathcal S''$.

\begin{theorem}\label{2.8} Consider probability measures $P_n$ and $P$ on a separable metric space ${\boldsymbol S}={\boldsymbol S}'\times {\boldsymbol S}''$.

(a) $P_n\Rightarrow P$ implies $P_n'\Rightarrow P'$ and $P_n''\Rightarrow P''$.
 
(b) $P_n\Rightarrow P$ if and only if $P_n(A'\times A'')\to P(A'\times A'')$ for each $P'$-continuity set $A'$ and each $P''$-continuity set $A''$.
 
(c) $P_n'\times P_n''\Rightarrow P$ if and only if $P_n'\Rightarrow P'$, $P_n''\Rightarrow P''$, and $P=P'\times P''$.
 \end{theorem}
 Proof. (a) Since $P'=P(\pi')^{-1}$, $P''=P(\pi'')^{-1}$ and the projections $\pi'$, $\pi''$ are continuous, it follows by the mapping theorem that $P_n\Rightarrow P$ implies $P_n'\Rightarrow P'$ and $P_n''\Rightarrow P''$.
 
 (b) Consider the $\pi$-system $\mathcal A$ of measurable rectangles $A'\times A''$: $A'\in \mathcal S'$ and $A''\in \mathcal S''$. Let $\mathcal A_P$ be the class of $A'\times A''\in \mathcal A$ such that 
 $P'(\partial A')=P''(\partial A'')=0$. Since 
 \[\partial (A'\cap B')\subset (\partial A')\cup (\partial B'),\qquad \partial (A''\cap B'')\subset (\partial A'')\cup (\partial B''),\]
 it follows that $\mathcal A_P$ is a  $\pi$-system:
 $$A'\times A'',B'\times B''\in \mathcal A_P\quad\Rightarrow \quad (A'\times A'')\cap (B'\times B'')\in \mathcal A_P.$$
 And since 
 $$\partial(A'\times A'')\subset((\partial A')\times {\boldsymbol S}'')\cup ({\boldsymbol S}'\times (\partial A'')),$$
 each set in $\mathcal A_P$ is a $P$-continuity set. Since $B'(x',r)$ in
 have disjoint boundaries for different values of $r$, and since the same is true of  the $B''(x'',r)$, there are arbitrarily small $r$ for which  $B(x,r)=B'(x',r)\times B''(x'',r)$ lies in $\mathcal A_P$. It follows that Theorem \ref{2.3} applies to $\mathcal A_P$: $P_n\Rightarrow P$ if and only if $P_n(A)\to P(A)$ for each $A\in\mathcal A_P$.

 The statement (c) is a consequence of (b).

\begin{exercise}\label{P2.7}The uniform distribution on the unit square and the uniform distribution on its diagonal have identical marginal distributions. Use this fact to demonstrate that the reverse to (a) in Theorem \ref{2.8} is false.
 
\end{exercise}
\begin{exercise}
Let $(X_n,Y_n)$ be a sequence of two-dimensional random vectors. 
 Show that if $(X_n,Y_n)\Rightarrow (X, Y)$, then besides $X_n\Rightarrow X$ and $Y_n\Rightarrow Y$, we have $X_n+Y_n\Rightarrow X+Y$.

 Give an example of $(X_n,Y_n)$ such that $X_n\Rightarrow X$ and $Y_n\Rightarrow Y$ but the sum $X_n+Y_n$ has no limit distribution.
\end{exercise}

\subsection{Weak convergence in $\boldsymbol R^k$ and $\boldsymbol R^\infty$ }\label{wcR}

Let $\boldsymbol R^k$ denote the $k$-dimensional Euclidean space with elements $x=(x_1,\ldots,x_k)$ and the ordinary metric 
\[\|x-y\|=\sqrt{(x_1-y_1)^2+\ldots+(x_k-y_k)^2}.\]
Denote by $\mathcal R^k$ the corresponding class of  $k$-dimensional Borel sets. Put $A_x=\{y:y_1\le x_1,\ldots,y_k\le x_k\}$, $x\in\boldsymbol R^k$. The probability measures on $(\boldsymbol R^k, \mathcal R^k)$ are completely determined by their distribution functions $F(x)=P(A_x)$ at the points of continuity $x\in\boldsymbol R^k$.

\begin{lemma}\label{Mtest} The Weierstrass M-test. Suppose that sequences of real numbers $x^n_i\to x_i$ converge for each i, and  for all $(n,i)$, $|x^n_i|\le M_i$, where $\sum_i M_i<\infty$. Then $\sum_i x_i<\infty$,  $\sum_i x^n_i<\infty$, and $\sum_i x^n_i\to\sum_i x_i$. 
\end{lemma}
Proof. The series of course converge absolutely, since $\sum_i M_i<\infty$. Now for any  $(n,i_0)$,
\[\Big|\sum_i x^n_i-\sum_i x_i\Big|\le \sum_{i\le i_0} |x^n_i-x_i|+2\sum_{i>i_0}M_i.\]
Given $\epsilon>0$, choose $i_0$ so that $\sum_{i>i_0}M_i<\epsilon/3$, and then choose $n_0$ so that $n>n_0$ implies $ |x^n_i-x_i|<{\epsilon\over 3i_0}$ for $i\le i_0$. Then  $n>n_0$ implies $|\sum_i x^n_i-\sum_i x_i|<\epsilon$.

\begin{lemma}\label{E1.2}
 Let $\boldsymbol R^\infty$ denote the space of the sequences $x=(x_1,x_2\ldots)$ of real numbers with metric
\[\rho(x,y)=\sum_{i=1}^\infty{1\wedge |x_i-y_i|\over2^i}.\]
Then $\rho(x^n,x)\to0$ if and only if $|x^n_i-x_i|\to0$ for each $i$.
\end{lemma}
Proof. If $\rho(x^n,x)\to0$, then  for each $i$ we have  $1\wedge |x^n_i-x_i|\to0$ and therefore  $|x^n_i-x_i|\to0$. The reverse implication holds by Lemma \ref{Mtest}.
\begin{definition}
 Let $\pi_k:\boldsymbol R^\infty\to\boldsymbol R^k$ be the natural projections $\pi_k(x)=(x_1,\ldots,x_k)$, $k=1,2,\ldots$, and let  $P$ be a probability measure on $(\boldsymbol R^\infty, \mathcal R^\infty)$. The probability measures $P\pi_k^{-1}$ defined on $(\boldsymbol R^k, \mathcal R^k)$ are called the {\it finite-dimensional distributions} of $P$.
\end{definition}
\begin{theorem}\label{p10}
 The space $\boldsymbol R^\infty$ is separable and complete. Let $P$ and $Q$ be two probability measures on $(\boldsymbol R^\infty, \mathcal R^\infty)$. 
 If $P\pi_k^{-1}=Q\pi_k^{-1}$  for each $k$, then $P=Q$.
\end{theorem}
Proof. Convergence in $\boldsymbol R^\infty$ implies coordinatewise convergence, therefore $\pi_k$ is continuous so that the sets
\[B_k(x,\epsilon)=\big\{ y\in\boldsymbol R^\infty:|y_i-x_i|<\epsilon,\ i=1,\ldots, k\big\}=\pi_k^{-1}\big\{ y\in\boldsymbol R^k:|y_i-x_i|<\epsilon,\ i=1,\ldots, k\big\}\]
are open. Moreover, $y\in B_k(x,\epsilon)$ implies $\rho(x,y)<\epsilon+2^{-k}$. Thus $B_k(x,\epsilon)\subset B(x,r)$ for $r>\epsilon+2^{-k}$. This means that the sets $B_k(x,\epsilon)$ form a base for the topology of $\boldsymbol R^\infty$. It follows that the space is separable: one countable, dense subset consists of those points having only finitely many nonzero coordinates, each of them rational.

If $(x^n)$ is a fundamental sequence, then each coordinate sequence $(x^n_i)$ is fundamental and hence converges to some $x_i$, implying $x^n\to x$. Therefore,  $\boldsymbol R^\infty$ is also complete.

Let $\mathcal A$ be  the class of finite-dimensional sets $\{x:\pi_k(x)\in H\}$ for some $k$ and some $H\in\mathcal R^k$. This class of cylinders is closed under finite intersections. To be able to apply Lemma \ref{PM.42} it remains to observe that  $\mathcal A$ generates $\mathcal R^\infty$: by separability each open set $G\subset\boldsymbol R^\infty$ is a countable union of sets in $\mathcal A$, since the sets $B_k(x,\epsilon)\in\mathcal A$ form a base.

\begin{theorem}\label{E2.4}
 Let $P_n,P$ be probability measures on $(\boldsymbol R^\infty, \mathcal R^\infty)$. Then $P_n\Rightarrow P$ if and only if $P_n\pi_k^{-1}\Rightarrow P\pi_k^{-1}$ for each $k$.
\end{theorem}
Proof. Necessity follows from the mapping theorem. Turning to sufficiency, let $\mathcal A$, again,  be  the class of finite-dimensional sets $\{x:\pi_k(x)\in H\}$ for some $k$ and some $H\in\mathcal R^k$. We proceed in three steps. 

Step 1. Show that $\mathcal A$ is a convergence-determining class.  This is proven using Theorem \ref{2.4}. Given $x$ and $\epsilon$, choose $k$ so that $2^{-k}<\epsilon/2$ and consider the collection of uncountably many finite-dimensional sets 
\[A_\eta=\{y:|y_i-x_i|<\eta, i=1,\ldots, k\}\mbox{ for }0<\eta<\epsilon/2.\]
We have $A_\eta\in\mathcal A_{x,\epsilon}$. On the other hand, $\partial A_\eta$ consists of the points $y$ such that $|y_i-x_i|\le\eta$ with equality for some $i$, hence these boundaries are disjoint. And since $\boldsymbol R^\infty$ is separable, Theorem \ref{2.4} applies.

Step 2. Show that $\partial (\pi_k^{-1}H)=\pi_k^{-1}\partial H$.

From the continuity of $\pi_k$ it follows that $\partial (\pi_k^{-1}H)\subset\pi_k^{-1}\partial H$. Using special properties of the projections we can prove inclusion in the other direction. If 
$x\in\pi_k^{-1}\partial H$, so that $\pi_kx\in\partial H$, then there are points $\alpha^{(u)}\in H$, $\beta^{(u)}\in H^c$ such that $\alpha^{(u)}\to \pi_kx$ and $\beta^{(u)}\to \pi_kx$ as $u\to\infty$. Since the points $(\alpha^{(u)}_1,\ldots,\alpha^{(u)}_k,x_{k+1},\ldots)$ lie in $\pi_k^{-1} H$ and converge to $x$, and since the points $(\beta^{(u)}_1,\ldots,\beta^{(u)}_k,x_{k+1},\ldots)$ lie in $(\pi_k^{-1} H)^c$ and converge to $x$, we conclude that $x\in \partial (\pi_k^{-1}H)$.

Step 3. Suppose that $P\pi_k^{-1}(\partial H)=0$ implies $P_n\pi_k^{-1}(H)\to P\pi_k^{-1}(H)$ and show that $P_n\Rightarrow P$.

If $A\in \mathcal A$ is a finite-dimensional $P$-continuity set, then we have  $A=\pi_k^{-1}H$ and
\[P\pi_k^{-1}(\partial H)=P(\pi_k^{-1}\partial H)=P(\partial \pi_k^{-1}H)=P(\partial A)=0.\]
Thus by assumption, $P_n(A)\to P(A)$ and according to step 1, $P_n\Rightarrow P$.

\subsection{Kolmogorov's extension theorem}

\begin{definition}\label{Kcon}
 We say that the system of finite-dimensional distributions $\mu_{t_1,\ldots,t_k}$ is consistent if the joint distribution functions
 \[F_{t_1,\ldots,t_k}(z_1,\ldots,z_k)=\mu_{t_1,\ldots,t_k}((-\infty,z_1]\times \ldots\times(-\infty,z_k])\]
 satisfy two consistency conditions
 
 (i) $F_{t_1,\ldots,t_k,t_{k+1}}(z_1,\ldots,z_k,\infty)= F_{t_1,\ldots,t_k}(z_1,\ldots,z_k)$,
 
 (ii) if $\pi$ is a permutation of $(1,\ldots,k)$, then $$F_{t_{\pi(1)},\ldots,t_{\pi(k)}}(z_{\pi(1)},\ldots,z_{\pi(k)})=F_{t_1,\ldots,t_k}(z_1,\ldots,z_k).$$

\end{definition}
\begin{theorem}\label{ket} 
Let  $\mu_{t_1,\ldots,t_k}$ be a consistent system of finite-dimensional distributions. Put $\Omega=\{\mbox{functions }\omega:[0,1]\to \mathbb R\}$ and $\mathcal F$ is the $\sigma$-algebra generated by the finite-dimensional sets $\{\omega: \omega(t_i)\in B_i, i=1,\ldots,n\}$, where $B_i$ are Borel subsets of $\mathbb R$. Then there is a unique probability measure $\mathbb P$ on $(\Omega,\mathcal F)$ such that a stochastic process defined by $X_t(\omega)=\omega(t)$ has the finite-dimensional distributions $\mu_{t_1,\ldots,t_k}$.
\end{theorem}
Without proof.  Kolmogorov's extension theorem does not directly imply the existence of the Wiener process because the $\sigma$-algebra $\mathcal F$ is not rich enough to ensure the continuity property for trajectories. However, it is used in the proof of Theorem \ref{13.6} establishing the existence of processes with cadlag trajectories.

\section{Tightness and Prokhorov's theorem}\label{secP}

\subsection{Tightness of probability measures}
Convergence of finite-dimensional distributions does not always imply weak convergence. This makes important the following concept of tightness. 

\begin{definition}
 A family of probability measures $\Pi$ on $({\boldsymbol S},\mathcal S)$ is called {\it tight} if for every $\epsilon$ there exists a compact set $K\subset {\boldsymbol S}$ such that $P(K)>1-\epsilon$ for all $P\in\Pi$. 
\end{definition}
\begin{lemma}\label{1.3}
 If ${\boldsymbol S}$ is separable and complete, then each probability measure $P$ on $({\boldsymbol S},\mathcal S)$ is tight. 
\end{lemma}
Proof. Separability: for each $k$ there is a sequence $A_{k,i}$ of open $1/k$-balls covering ${\boldsymbol S}$. Choose $n_k$ large enough that  $P(B_k)>1-\epsilon 2^{-k}$ where $B_k=A_{k,1}\cup\ldots\cup A_{k,n_k}$. Completeness: the totally bounded set $B_{1}\cap B_2\cap\ldots$ has compact closure $K$. But clearly $P(K^c)\le\sum_k P(B_k^c)<\epsilon$.
\begin{exercise}
 Check whether the following sequence of distributions on $\boldsymbol R$ 
 \[
P_n(A)=(1-n^{-1})1_{\{0\in A\}}+n^{-1}1_{\{n^2\in A\}},\qquad n\ge1,
\]
is tight or it ``leaks" towards infinity. Notice that the corresponding mean value is $n$.
\end{exercise}
\begin{definition}
 A family of probability measures $\Pi$ on $({\boldsymbol S},\mathcal S)$ is called {\it relatively compact} if any sequence of its elements contains a weakly convergent subsequence. The limiting probability measures might be different for different subsequences and lie outside $\Pi$.
\end{definition}

\begin{definition}\label{p72}
 Let $\boldsymbol P$ be the space of probability measures on $({\boldsymbol S},\mathcal S)$. The {\it Prokhorov distance} $\pi(P,Q)$ between $P,Q\in \boldsymbol P$ is defined as the infimum of those positive $\epsilon$ for which 
 \[P(A)\le Q(A^\epsilon)+\epsilon, \quad Q(A)\le P(A^\epsilon)+\epsilon, \quad\mbox{for all }A\in\mathcal S.\]
\end{definition}
\begin{lemma}
 The Prokhorov distance $\pi$ is a metric on $\boldsymbol P$.
\end{lemma}
Proof. Obviously $\pi(P,Q)=\pi(Q,P)$ and $\pi(P,P)=0$. If  $\pi(P,Q)=0$, then for any $F\in\mathcal S$ and  $\epsilon>0$, $P(F)\le Q(F^\epsilon)+\epsilon$. For closed $F$ letting $\epsilon\to0$ gives $P(F)\le Q(F)$. By symmetry, we have $P(F)= Q(F)$ implying $P= Q$.

To verify the triangle inequality notice that if $\pi(P,Q)<\epsilon_1$ and $\pi(Q,R)<\epsilon_2$, then
 \[P(A)\le Q(A^{\epsilon_1})+\epsilon_1\le R((A^{\epsilon_1})^{\epsilon_2})+\epsilon_1+\epsilon_2\le R(A^{\epsilon_1+\epsilon_2})+\epsilon_1+\epsilon_2.\]
Thus, using the symmetric relation we obtain $\pi(P,R)<\epsilon_1+\epsilon_2$. Therefore, $\pi(P,R)\le \pi(P,Q)+\pi(Q,R)$.
\begin{theorem}\label{6.8}
Suppose ${\boldsymbol S}$ is a complete separable metric space.  Then weak convergence is equivalent to $\pi$-convergence, $(\boldsymbol P,\pi)$ is separable and complete, and $\Pi\subset \boldsymbol P$ is relatively compact iff its $\pi$-closure is $\pi$-compact.
 \end{theorem}
Without proof.

\begin{theorem}\label{2.6}
A necessary and sufficient condition for $P_n\Rightarrow P$ is that each subsequence $P_{n'}$ contains a further subsequence $P_{n''}$ converging weakly to $P$.
 \end{theorem}
 Proof. The necessity is easy but useless. As for sufficiency, if $P_n\nRightarrow P$, then $\int_{{\boldsymbol S}}f(x)P_n(dx)\nrightarrow\int_{{\boldsymbol S}}f(x)P(dx)$ for some bounded, continuous $f$. But then, for some $\epsilon>0$ and some subsequence $P_{n'}$,
 \[\Big|\int_{{\boldsymbol S}}f(x)P_{n'}(dx)-\int_{{\boldsymbol S}}f(x)P(dx)\Big|\ge\epsilon\quad \mbox{ for all }n',\]
 and no further subsequence can converge weakly to $P$.
 
\begin{theorem}\label{5.1}
Prokhorov's theorem, the direct part. If a family of probability measures $\Pi$ on $({\boldsymbol S},\mathcal S)$ is tight, then it is relatively compact. 
\end{theorem}
Proof. See the next subsection.

\begin{theorem}\label{5.2}
Prokhorov's theorem, the reverse part. Suppose ${\boldsymbol S}$ is a complete separable metric space. If $\Pi$ is relatively compact, then it is tight.
\end{theorem}
Proof. Consider open sets $G_n\uparrow {\boldsymbol S}$. For each $\epsilon$ there is an $n$ such that $P(G_n)>1-\epsilon$ for all $P\in\Pi$. To show this we assume the opposite: $P_n(G_n)\le1-\epsilon$ for some $P_n\in\Pi$. By the assumed relative compactness, $P_{n'}\Rightarrow Q$ for some subsequence and some probability measure $Q$. Then $$Q(G_n)\le\liminf_{n'}P_{n'}(G_n)\le\liminf_{n'}P_{n'}(G_{n'})\le1-\epsilon$$ 
which is impossible since $G_n\uparrow {\boldsymbol S}$.

If $A_{ki}$ is a sequence of open balls of radius $1/k$ covering ${\boldsymbol S}$ (separability), so that ${\boldsymbol S}=\cup_iA_{k,i}$ for each $k$. From the previous step, it follows that there is an $n_k$ such that  $P(\cup_{i\le n_k}A_{k,i})>1-\epsilon2^{-k}$ for all $P\in\Pi$. Let $K$ be the closure of the totally bounded set $\cap_{k\ge1}\cup_{i\le n_k}A_{k,i}$, then $K$ is compact (completeness) and $P(K)>1-\epsilon$ for all $P\in\Pi$.

%\begin{theorem}\label{p58}
% Let $P$ and $P_n$ be probability measures on $(S,\mathcal S)$. Suppose $\mathcal A$ is a separating class and $P_n(A)\to P(A)$ for all $A\in \mathcal A$ such that $P(\partial A)=0$. If, moreover, $P_n$ is tight, then  $P_n\Rightarrow P$.
%
%\end{theorem}
%Proof. Tightness implies relative compactness which in turn implies that each subsequence $(P_{n'})\subset (P_n)$ contains a further  subsequence $(P_{n''})\subset (P_{n'})$ converging weekly to some $Q$. From $P_{n}(A)\to P(A)$ and $P_{n''}(A)\to Q(A)$ we get $P(A)=Q(A)$ for all $A\in \mathcal A$ such that $P(\partial A)=Q(\partial A)=0$. Since $\mathcal A$ is a separating class, we conclude $Q=P$. It remains to apply Theorem \ref{2.6}.

\subsection{Proof of Prokhorov's theorem}

This subsection contains a proof of the direct half of Prokhorov's theorem. Let $(P_n)$ be a sequence in the tight family $\Pi$. We are to find a subsequence $(P_{n'})$ and a probability measure $P$ such that $P_{n'}\Rightarrow P$. The proof, like that of Helly's selection theorem will depend on a diagonal argument.

Choose compact sets $K_1\subset K_2\subset\ldots$ such that $P_n(K_i)>1-i^{-1}$ for all $n$ and $i$. The set $K_\infty=\cup_iK_i$ is separable: compactness = each open cover has a finite subcover, separability = each open cover has a countable subcover. Hence, by Theorem \ref{M3'}, there exists a countable class $\mathcal A$ of open sets with the following property: if $G$ is open and $x\in K_\infty\cap G$, then $x\in A\subset A^-\subset G$ for some $A\in\mathcal A$. Let $\mathcal H$ consist of $\emptyset$ and the finite unions of sets of the form $ A^-\cap K_i$ for $A\in\mathcal A$ and $i\ge1$.

Consider the countable class $\mathcal H=(H_j)$. For $(P_n)$ there is  a subsequence $(P_{n_1})$ such that  $P_{n_1}(H_1)$ converges as $n_1\to\infty$. For $(P_{n_1})$  there is  a further subsequence $(P_{n_2})$ such that  $P_{n_2}(H_2)$ converges as $n_2\to\infty$. Continuing in this way we get a collection of indices $(n_{1k})\supset(n_{2k})\supset\ldots$ such that $P_{n_{jk}}(H_j)$ converges as $k\to\infty$ for each $j\ge1$. Putting $n'_j=n_{jj}$ we find a subsequence $(P_{n'})$ for which the limit 
\[\alpha(H)=\lim_{n'}P_{n'}(H)\mbox{ exists  for each }H\in\mathcal H.\]
Furthermore, for open sets $G\subset {\boldsymbol S}$ and arbitrary sets $M\subset {\boldsymbol S}$ define
\[\beta(G)=\sup_{H\subset G}\alpha(H),\quad \gamma(M)=\inf_{G\supset M}\beta(G).\]
Our objective is to construct on $({\boldsymbol S},\mathcal S)$ a probability measure $P$ such that $P(G)=\beta(G)$ for all open sets  $G$.
If there does exist such a $P$, then the proof will be complete: if $H\subset G$, then 
\[\alpha(H)=\lim_{n'}P_{n'}(H)\le \liminf_{n'}P_{n'}(G),\]
whence $P(G)\le\liminf_{n'}P_{n'}(G)$, and therefore $P_{n'}\Rightarrow P$. 
The construction of the probability measure $P$ is divided in seven steps.

 Step 1: if $F\subset G$, where $F$ is closed and $G$ is open, and if $F\subset H$, for some $H\in\mathcal H$, then $F\subset H_0\subset G$, for some $H_0\in\mathcal H$.

Since $F\subset K_{i_0}$ for some $i_0$, the closed set $F$ is compact. For each $x\in F$, choose an $A_x\in\mathcal A$ such that $x\in A_x\subset A_x^-\subset G$. The sets $A_x$ cover the compact $F$, and there is a finite subcover $A_{x_1},\ldots,A_{x_k}$. We can take $H_0=\cup_{j=1}^k(A_{x_j}^-\cap K_{i_0})$.

Step 2: $\beta$ is finitely subadditive on the open sets.

Suppose that $H\subset G_1\cup G_2$, where $H\in\mathcal H$ and $G_1, G_2$ are open. Define
\begin{align*}
 F_1&=\big\{x\in H:\rho(x,G_1^c)\ge \rho(x,G_2^c)\big\},\\
F_2&=\big\{x\in H:\rho(x,G_2^c)\ge \rho(x,G_1^c)\big\},
\end{align*}
so that $H=F_1\cup F_2$ with $F_1\subset G_1$ and $F_2\subset G_2$. According to Step 1, since $F_i\subset H$, we have $F_i\subset H_i\subset G_i$ for some $H_i\in\mathcal H$.

The function $\alpha(H)$ has these three properties 
\begin{align*}
 \alpha(H_1)&\le\alpha(H_2)\qquad\qquad\quad \mbox{ if }H_1\subset H_2,\\
 \alpha(H_1\cup H_2)&= \alpha(H_1)+\alpha(H_2)\quad \mbox{ if }H_1\cap H_2=\emptyset,\\
 \alpha(H_1\cup H_2)&\le \alpha(H_1)+\alpha(H_2).
\end{align*}
It follows first,
\[\alpha(H)\le \alpha(H_1\cup H_2)\le\alpha(H_1)+\alpha(H_2)\le \beta(G_1)+\beta(G_2),\]
and then
$$\beta(G_1\cup G_2)=\sup_{H\subset G_1\cup G_2}\alpha(H)\le \beta(G_1)+\beta(G_2).$$

Step 3: $\beta$ is countably subadditive on the open sets.

If $H\subset \cup_n G_n$, then, since $H$ is compact, $H\subset \cup_{n\le n_0} G_n$ for some $n_0$, and finite subadditivity imples
\[\alpha(H)\le \sum_{n\le n_0}\beta(G_n)\le \sum_{n}\beta(G_n).\]
Taking the supremum over $H$ contained in $\cup_{n} G_n$ gives $\beta(\cup_{n}G_n)\le \sum_{n}\beta(G_n)$.

Step 4: $\gamma$ is an outer measure.

Since $\gamma$ is clearly monotone and satisfies $\gamma(\emptyset)=0$, we need only prove that it is countably subadditive. Given a positive $\epsilon$ and arbitrary $M_n\subset {\boldsymbol S}$, choose open sets $G_n$ such that  $M_n\subset G_n$ and $\beta(G_n)<\gamma(M_n)+\epsilon/2^n$. Apply Step 3
$$\gamma(\bigcup_{n}M_n)\le \beta(\bigcup_{n}G_n)\le \sum_{n}\beta(G_n)\le \sum_{n}\gamma(M_n)+\epsilon,$$
 and let $\epsilon\to0$ to get
 $\gamma(\bigcup_{n}M_n)\le \sum_{n}\gamma(M_n)$.

Step 5: $\beta(G)\ge\gamma(F\cap G)+\gamma(F^c\cap G)$ for $F$ closed and $G$ open.

Choose $H_3,H_4\in\mathcal H$ for which \begin{align*}
&H_3\subset F^c\cap G\quad \mbox{ and }\quad \alpha(H_3)>\beta(F^c\cap G)-\epsilon,\\
&H_4\subset H_3^c\cap G\quad \mbox{ and }\quad \alpha(H_4)>\beta(H_3^c\cap G)-\epsilon.
\end{align*}
Since $H_3$ and $H_4$ are disjoint and are contained in $G$, it follows from the properties of the functions $\alpha,\beta,$ and $\gamma$ that
\begin{align*}
\beta(G)\ge \alpha(H_3\cup H_4)=\alpha(H_3)+\alpha(H_4)&>\beta(F^c\cap G)+\beta(H_3^c\cap G)-2\epsilon\\
&\ge \gamma(F^c\cap G)+ \gamma(F\cap G)-2\epsilon.
\end{align*}
Now it remains to  let $\epsilon\to0$.

Step 6: if $F\subset {\boldsymbol S}$ is closed, then $F$ is in the class $\mathcal M$ of $\gamma$-measurable sets.

By Step 5, $\beta(G)\ge\gamma(F\cap L)+\gamma(F^c\cap L)$ if $F$ is closed, $G$ is open, and $G\supset L$. Taking the infimum over these $G$ gives $\gamma(L)\ge\gamma(F\cap L)+\gamma(F^c\cap L)$ confirming that $F$ is $\gamma$-measurable.

Step 7: $\mathcal S\subset\mathcal M$, and the restriction $P$ of $\gamma$ to  $\mathcal S$ is a probability measure satisfying $P(G)=\gamma(G)=\beta(G)$ for all open sets $G\subset {\boldsymbol S}$.

Since each closed set lies in $\mathcal M$ and $\mathcal M$ is a $\sigma$-algebra, we have $\mathcal S\subset\mathcal M$.
To see that the  $P$ is a probability measure, observe that each $K_i$ has a finite covering by $\mathcal A$-sets and therefore $K_i\in\mathcal H$. Thus
\[1\ge P({\boldsymbol S})=\beta({\boldsymbol S})\ge \sup_i\alpha(K_i)\ge \sup_i(1-i^{-1})=1.\]

\subsection{Skorokhod's representation theorem}
\begin{theorem}\label{6.7}
 Suppose that $P_n\Rightarrow P$ and $P$ has a separable support. Then there exist random elements $X_n$ and $X$, defined on a common probability space $(\Omega,\mathcal F,\mathbb P)$, such that 
 $P_n$ is the probability distribution of $X_n$, $P$ is the probability distribution of $X$, and $X_n(\omega)\to X(\omega)$ for every $\omega$.

\end{theorem}
Proof. We split the proof in four steps.

Step 1: show that for each $\epsilon$, there is a finite $\mathcal S$-partition $B_0,B_1,\ldots, B_k$ of ${\boldsymbol S}$ such that 
\[0<P(B_0)<\epsilon,\quad P(\partial B_i)=0,\quad {\rm diam} (B_i)<\epsilon,\quad i=1,\ldots,k.\]

Let $M$ be a separable $\mathcal S$-set for which $P(M)=1$. For each $x\in M$, choose $r_x$ so that $0<r_x<\epsilon/2$ and $P(\partial B(x,r_x))=0$. Since 
$M$ is a separable, it can be covered by a countable subcollection $A_1,A_2,\ldots$ of the balls $B(x,r_x)$. Choose $k$ so that $P(\cup_{i=1}^k A_i)>1-\epsilon$. Take
\[B_0=\big(\bigcup_{i=1}^k A_i\big)^c,\quad B_1=A_1,\quad B_i=A_1^c\cap\ldots\cap A_{i-1}^c\cap A_i,\]
and notice that $\partial B_i\subset \partial A_1\cup\ldots\cup\partial A_k$.

Step 2: definition of $n_j$. 

Take $\epsilon_j=2^{-j}$. By step 1, there are $\mathcal S$-partitions $B_0^j,B_1^j,\ldots, B_k^j$ such that 
\[0<P(B_0^j)<\epsilon_j,\quad P(\partial B_i^j)=0,\quad {\rm diam} (B_i^j)<\epsilon_j,\quad i=1,\ldots,k_j.\]
If some $P(B_i^j)=0$, we redefine these partitions by amalgamating such $B_i^j$ with $B_0^j$, so that $P(\cdot|B_i^j)$ is well defined for $i\ge1$. By the assumption $P_n\Rightarrow P$, there is for each $j$ an $n_j$ such that 
\[P_n(B_i^j)\ge (1-\epsilon_j)P(B_i^j),\quad i=0,1,\ldots, k_j,\quad n\ge n_j.\]
Putting $n_0=1$, we can assume $n_0< n_1<\cdots$.

Step 3: construction of  $X, Y_n,Y_{ni},Z_n,\xi$. 

Define $m_n=j$ for $n_j\le n<n_{j+1}$ and write $m$ instead of $m_n$. By Theorem \ref{ket} we can find an $(\Omega,\mathcal F,\mathbb P)$ supporting random elements $X, Y_n,Y_{ni},Z_n$ of ${\boldsymbol S}$ and a random variable $\xi$, all independent of each other and having distributions satisfying: $X$ has distribution $P$, $Y_n$ has distribution $P_n$, 
\begin{align*}
 &\mathbb P(Y_{ni}\in A)=P_n(A|B_i^m),  \quad  \mathbb P(\xi\le \epsilon)=\epsilon,\\
 &\epsilon_m\mathbb P(Z_n\in A)=\sum_{i=0}^{k_m}P_n(A|B^m_i)\Big(P_n(B^m_i)-(1-\epsilon_m)P(B^m_i)\Big).
\end{align*}
Note that $\mathbb P(Y_{ni}\in B_i^m)=1$.

Step 4: construction of $X_n$. 

Put $X_n=Y_n$ for $n<n_1$. For $n\ge n_{1}$, put 
\[X_n=1_{\{\xi\le1-\epsilon_m\}}\sum_{i=0}^{k_m}1_{\{X\in B_i^m\}}Y_{ni}+1_{\{\xi>1-\epsilon_m\}}Z_n.\]
By step 3, we $X_n$ has distribution $P_n$ because
\begin{align*}
 \mathbb P(X_n\in A)&=(1-\epsilon_m)\sum_{i=0}^{k_m}\mathbb P(X\in B_i^m,Y_{ni}\in A)+\epsilon_m\mathbb P(Z_n\in A)\\
 &=(1-\epsilon_m)\sum_{i=0}^{k_m}\mathbb P(X\in B_i^m)P_n(A|B_i^m)\\
 &\qquad+\sum_{i=0}^{k_m}P_n(A|B^m_i)\Big(P_n(B^m_i)-(1-\epsilon_m)P(B^m_i)\Big)\\
 &=P_n(A).
\end{align*}
Let 
$$E_j=\{X\notin B_0^j;\ \xi\le1-\epsilon_j \}\mbox{ and }E=\liminf_jE_j=\bigcup_{j=1}^\infty\bigcap_{i=j}^\infty E_i.$$ 
Since $\mathbb P(E_j^c)<2\epsilon_j$, by the Borel-Cantelli lemma, $\mathbb P(E^c)=\mathbb P(E_j^c \mbox{ i.o.})=0$ implying $\mathbb P(E)=1$. If $\omega\in E$, then both $X_n(\omega)$ and $X(\omega)$ lie in the same $B_i^m$ having diameter less than $\epsilon_m$. Thus, $\rho(X_n(\omega),X(\omega))<\epsilon_m$ and $X_n(\omega)\to X(\omega)$ for $\omega\in E$. It remains to redefine $X_n$ as $X$ outside $E$.

\begin{corollary} The mapping theorem. Let $h:{\boldsymbol S}\to {\boldsymbol S}'$ be a continuous mapping  between two metric spaces. If $P_n\Rightarrow P$ on ${\boldsymbol S}$ and $P$ has a separable support, then $P_nh^{-1}\Rightarrow Ph^{-1}$ on ${\boldsymbol S}'$.
 
\end{corollary}
Proof. Having $X_n(\omega)\to X(\omega)$ we get $h(X_n(\omega))\to h(X(\omega))$ for every $\omega$. It follows, by Corollary \ref{3.1} that $h(X_n)\Rightarrow h(X)$ which is equivalent to $P_nh^{-1}\Rightarrow Ph^{-1}$.

\section{Functional Central Limit Theorem on $\boldsymbol C=\boldsymbol C[0,1]$}\label{secC}

\subsection{Weak convergence in $\boldsymbol C$}
\begin{definition}
  An element of the set $\boldsymbol C=\boldsymbol C[0,1]$ is a continuous function $x=x(t)$. The distance between points in $\boldsymbol C$ is measured by the uniform metric
\[\rho(x,y)=\|x-y\|=\sup_{0\le t\le1}|x(t)-y(t)|.\]
Denote by $\mathcal C$ the Borel $\sigma$-algebra of subsets of $\boldsymbol C$.
\end{definition}
\begin{exercise}
Draw a picture for an open ball $B(x,r)$ in $\boldsymbol C$.\\
For any real number $a$ and $t\in[0,1]$ the set $\{x:x(t)<a\}$ is an open subset of $\boldsymbol C$.
\end{exercise}
\begin{example}\label{E1.3}
Convergence $\rho(x_n,x)\to0$ means uniform convergence of continuous functions, it is stronger than pointwise convergence. Consider the function  $z_n(t)$ that increases linearly from 0 to 1 over $[0,n^{-1}]$, decreases linearly from 1 to 0 over $[n^{-1},2n^{-1}]$, and equals 0 over $[2n^{-1},1]$. Despite $z_n(t)\to0$ for any $t$ we have $\|z_n\|=1$ for all $n$.
\end{example}
\begin{theorem}\label{p11} 
The space $\boldsymbol C$ is separable and complete.
 \end{theorem}
 Proof. Separability. Let $L_k$ be the set of polygonal functions that are linear over each subinterval $[{i-1\over k},{i\over k}]$ and have rational values at the end points. We will show that the countable set $\cup_{k\ge 1}L_k$ is dense in $\boldsymbol C$. For given $x\in \boldsymbol C$ and $\epsilon>0$, choose $k$ so that 
 \[|x(t)-x(i/k)|<\epsilon\quad\mbox{for all }t\in[{(i-1)/ k},{i/k}],\quad 1\le i\le k\]
which is possible by uniform continuity. Then choose $y\in L_k$ so that $|y(i/k)-x(i/k)|<\epsilon$ for each $i$. It remains to draw a picture with trajectories over an interval $[{i-1\over k},{i\over k}]$ and check that $\rho(x,y)\le3\epsilon$.

Completeness. Let $(x_n)$ be a fundamental sequence so that
\[\epsilon_n=\sup_{m> n}\sup_{0\le t\le1}|x_n(t)-x_m(t)|\to0,\quad n\to\infty.\]
Then for each $t$, the sequence $(x_n(t))$ is fundamental on $\boldsymbol R$ and hence has a limit $x(t)$. Letting $m\to\infty$ in the inequality $|x_n(t)-x_m(t)|\le \epsilon_n$ gives $|x_n(t)-x(t)|\le \epsilon_n$. Thus $x_n$ converges uniformly to $x\in \boldsymbol C$.

\begin{definition}
 Convergence of finite-dimensional distributions $X^n\stackrel{\rm fdd}{\longrightarrow}X$ means that for all $t_1,\ldots,t_k$
$$(X^n_{t_1},\ldots,X^n_{t_k})\Rightarrow (X_{t_1},\ldots,X_{t_k}).$$
\end{definition}
\begin{exercise}
The projection $\pi_{t_1,\ldots,t_k}:\boldsymbol C\to \boldsymbol R^k$ defined by $\pi_{t_1,\ldots,t_k}(x)=(x(t_1),\ldots,x(t_k))$ is a continuous map.
\end{exercise}
\begin{example} By the mapping theorem, if $X^n\Rightarrow X$, then $X^n\stackrel{\rm fdd}{\longrightarrow}X$. The reverse in not true.
 Consider $z_n(t)$ from Example \ref{E1.3} and put $X^n=z_n$, $X=0$ so that  $X^n\stackrel{\rm fdd}{\longrightarrow}X$. 
Take $h(x)=\sup_tx(t)$. It satisfies $|h(x)-h(y)|\le\rho(x,y)$ and therefore is a continuous function on $\boldsymbol C$. 
Since $h(z_n)\equiv1$, we have  $h(X^n)\nRightarrow h(X)$, and according to the mapping theorem $X^n\nRightarrow X$. 
\end{example}
\begin{definition}
 Define a modulus of continuity of a function $x:[0,1]\to \boldsymbol R$ by
\[w_x(\delta)=w(x,\delta)=\sup_{|s-t|\le\delta}|x(s)-x(t)|,\quad \delta\in(0,1].\]
\end{definition}
For any $x:[0,1]\to \boldsymbol R$ its modulus of continuity $w_x(\delta)$ is non-decreasing over $\delta$. Clearly, $x\in \boldsymbol C$ if and only if $w_x(\delta)\to0$ as $\delta\to0$. The limit $j_x=\lim_{\delta\to0}w_x(\delta)$ is the absolute value of the largest jump of $x$.
\begin{exercise}
Show that for any fixed $\delta\in(0,1]$ we have $|w_x(\delta)-w_y(\delta)|\le2\rho(x,y)$ implying that $w_x(\delta)$ is a continuous function on $\boldsymbol C$.
\end{exercise}
\begin{example}
 For $z_n\in \boldsymbol C$ defined in Example \ref{E1.3} we have $w(z_n,\delta)=1$ for $n\ge \delta^{-1}$.
\end{example}

\begin{exercise}
Given a probability measure $P$ on the measurable space $(\boldsymbol C,\mathcal C)$ there exists a random process $X$ on a probability space $(\Omega,\mathcal F,\mathbb P)$ such that $\mathbb P(X\in A)=P(A)$ for any $A\in\mathcal C$.
\end{exercise}

\begin{theorem} \label{7.5}
Let $P_n, P$ be probability measures on $(\boldsymbol C,\mathcal C)$. Suppose $P_{n}\pi^{-1}_{t_1,\ldots,t_k}\Rightarrow P\pi^{-1}_{t_1,\ldots,t_k}$ holds for all tuples  $(t_1,\ldots,t_k)\subset[0,1]$. If for every positive $\epsilon$ 
\[
(i) \quad \quad \quad \quad \lim_{\delta\to0}\limsup_{n\to\infty}P_n(x: w_x(\delta)\ge\epsilon)=0,\quad \quad \quad \quad %\label{sCp} 
\]
then $P_n\Rightarrow P$.
\end{theorem}

Proof. The proof is given in terms of convergence in distribution using Theorem \ref{3.2}.

For $u=1,2,\ldots$, define $M_u:\boldsymbol C\to \boldsymbol C$ in the following way. Let $(M_ux)(t)$ agree with $x(t)$ at the points $0,1/u,2/u,\ldots,1$ and be defined by linear interpolation between these points. Observe that $\rho(M_ux,x)\le2w_x(1/u)$.

Further, for a vector $\alpha=(\alpha_0,\alpha_1,\ldots,\alpha_u)$ define $(L_u\alpha)(t)$ as an element of $\boldsymbol C$ such that it has values $\alpha_i$ at points $t=i/n$ and is linear in between. Clearly, $\rho(L_u\alpha,L_u\beta)=\max_i|\alpha_i-\beta_i|$, so that $L_u:\boldsymbol R^{u+1}\to \boldsymbol C$ is continuous.

Let $t_i=i/u$. Observe that $M_u=L_u\pi_{t_0,\ldots,t_u}$. Since $\pi_{t_0,\ldots,t_u}X^n\Rightarrow\pi_{t_0,\ldots,t_u}X$ and $L_u$ is continuous, the mapping theorem gives $M_uX^n\Rightarrow M_uX$ as $n\to\infty$. 
Since
\[\limsup_{u\to\infty}\rho(M_uX,X)\le2\limsup_{u\to\infty}w(X,1/u)=0,\]
we have $M_uX\to X$ in probability and therefore $M_uX\Rightarrow X$. 

Finally, due to $\rho(M_uX^n,X^n)\le2w(X^n,1/u)$ and condition (i) we have 
\[\limsup_{u\to\infty}\limsup_{n\to\infty}\mathbb P\big(\rho(M_uX^n,X^n)\ge\epsilon\big)\le\limsup_{u\to\infty}\limsup_{n\to\infty}\mathbb P(2w(X^n,1/u)\ge\epsilon)=0.\]
It remains to apply Theorem \ref{3.2}.

\begin{lemma}\label{p12}
Let $P$ and $Q$ be two probability measures on $(\boldsymbol C,\mathcal C)$. If $P\pi^{-1}_{t_1,\ldots,t_k}=Q\pi^{-1}_{t_1,\ldots,t_k}$ for all $0\le t_1<\ldots<t_k\le1$, then 
$P=Q$.
 \end{lemma}
Proof.  Denote by $\mathcal C_f$ the collection of {\it cylinder sets} of the form
\[\pi^{-1}_{t_1,\ldots,t_k}(H)=\{y\in \boldsymbol C: (y(t_1),\ldots,y(t_k))\in H\},\qquad \qquad \qquad (*)\]
where $0\le t_1<\ldots<t_k\le1$ and a Borel subset $H\subset \boldsymbol R^k$. Due to the continuity of the projections we have $\mathcal C_f\subset\mathcal C$. 

It suffices to check, using Lemma \ref{PM.42}, that  $\mathcal C_f$ is a separating class.  Clearly, $\mathcal C_f$ is closed under formation of finite intersections. To show that $\sigma(\mathcal C_f)=\mathcal C$, observe that a closed ball centered at $x$ of radius $a$ can be represented as $\cap_r(y:|y(r)-x(r)|\le a)$, where $r$ ranges over rationals in [0,1]. It follows that $\sigma(\mathcal C_f)$ contains all closed balls, hence the open balls, and hence the $\sigma$-algebra generated by  the open balls.  By separability,  the $\sigma$-algebra generated by  the open balls, the so-called {\it ball  $\sigma$-algebra},  coincides with the Borel  $\sigma$-algebra generated by the open sets.

\begin{figure}
\centering
\includegraphics[height=6cm]{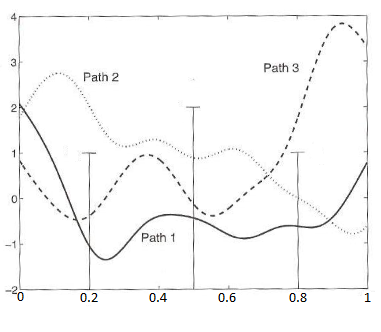}
\caption{Cylinder sets.}
\label{fcy}
\end{figure}

\begin{exercise}
Which of the three paths on Figure \ref{fcy} belong to the cylinder set $(*)$ with $k=3$, $t_1=0.2, t_2=0.5, t_3=0.8$, and $H=[-2,1]\times[-2,2]\times[-2,1]$.
\end{exercise}

\begin{theorem}\label{7.1}
Let $P_n$ be probability measures on $(\boldsymbol C,\mathcal C)$. If their finite-dimensional distributions converge weakly 
$P_{n}\pi^{-1}_{t_1,\ldots,t_k}\Rightarrow \mu_{t_1,\ldots,t_k}$, and if $P_n$ is tight, then 

(a) there exists a probability measure  $P$ on  $(\boldsymbol C,\mathcal C)$ with $P\pi^{-1}_{t_1,\ldots,t_k}=\mu_{t_1,\ldots,t_k}$, and

(b) $P_n\Rightarrow P$. 
\end{theorem}
Proof. Tightness implies relative compactness which in turn implies that each subsequence $(P_{n'})\subset (P_n)$ contains a further  subsequence $(P_{n''})\subset (P_{n'})$ converging weakly to some probability measure $P$. By the mapping theorem 
$P_{n''}\pi^{-1}_{t_1,\ldots,t_k}\Rightarrow P\pi^{-1}_{t_1,\ldots,t_k}$. Thus by hypothesis,
$P\pi^{-1}_{t_1,\ldots,t_k}=\mu_{t_1,\ldots,t_k}$. Moreover, by Lemma \ref{p12}, the limit $P$ must be the same for all converging subsequences, thus applying Theorem \ref{2.6} we may conclude that $P_n\Rightarrow P$.

\subsection{Wiener measure and Donsker's theorem}

\begin{definition}\label{p87}
Let $\xi_i$ be a sequence of r.v. defined on the same probability space $(\Omega,\mathcal{F},\mathbb P)$. 
Put $S_n=\xi_1+\ldots+\xi_n$ and let $X^n_t(\omega)$ as a function of $t$ be the element of $\boldsymbol C$ defined by linear interpolation between its values $X^n_{i/n}(\omega)={S_i(\omega)\over\sigma\sqrt n}$ at the points $t=i/n$. 
%Figure \ref{frw}  gives an illustration of this construction for the symmetric simple random walk with iid $\xi_i$ such that $\mathbb P(\xi_i=1)=\mathbb P(\xi_i=-1)=1/2$. 
\end{definition}

\begin{theorem} \label{8.2'}Let $X^n=(X^n_t:0\le t\le1)$ be defined by Definition \ref{p87} and let $P_n$ be the probability distribution of $X^n$. If $\xi_i$ are iid with zero mean and finite variance $\sigma^2$, then 

(a) $P_{n}\pi^{-1}_{t_1,\ldots,t_k}\Rightarrow \mu_{t_1,\ldots,t_k}$, where 
$\mu_{t_1,\ldots,t_k}$ are Gaussian distributions on $\boldsymbol R^k$ satisfying 
\[
\mu_{t_1,\ldots,t_k}\big\{(x_1,\ldots,x_k): x_i-x_{i-1}\le\alpha_i,i=1,\ldots,k\big\}
%=\prod_{i=1}^k{1\over\sqrt{2\pi(t_i-t_{i-1})}}\int_{-\infty}^{\alpha_i}e^{-{z^2\over2(t_i-t_{i-1})}}dz
=\prod_{i=1}^k\Phi\Big({\alpha_i\over\sqrt{t_i-t_{i-1}}}\Big), \mbox{ where }x_0=0,
\]

(b) the sequence $(P_n)$ of probability measures on $(\boldsymbol C,\mathcal C)$ is tight.
\end{theorem}
Proof. The claim (a) follows from the classical CLT and independence of increments of $S_n$. 
For example, if $0\le s\le t\le1$, then
\begin{align*}
 (X^n_s,X^n_t-X^n_s)&={1\over\sigma\sqrt{n}}(S_{\lfloor ns\rfloor},S_{\lfloor nt\rfloor}-S_{\lfloor ns\rfloor})+ \epsilon^n_{s,t},\\
 \epsilon^n_{s,t}&={1\over\sigma\sqrt{n}}(\{ ns\}\xi_{\lfloor ns\rfloor+1},\{ nt\}\xi_{\lfloor nt\rfloor+1}-\{ ns\}\xi_{\lfloor ns\rfloor+1}),
\end{align*}
where $\{ nt\}$ stands for the fractional part of $nt$. By the classical CLT and Theorem \ref{2.8}c, ${1\over\sigma\sqrt{n}}(S_{\lfloor ns\rfloor},S_{\lfloor nt\rfloor}-S_{\lfloor ns\rfloor})$ has $\mu_{s,t}$ as a limit distribution. Applying Corollary \ref{3.1} to $\epsilon^n_{s,t}$, we derive $P_{n}\pi^{-1}_{s,t}\Rightarrow \mu_{s,t}$. 

The proof of (b) is postponed until the next subsection.

\begin{definition}
Wiener measure $\mathbb W$ is a probability measure on $\boldsymbol C$ with $\mathbb W\pi^{-1}_{t_1,\ldots,t_k}= \mu_{t_1,\ldots,t_k}$ given by the formula in Theorem \ref{8.2'} part (a). The standard Wiener process $W$ is the random element on $(\boldsymbol C,\mathcal C,\mathbb W)$ defined by $W_t(x)=x(t)$.
\end{definition}
The existence of $\mathbb W$ follows from  Theorems \ref{7.1} and \ref{8.2'}.

\begin{theorem} \label{8.2}Let $X^n=(X^n_t:0\le t\le1)$ be defined by Definition \ref{p87}. If $\xi_i$ are iid with zero mean and finite variance $\sigma^2$, then $X^n$ converges in distribution to the standard Wiener process.
\end{theorem}
Proof 1. This is a corollary of Theorems \ref{7.1} and \ref{8.2'}.\\

\noindent Proof 2. An alternative proof is based on Theorem \ref{7.5}. We have to verify that  condition (i) of Theorem \ref{7.5} holds under the assumptions of Theorem \ref{8.2'}. To this end take $t_j=j\delta$, $j=0,\ldots, \delta^{-1}$ assuming $n\delta>1$. Then
\begin{align*}
\mathbb P(w(X^n,\delta)\ge3\epsilon)&\le \sum_{j=1}^{1/\delta}\mathbb P\Big(\sup_{t_{j-1}\le s\le t_j} |X^n_s-X^n_{t_{j-1}}| \ge \epsilon\Big)\\
&= \sum_{j=1}^{1/\delta}\mathbb P\Big(\max_{(j-1)n\delta\le k\le jn\delta} {|S_k-S_{(j-1)n\delta}|\over\sigma\sqrt {n}} \ge \epsilon\Big)= \sum_{j=1}^{1/\delta}\mathbb P\Big(\max_{k\le n\delta} |S_k|\ge \epsilon\sigma\sqrt {n}\Big)\\
&=\delta^{-1}\mathbb P\Big(\max_{k\le n\delta} |S_k|\ge \epsilon\sigma\sqrt {n}\Big)\le 3\delta^{-1}\max_{k\le n\delta}\mathbb P\Big( |S_k|\ge \epsilon\sigma
\sqrt {n}/3\Big),
\end{align*}
where the last is Etemadi's inequality:
\[\mathbb P\Big(\max_{k\le n} |S_k|\ge \alpha\Big)\le 3\max_{k\le n}\mathbb P\Big( |S_k|\ge \alpha/3\Big).\]
Remark: compare this with Kolmogorov's inequality
$\mathbb P(\max_{k\le n} |S_k|\ge \alpha)\le {n\sigma^2\over \alpha^2}$.

It suffices to check that assuming $\sigma=1$,
\[\lim_{\lambda\to\infty}\limsup_{n\to\infty}\lambda^2\max_{k\le n}\mathbb P\Big( |S_k|\ge \epsilon\lambda\sqrt n\Big)=0.\]
Indeed, by the classical CLT, 
$$\mathbb P(|S_k|\ge \epsilon\lambda\sqrt k)<4(1-\Phi(\epsilon\lambda))\le {6\over \epsilon^4\lambda^4}$$
for sufficiently large $k\ge k(\lambda\epsilon)$.
It follows,
\[\limsup_{n\to\infty}\lambda^2\max_{k(\lambda\epsilon)\le k\le n}\mathbb P\Big( |S_k|\ge \epsilon\lambda\sqrt n\Big)\le\limsup_{n\to\infty}\lambda^2\max_{k\ge k(\lambda\epsilon)}\mathbb P\Big( |S_k|\ge \epsilon\lambda\sqrt k\Big)\le{6\over \epsilon^4\lambda^2}.\]
On the other hand, by Chebyshev's inequality,
\[\limsup_{n\to\infty}\lambda^2\max_{k\le k(\lambda\epsilon)}\mathbb P\Big( |S_k|\ge \epsilon\lambda\sqrt n\Big)\le\limsup_{n\to\infty}{\lambda^2k(\lambda\epsilon)
\over \epsilon^2\lambda^2n}=0\]
finishing the proof of (i) of Theorem \ref{7.5}.

\begin{example}
We show that  $h(x)=\sup_tx(t)$ is a continuous mapping from $\boldsymbol C$ to $\boldsymbol R$. Indeed, if $h(x)\ge h(y)$, then there are $t_i$ such that
\[0\le h(x)-h(y)=x(t_1)-y(t_2)\le x(t_1)-y(t_1)\le \|x-y\|.\]
Thus, we have $|h(x)-h(y)|\le \rho(x,y)$ and continuity follows.
\end{example}

\begin{example}
Turning to the symmetric simple random walk, put $M_n=\max(S_0,\ldots,S_n)$. As we show later in Theorem \ref{(9.10)}, for any $b\ge0$,
\[\mathbb P(M_n\le b\sqrt n)\to{2\over\sqrt{2\pi}}\int_0^b e^{-u^2/2}du.\]
From $h(X^n)\Rightarrow h(W)$ with  $h(x)=\sup_tx(t)$ we conclude that $\sup_{0\le t\le1}W_t$ is distributed as $|W_1|$. The same limit holds for $M_n=\max({S_0\over\sigma},{S_1-\mu\over\sigma}\ldots,{S_n-n\mu\over\sigma})$ for sums of iid r.v. with mean $\mu$ and standard deviation $\sigma$. For this reason the functional CLT is also called an {\it invariance principle}: the general limit can be computed via the simplest relevant case.
\end{example}

\begin{exercise}
 Check if the following functionals are continuous on $\boldsymbol C$: 
 $$\sup_{\{0\le s,t\le1\}}|x(t)-x(s)|,\qquad \int_0^1 x(t)dt.$$
\end{exercise}

\subsection{Tightness in $\boldsymbol C$}

\begin{theorem}\label{7.2} The Arzela-Ascoli theorem. The set $A\subset \boldsymbol C$ is relatively compact if and only if 
\begin{align*}
(i)& \quad\sup_{x\in A}|x(0)|<\infty, %\label{s0} 
\\
(ii)& \quad \lim_{\delta\to0}\sup_{x\in A} w_x(\delta)=0. %\label{sCp} 
\end{align*}
 \end{theorem}
 Proof. {\it Necessity}. If the closure of $A$ is compact, then (i) obviously must hold. For a fixed $x$ the function $w_x(\delta)$ monotonely converges to zero as $\delta\downarrow0$. Since for each $\delta$ the function $w_x(\delta)$ is continuous in $x$ this convergence is uniform over $x\in K$ for any compact $K$. It remains to see that taking $K$ to be the closure of $A$ we obtain (ii).
 
{\it Sufficiency}. Suppose now that (i) and (ii) hold. For a given $\epsilon>0$, choose $n$ large enough for $\sup_{x\in A} w_x(1/n)<\epsilon$. Since
\[|x(t)|\le |x(0)|+\sum_{i=1}^n |x(ti/n)-x(t(i-1)/n)|\le  |x(0)|+n\sup_{x\in A} w_x(1/n),\]
we derive $\alpha:=\sup_{x\in A}\|x\|<\infty$. The idea is to use this and (ii) to prove that $A$ is totally bounded, since $\boldsymbol C$ is complete, it will follow that $A$ is relatively compact.
In other words, we have to find a finite $B_\epsilon\subset \boldsymbol C$ forming a $2\epsilon$-net for $A$. 

\begin{figure}
\centering
\includegraphics[width=12cm,height=4cm]{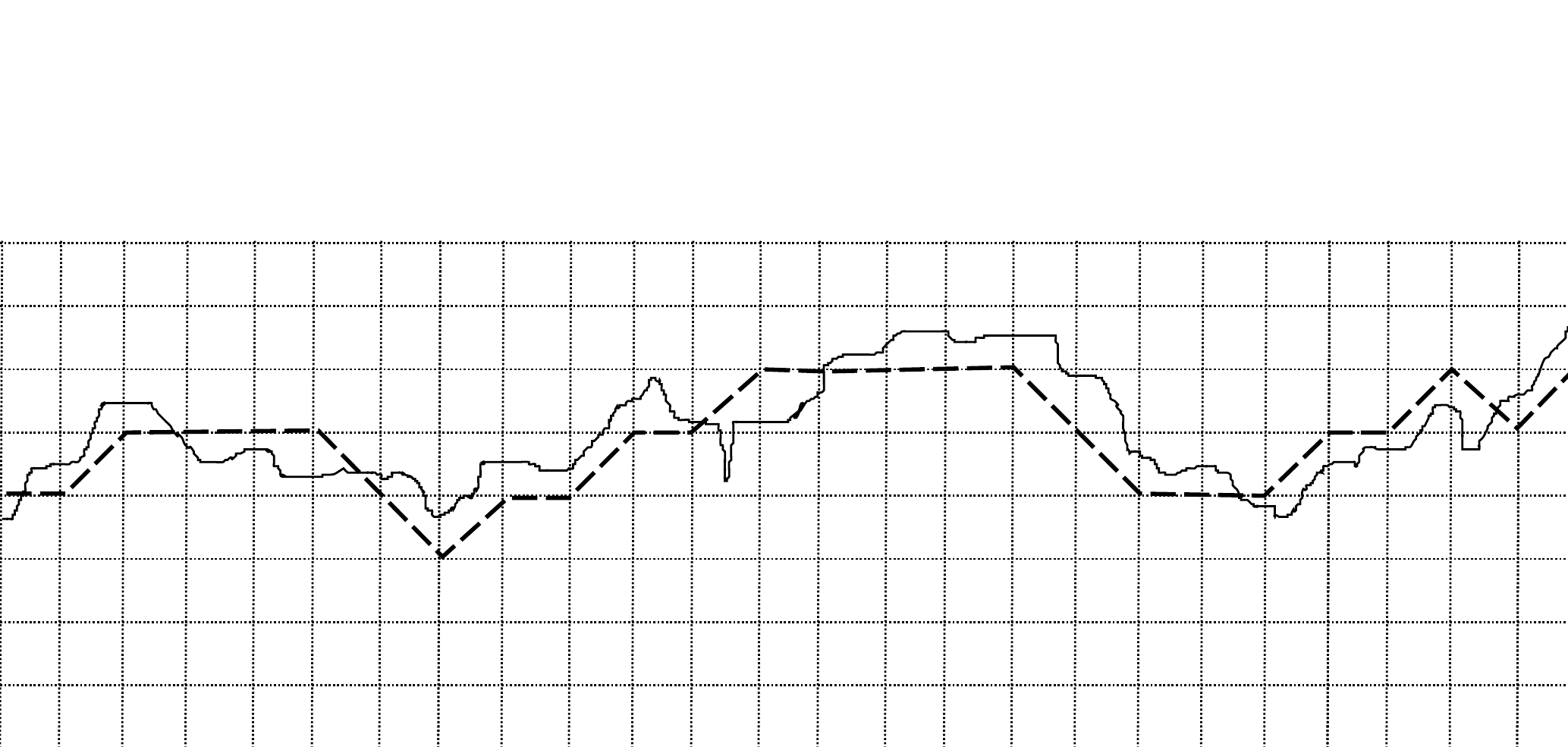}
\caption{The Arzela-Ascoli theorem: constructing a $2\epsilon$-net.}
\label{aras}
\end{figure}

Let $-\alpha=\alpha_0<\alpha_1<\ldots<\alpha_k=\alpha$ be such that $\alpha_j-\alpha_{j-1}\le\epsilon$. Then $B_\epsilon$ can be taken as a set of the continuous polygonal functions $y:[0,1]\to[-\alpha,\alpha]$ that linearly connect the pairs of points $({i-1\over n}, \alpha_{j_{i-1}}), ({i\over n}, \alpha_{j_i})$. See Figure \ref{aras}.
Let $x\in A$. It remains to show that there is a $y\in B_\epsilon$  such that $\rho(x,y)\le2\epsilon$. Indeed, since $|x(i/n)|\le\alpha$, there is a $y\in B_\epsilon$  such that $|x(i/n)-y(i/n)|<\epsilon$ for all $i=0,1,\ldots,n$. Both $y(i/n)$ and $y((i-1)/n)$ are within $2\epsilon$ of $x(t)$ for $t\in[(i-1)/n,i/n]$. Since $y(t)$ is a convex combination of $y(i/n)$ and $y((i-1)/n)$, it too is within $2\epsilon$ of $x(t)$. Thus $\rho(x,y)\le2\epsilon$ and $B_\epsilon$ is a $2\epsilon$-net for $A$. 

\begin{exercise}
 Draw a curve $x\in A$ (cf Figure \ref{aras}) for which you can not find a $y\in B_\epsilon$  such that $\rho(x,y)\le\epsilon$.
\end{exercise}
The next theorem explains the nature of  condition (i) in Theorem \ref{7.5}.
\begin{theorem}\label{7.3}
Let $P_n$ be probability measures on $(\boldsymbol C,\mathcal C)$. The sequence  $(P_n)$ is tight if and only if the following two conditions hold:
\begin{align*}
(i)& \quad \lim_{a\to\infty}\limsup_{n\to\infty}P_n(x: |x(0)|\ge a)=0, %\label{s0} 
\\
(ii)& \quad \lim_{\delta\to0}\limsup_{n\to\infty}P_n(x: w_x(\delta)\ge\epsilon)=0,\mbox{ for each positive }\epsilon. %\label{sCp} 
\end{align*}
\end{theorem}
Proof. Suppose $(P_n)$ is tight. Given a positive $\eta$, choose a compact $K$ such that $P_n(K)>1-\eta$ for all $n$. By the Arzela-Ascoli theorem we have $K\subset (x: |x(0)|\le a)$ for large enough $a$ and $K\subset (x: w_x(\delta)\le\epsilon)$ for small enough $\delta$. Hence the necessity.

According to condition (i), for each positive $\eta$, there exist large  $a_\eta$ and $n_\eta$ such that
\[P_n(x: |x(0)|\ge a_\eta)\le\eta,\quad n\ge n_\eta,\]
and condition (ii) implies that for each positive $\epsilon$ and $\eta$, there exist a small $\delta_{\epsilon,\eta}$ and a large $n_{\epsilon,\eta}$ such that
\[P_n(x: w_x(\delta_{\epsilon,\eta})\ge\epsilon)\le\eta,\quad n\ge n_{\epsilon,\eta},\]
Due to Lemma \ref{1.3} for any finite $k$ the measure $P_k$ is tight, and so by the necessity there is a $a_{k,\eta}$ such that $P_k(x: |x(0)|\ge a_{k,\eta})\le\eta$, and there is a $\delta_{k,\epsilon,\eta}$ such that $P_k(x: w_x(\delta_{k,\epsilon,\eta})\ge\epsilon)\le\eta$. 

Thus in proving sufficiency, we may put $n_\eta=n_{\epsilon,\eta}=1$ in the above two conditions. Fix an arbitrary small positive $\eta$. Given the two improved conditions, we have $P_n(B)\ge1-\eta$ and $P_n(B_{k})\ge1-2^{-k}\eta$ with $B=(x: |x(0)|< a_{\eta})$ and $B_{k}=(x: w_x(\delta_{1/k,2^{-k}\eta})<1/k)$. If $K$ is the closure of intersection of $B\cap B_1\cap B_{2}\cap\ldots$, then $P_n(K)\ge1-2\eta$. To finish the proof observe that $K$ is compact by the Arzela-Ascoli theorem.

\begin{example} 
 Consider the Dirac probability measure $P_n$ concentrated on the point $z_n\in \boldsymbol C$ from Example \ref{E1.3}. Referring to Theorem \ref{7.3} verify that the sequence $(P_n)$ is not tight. 
\end{example}

\noindent {\bf Proof of Theorem \ref{8.2'} part b}. The stated tightness follows from Theorem \ref{7.3}. Indeed, condition (i) in Theorem \ref{7.3} is trivially fulfilled as $X^n_0\equiv0$. Furthermore, condition (i) of Theorem \ref{7.5} (established in the proof 2 of Theorem \ref{8.2}) translates into  (ii) in Theorem \ref{7.3}.
%, and the one-dimensional convergence  $\{P_n\pi^{-1}_0\}\to\mu_0$ implies that $\{P_n\pi^{-1}_0\}$ is tight which in turn gives

\section{Applications of the functional CLT}
\subsection{The minimum and maximum of the Brownian path}
\begin{theorem}\label{(9.10)} Consider the standard Wiener process $W=(W_t, 0\le t\le 1)$ and let 
$$m=\inf_{0\le t\le1} W_t,\qquad M=\sup_{0\le t\le1} W_t.$$ 
If $a\le0\le b$ and $a\le a'<b'\le b$, then 
 \begin{align*}
 \mathbb P(&a<m\le M<b;\ a'<W_1<b')\\
  &=\sum_{k=-\infty}^\infty\Big(\Phi(2k(b-a)+b')-\Phi(2k(b-a)+a')\Big)\\
  &\qquad -\sum_{k=-\infty}^\infty\Big(\Phi(2k(b-a)+2b-a')-\Phi(2k(b-a)+2b-b')\Big),
\end{align*}
so that with $a'=a$ and $b'=b$ we get
\begin{align*}
 \mathbb P(&a<m\le M<b)=\sum_{k=-\infty}^\infty(-1)^k\Big(\Phi(k(b-a)+b)-\Phi(k(b-a)+a)\Big).
\end{align*}
\end{theorem}
Proof. Let $S_n$ be the symmetric simple random walk and put $m_n=\min(S_0,\ldots,S_n)$, $M_n=\max(S_0,\ldots,S_n)$. Since the mapping of $\boldsymbol C$ into $\boldsymbol R^3$ defined by
\[x\to\big(\inf_tx(t),\sup_tx(t),x(1)\big)\]
is continuous, the functional CLT entails $n^{-1/2}(m_n,M_n,S_n)\Rightarrow(m,M,W_1)$. The theorem's main  statement will be obtained in two steps. 

Step 1: show that for integers satisfying $i<0<j$ and $i\le i'<j'\le j$,
 \begin{align*}
 \mathbb P(i<m_n\le M_n<j;\ &i'<S_n<j')\\
 &=\sum_{k=-\infty}^\infty\mathbb P(2k(j-i)+i'<S_n<2k(j-i)+j')\\
 &\quad -\sum_{k=-\infty}^\infty\mathbb P(2k(j-i)+2j-j'<S_n<2k(j-i)+2j-i').
\end{align*}
In other words, we have to show that for  $i<0< j$, $i< l< j$
 \begin{align*}
 \mathbb P(i<m_n\le M_n<j;\ S_n=l)&=\sum_{k=-\infty}^\infty\mathbb P(S_n=2k(j-i)+l)\\
& \qquad -\sum_{k=-\infty}^\infty\mathbb P(S_n=2k(j-i)+2j-l).\qquad\qquad(*)
\end{align*}
Observe that here both series are just finite sums as $|S_n|\le n$.

Equality $(*)$ is proved by induction on $n$. For $n=1$, if $j>1$, then
\begin{align*}
 \mathbb P(i<m_1&\le M_1<j;\ S_1=1)=\mathbb P(S_1=1)\\
 &=\sum_{k=-\infty}^\infty\mathbb P(S_1=2k(j-i)+1) -\sum_{k=-\infty}^\infty\mathbb P(S_1=2k(j-i)+2j-1),
 \end{align*}
and if $i<-1$, then
\begin{align*}
 \mathbb P(i<m_1&\le M_1<j;\ S_1=-1)=\mathbb P(S_1=-1)\\
 &=\sum_{k=-\infty}^\infty\mathbb P(S_1=2k(j-i)-1) -\sum_{k=-\infty}^\infty\mathbb P(S_1=2k(j-i)+2j+1).
 \end{align*}
 Assume as induction hypothesis that the statement holds for $(n-1,i,j,l)$ with all relevant triplets $(i,j,l)$. Conditioning on the first step of the random walk, we get 
 \begin{align*}
 \mathbb P(i<m_n\le M_n<j;\ S_n=l)&={1\over2} \cdot \mathbb P(i-1<m_{n-1}\le M_{n-1}<j-1;\ S_{n-1}=l-1)\\
&+{1\over2} \cdot \mathbb P(i+1<m_{n-1}\le M_{n-1}<j+1;\ S_{n-1}=l+1),
\end{align*}
which together with the induction hypothesis yields the stated equality $(*)$
\begin{align*}
2 \mathbb P(&i<m_n\le M_n<j;\ S_n=l)\\
 &=\sum_{k=-\infty}^\infty\Big(\mathbb P(S_{n-1}=2k(j-i)+l-1)+\mathbb P(S_{n-1}=2k(j-i)+l+1)\Big) 
 \\
 &\quad -\sum_{k=-\infty}^\infty\Big(\mathbb P(S_{n-1}=2k(j-i)+2j-l+1)+\mathbb P(S_{n-1}=2k(j-i)+2j-l-1)\Big)
\\
 &=2\sum_{k=-\infty}^\infty\mathbb P(S_n=2k(j-i)+l) -2\sum_{k=-\infty}^\infty\mathbb P(S_n=2k(j-i)+2j-l).
\end{align*}

Step 2: show that for $c>0$ and $a<b$, 
 \begin{align*}
 \sum_{k=-\infty}^\infty\mathbb P\Big(2k\lfloor c\sqrt n\rfloor+\lfloor a\sqrt n\rfloor<&S_n<2k\lfloor c\sqrt n\rfloor+\lfloor b\sqrt n\rfloor\Big)\\
 & \to
\sum_{k=-\infty}^\infty\Big(\Phi(2kc+b)-\Phi(2kc+a)\Big),\quad n\to\infty.
\end{align*}
This is obtained using the CLT. The interchange of the limit with the summation over $k$ follows from
 \begin{align*}
\lim_{k_0\to\infty} \sum_{|k|>k_0}\mathbb P\Big(2k\lfloor c\sqrt n\rfloor+\lfloor a\sqrt n\rfloor<S_n<2k\lfloor c\sqrt n\rfloor+\lfloor b\sqrt n\rfloor\Big)=0,
\end{align*}
which in turn can be justified by the following series form of Scheffe's theorem. If $\sum_ks_{kn}=\sum_ks_{k}=1$, the terms being nonnegative, and if $s_{kn}\to s_k$ for each $k$, then $\sum_kr_ks_{kn}\to\sum_kr_ks_{k}$ provided $r_k$ is bounded. To apply this in our case we should take
\[s_{kn}=\mathbb P\Big(2k\lfloor \sqrt n\rfloor-\lfloor \sqrt n\rfloor<S_n\le2k\lfloor \sqrt n\rfloor+\lfloor \sqrt n\rfloor\Big), \quad s_k=\Phi(2k+1)-\Phi(2k-1).\]
\begin{corollary}\label{(9.14)} Consider the standard Wiener process $W$. If $a\le0\le b$, then 
  \begin{align*}
 \mathbb P(\sup_{0\le t\le1}W_t<b)&=2\Phi(b)-1,\\
\mathbb P(\inf_{0\le t\le1}W_t>a)&=1-2\Phi(a),\\
 \mathbb P(\sup_{0\le t\le1} |W_t|<b)&=2\sum_{k=-\infty}^\infty\Big\{\Phi((4k+1)b)-\Phi((4k-1)b)\Big\}.
\end{align*}
\end{corollary}

\subsection{The arcsine law}
\begin{lemma}\label{p247}
For $x\in \boldsymbol C$ and a Borel measurable, bounded $v:\boldsymbol R\to\boldsymbol R$, put $h(x)=\int_0^1v(x(t))dt$. If $v$ is continuous except on a set $D_v$ with 
$\lambda(D_v)=0$, where $\lambda$ is the Lebesgue measure, then $h$ is $\mathcal C$-measurable and is continuous except on a set of Wiener measure 0.
\end{lemma}
Proof. Since both mappings $x\to x(t)$ and $t\to x(t)$ are continuous, the mapping $(x,t)\to x(t)$ is continuous in the product topology and therefore Borel measurable.  It follows that the mapping $\psi(x,t)=v(x(t))$ is also measurable. Since $\psi$ is bounded, $h(x)=\int_0^1\psi(x,t)dt$ is $\mathcal C$-measurable, see Fubini's theorem.

Let $E=\{(x,t):x(t)\in D_v\}$. If $\mathbb W$ is Wiener measure on $(\boldsymbol C, \mathcal C)$, then by the hypothesis $\lambda(D_v)=0$,
\begin{align*}
&\mathbb W\{x:(x,t)\in E\}=\mathbb W\{x:x(t)\in D_v\}=0\mbox{ for each }t\in[0,1].
\end{align*}
It follows by Fubini's theorem applied to the measure $\mathbb W\times\lambda$ on $\boldsymbol C\times[0,1]$ that $\lambda\{t:(x,t)\in E\}=0$ for all $x$ outside a set $A_v\in\mathcal C$ satisfying $\mathbb W(A_v)=0$. Suppose that $\|x_n-x\|\to0$. If $x\notin A_v$, then $x(t)\notin D_v$ for almost all $t$ and hence $v(x_n(t))\to v(x(t))$ for almost all $t$. It follows by the bounded convergence theorem that
\[\mbox{if }x\notin A_v\mbox{ and } \|x_n-x\|\to0,\quad \mbox{ then }\int_0^1v(x_n(t))dt\to\int_0^1v(x(t))dt.\]

\begin{exercise}
Let $W$ be a standard Wiener process and  $t_0\in(0,1)$. Put $W'_s={W_t-W_{t_0}\over \sqrt{1-t_0}}$ for $s={t-t_0\over 1-t_0}$, $t\in[t_0,1]$. Using the Donsker invariance principle show that $(W'_s, 0\le s\le 1)$ is also distributed as a standard Wiener process.
\end{exercise}

\begin{lemma}\label{M15}
Each of the following three mappings  $h_i:\boldsymbol C\to\boldsymbol R$
\begin{align*}
 h_1(x)&=\sup\{t:x(t)=0, t\in[0,1] \},\\
 h_2(x)&=\lambda\{t: x(t)>0,t\in[0,1]\},\\
 h_3(x)&=\lambda\{t: x(t)>0,t\in[0,h_1(x)]\}
\end{align*}
is $\mathcal C$-measurable %. Moreover, let $P$ be a probability measure on $C$ such that $P\pi^{-1}_t$ is absolutely continuous with respect to Lebesgue measure for almost all $t$. Then each of the mappings is continuous except on a set of $P$-measure 0. In particular, each of the mappings is 
and continuous except on a set of Wiener measure 0.
\end{lemma}
Proof. Using the previous lemma with $v(z)=1_{\{z\in(0,\infty)\}}$ we obtain the assertion for $h_2$.

Turning to $h_1$, observe that 
$$\{x:h_1(x)<\alpha\}=\{x:x(t)>0, t\in[\alpha,1]\}\cup\{x:x(t)<0, t\in[\alpha,1]\}$$ 
is open and hence $h_1$ is measurable. If $h_1$ is discontinuous at $x$, then there exist $0<t_0<t_1<1$ such that  $x(t_1)=0$ and 
\[ \mbox{either }\quad x(t)>0 \mbox{ for all } t\in[t_0,1]\setminus\{t_1\}\quad\mbox{ or }\quad  x(t)<0 \mbox{ for all } 
 t\in[t_0,1]\setminus\{t_1\}.\]
That $h_1$ is continuous except on a set of Wiener measure 0 will therefore follow if we show that, for each $t_0$, the random variables
$$M_0=\sup\{W_t,t\in[t_0,1]\}\quad\mbox{ and }\quad \inf\{W_t,t\in[t_0,1]\}$$ 
have continuous distributions. By the last exercise and Theorem \ref{(9.10)}, $M'=M_0-W_{t_0}$ has a continuous distribution. Because $M'$ and $W_{t_0}$ are independent, we conclude that their sum also has a continuous distribution. The infimum is treated the same way.

Finally, for $h_3$, use the representation
\[h_3(x)=\psi(x,h_1(x)),\quad \mbox{ where }\psi(x,t)=\int_0^tv(x(u))du \quad \mbox{ with } v(z)=1_{\{z\in(0,\infty)\}}.\]

\begin{theorem}\label{(9.23)} Consider the standard Wiener process $W$ and let 

$T=h_1(W)$ be the time at which $W$ last passes through 0, 

$U=h_2(W)$ be the total amount of time $W$ spends above 0, and 

$V=h_3(W)$ be the total amount of time $W$ spends above 0 in the interval $[0,T]$. 

\noindent so that
$$U=V+(1-T)1_{\{W_1\ge0\}}.$$
Then the triplet $(T,V,W_1)$ has the joint density
\[f(t,v,z)=1_{\{0<v<t<1\}}g(t,z),\quad g(t,z)={1\over2\pi }{|z|\over t^{3/2}(1-t)^{3/2}}e^{-{z^2\over 2(1-t)}}.\]
In particular, the conditional distribution of $V$ given $(T,W_1)$ is uniform on $[0,T]$, and 
\[\mathbb P(T\le t)=\mathbb P(U\le t)={2\over\pi }\arcsin(\sqrt t), \quad 0<t<1.\] 
\end{theorem}
Proof. The main idea is to apply the invariance principle via the symmetric simple random walk $S_n$. We will use three properties of $S_n$ and its path functionals $(T_n,U_n,V_n)$. First, we need the local limit theorem for $p_n(i)=\mathbb P(S_n=i)$ similar to that of Example \ref{E3.4}: 
\[\mbox{if }{i\over\sqrt n}\to z,\mbox{ with } n-i \mbox{ being even}, \mbox{ then }{\sqrt n\over 2}p_n(i)\to{1\over\sqrt{2\pi}}e^{-z^2/2}.\]
Second, we need the fact that
\[\mathbb P(S_1\ge1,\ldots, S_{n-1}\ge1,S_n=i)={i\over n}p_n(i),\quad i\ge1.\]
The third fact we need is that if $S_{2n}=0$, then $U_{2n}=V_{2n}$ and
\[\mathbb P(V_{2n}=2j|S_{2n}=0)={1\over n+1},\quad j=0,1,\ldots,n.\]

Using these three facts we obtain that  for $0\le2j\le2k<n$ and $i\ge1$,
\begin{align*}
 \mathbb P(T_n=2k,&V_{2n}=2j,S_n=i)\\
& = \mathbb P(S_{2k}=0,V_{2n}=2j,S_{2k+1}\ge1,\ldots, S_{n-1}\ge1,S_n=i)\\
 &=\mathbb P(S_{2k}=0)\mathbb P(V_{2k}=2j|S_{2k}=0)\mathbb P(S_{2k+1}\ge1,\ldots, S_{n-1}\ge1,S_{n}=i|S_{2k}=0)\\
 &=p_{2k}(0){1\over k+1}{i\over n-2k}p_{n-2k}(i).
\end{align*}
We apply Theorem \ref{3.3} to the three-dimensional lattice of points $({2k\over n},{2j\over n},{i\over \sqrt n})$ for which $i\equiv n \mbox{ (mod 2)}$. The volume of the corresponding cell is ${2\over n}\cdot{2\over n}\cdot{2\over \sqrt n}=8n^{-5/2}$. If
\[{2k\over n}\to t,\quad{2j\over n}\to v,\quad{i\over\sqrt  n}\to z,\quad0<v<t<1,\quad z>0,\]
then 
\begin{align*}
 {n^{5/2}\over8}\mathbb P(T_n=2k,&V_{2n}=2j,S_n=i)\\
 &={\sqrt n\over\sqrt {2k}}{\sqrt {2k}\over2}p_{2k}(0){n\over 2(k+1)}{i\over \sqrt n}{n\over n-2k}{\sqrt n\over\sqrt {n-2k}}{\sqrt {n-2k}\over2}p_{n-2k}(i)\\
 &\to {1\over\sqrt{2\pi}}{1\over t^{3/2}}z{1\over (1-t)^{3/2}}{1\over\sqrt{2\pi}}e^{-{z^2\over2(1-t)}}=g(t,z).
\end{align*}
The same result holds for negative $z$ by symmetry.

The joint density of $(T,W_1)$ is $tg(t,z)1_{\{0<t<1\}}$, hence the marginal density for $T$ equals
\begin{align*}
f_T(t)&=\int_{-\infty}^\infty tg(t,z)dz=\int_{0}^\infty e^{-{z^2\over 2(1-t)}} {zdz\over \pi(1-t)^{3/2}t^{1/2}}={1\over\pi(1-t)^{1/2}t^{1/2}}
\end{align*}
implying
\[\mathbb P(T\le t)={2\over\pi }\arcsin(\sqrt t), \quad 0<t<1.\] 
Notice also that
\begin{align*}
G(u)&:=\int_{-\infty}^\infty \int_{u}^1g(t,z)dtdz= \int_{u}^1\int_{0}^\infty e^{-{z^2\over 2(1-t)}} {zdz\over \pi(1-t)^{3/2}t^{3/2}}dt\\
&=\int_{u}^1{dt\over\pi(1-t)^{1/2}t^{3/2}}=-{2\over\pi}\int_{u}^1{dt^{-1/2}\over\sqrt{1-t}}={2\over\pi}\int_1^{u^{-1/2}}{ydy\over\sqrt{y^2-1}}={2\over\pi}\sqrt{u^{-1}-1}.
\end{align*}
If $(T,W_1)=(t,z)$, then $U$ is distributed uniformly over $[1-t,1]$ for $z\ge0$, and  uniformly over $[0,t]$ for $z<0$:
\begin{align*}
\mathbb P(U\le u|T=t,W_1=z)&={u-1+t\over t}1_{\{u\in[1-t,1], z\ge0\}}+{u\over t}1_{\{u\in[0,t], z<0\}}+1_{\{u\in(t,1], z<0\}}.
\end{align*}
Thus the marginal distribution function of $U$ equals
\begin{align*}
\mathbb P(U\le u)&=\mathbb E\Big({u-1+T\over T}1_{\{u\in[1-T,1], W_1\ge0\}}+{u\over T}1_{\{u\in[0,T], W_1<0\}}+1_{\{u\in(T,1], W_1<0\}}\Big)\\
&= \int_0^\infty\int _{1-u}^1(u-1+t)g(t,z)dtdz+ \int_{-\infty}^0\int _{u}^1ug(t,z)dtdz+ \int_{-\infty}^0\int _0^{u}tg(t,z)dtdz\\
&={1\over2}\int _{1-u}^1f_T(t)dt+{u-1\over2}G(1-u)+{u\over2} G(u)+{1\over2}\int _0^{u}f_T(t)dt\\
&={1\over2}\mathbb P(T>1- u)+{1\over2}\mathbb P(T\le u)={2\over\pi }\arcsin(\sqrt u).
\end{align*}

\subsection{The Brownian bridge}
\begin{definition}
 The transformed standard Wiener process $W^\circ_t=W_t-tW_1$, $t\in[0,1]$, is called the standard  Brownian bridge.
\end{definition}
\begin{exercise}
 Show that the standard  Brownian bridge $W^\circ$ is a Gaussian process with zero mean and covariance $\mathbb E(W^\circ_sW^\circ_t)=s(1-t)$ for $s\le t$.

\end{exercise}
\begin{example}
Define $h:\boldsymbol C\to \boldsymbol C$ by $h(x(t))=x(t)-tx(1)$. This is a continuous mapping since $\rho(h(x),h(y))\le2\rho(x,y)$, and $h(X^n)\Rightarrow W^\circ$ by Theorem \ref{8.2}.
 \end{example}

\begin{theorem}\label{(9.32)}
 Let $P_\epsilon$ be the probability measure on $(\boldsymbol C,\mathcal C)$ defined by 
 \[P_\epsilon(A)=\mathbb P(W\in A|0\le W_1\le\epsilon),\quad A\in\mathcal C.\]
 Then $P_\epsilon\Rightarrow \mathbb W^\circ$ as $\epsilon\to0$, where $\mathbb W^\circ$ is the distribution of the Brownian bridge $W^\circ$.
\end{theorem}
Proof. We will prove that for every closed $F\in\mathcal C$
\[\limsup_{\epsilon\to0}\mathbb P(W\in F|0\le W_1\le\epsilon)\le \mathbb P(W^\circ\in F).\]
Using $W^\circ_t=W_t-tW_1$ we get $\mathbb E(W^\circ_tW_1)=0$ for all $t$. From the normality we conclude that $W_1$ is independent of each $(W^\circ_{t_1},\ldots,W^\circ_{t_k})$. Therefore,
\[\mathbb P(W^\circ\in A,W_1\in B)=\mathbb P(W^\circ\in A)\mathbb P(W_1\in B),\quad A\in \mathcal C_f, B\in \mathcal R,\]
and since $ \mathcal C_f$, the collection of finite-dimensional sets, see the proof of Lemma \ref{p12}, is a separating class, it follows 
\[\mathbb P(W^\circ\in A|0\le W_1\le\epsilon)=\mathbb P(W^\circ\in A),\quad A\in \mathcal C, \epsilon >0.\]
Observe that  $\rho(W,W^\circ)=|W_1|$. Thus, 
$$\{|W_1|\le\delta\}\cap\{W\in F\}\subset \{W^\circ\in F_\delta\},\quad \mbox{where }F_\delta=\{x:\rho(x,F)\le\delta\}.$$
Therefore, if $\epsilon<\delta$
\[\mathbb P(W\in F|0\le W_1\le\epsilon)\le\mathbb P(W^\circ\in F_\delta|0\le W_1\le\epsilon)=\mathbb P(W^\circ\in F_\delta),\]
leading to the required result
\[\limsup_{\epsilon\to0}\mathbb P(W\in F|0\le W_1\le\epsilon)\le\limsup_{\delta\to0}\mathbb P(W^\circ\in F_\delta)= \mathbb P(W^\circ\in F).\]

\begin{theorem} \label{(9.39)}
Distribution functions for several functionals of the Brownian bridge:
\begin{align*}
\mathbb P\big(a<\inf_t W^\circ_t\le\sup_t W^\circ_t\le b\big)&=\sum_{k=-\infty}^\infty \Big(e^{-2k^2(b-a)^2}-e^{-2(b+k(b-a))^2}\Big),\quad a<0<b,\\
 \mathbb P\big(\sup_t |W^\circ_t|\le b\big)&=1+2\sum_{k=1}^\infty (-1)^ke^{-2k^2b^2},\quad b>0,\\
 \mathbb P\big(\sup_t W^\circ_t\le b\big)&=\mathbb P\big(\inf_t W^\circ_t>- b\big)=1-e^{-2b^2},\quad b>0,\\
 \mathbb P\big( h_2(W^\circ)\le u\big)&=u,\quad u\in[0,1].
\end{align*}
\end{theorem}
Proof. The main idea of the proof is the following. Suppose that $h:\boldsymbol C\to\boldsymbol R^k$ is a measurable mapping and that the set $D_h$ of its discontinuities satisfies $\mathbb W^\circ(D_h)=0$. It follows by Theorem \ref{(9.32)} and the mapping theorem that 
\[\mathbb P(h(W^\circ)\le\alpha)=\lim_{\epsilon\to0}\mathbb P(h(W)\le\alpha|0\le W_1\le\epsilon).\]
Using either this or alternatively,
\[\mathbb P(h(W^\circ)\le\alpha)=\lim_{\epsilon\to0}\mathbb P(h(W)\le\alpha|-\epsilon\le W_1\le0)\]
one can find explicit forms for distributions connected with $W^\circ$.

Turning to Theorem \ref{(9.10)} with $a<0<b$ and $a'=0,b'=\epsilon$ we get
 \begin{align*}
 \mathbb P(a<m\le M<b;&\ 0<W_1<\epsilon)\\
 &=\sum_{k=-\infty}^\infty\Big(\Phi(2k(b-a)+\epsilon)-\Phi(2k(b-a))\Big)\\
  &\qquad-\sum_{k=-\infty}^\infty\Big(\Phi(2k(b-a)+2b)-\Phi(2k(b-a)+2b-\epsilon)\Big).
\end{align*}
This implies the first statement as
\[{\Phi(z+\epsilon)-\Phi(z)\over\epsilon}\to{e^{-z^2/2}\over\sqrt{2\pi}}.\]

As for the last statement, we need to show,  in terms of  $U=h_2(W)$, that
\[\lim_{\epsilon\to0}\mathbb P(U\le u|-\epsilon\le W_1\le0)=u,\]
or, in terms of  $V=h_3(W)$, that
\[\lim_{\epsilon\to0}\mathbb P(V\le u|-\epsilon\le W_1\le0)=u.\]
Recall that the distribution of $V$ for given $T$ and $W_1$ is uniform on $(0,T)$, in other words, $L=V/T$ is uniformly distributed on $(0,1)$ and is independent of $(T,W_1)$. Thus,
\begin{align*}
\mathbb P(V\le u|-\epsilon\le W_1\le0)&=\mathbb P(TL\le u|-\epsilon\le W_1\le0)\\
&=\int_0^1\mathbb P(T\le u/s|-\epsilon\le W_1\le0)ds\\
&=u+\int_u^1\mathbb P(T\le u/s|-\epsilon\le W_1\le0)ds.
\end{align*}
It remains to see that
\begin{align*}
\int_u^1\mathbb P(T\le u/s&|-\epsilon\le W_1\le0)ds={1\over\Phi(\epsilon)-\Phi(0)}\int_u^1\mathbb P(T\le u/s;-\epsilon\le W_1\le0)ds\\
&\le c\epsilon^{-1}\int_u^1\mathbb P(T\le r;-\epsilon\le W_1\le0){dr\over r^2}\le
 c\epsilon^{-1}u^{-2}\int_u^1\int_0^{r}\int_0^{\epsilon}tg(t,z)dzdtdr\\
& \le
 c_1\epsilon u^{-2}\int_0^{1}{dt\over t^{1/2}(1-t)^{1/2}}\to0,\quad \epsilon\to0.
\end{align*}

\section{The space $\boldsymbol D=\boldsymbol D[0,1]$}\label{secD}

\subsection{Cadlag functions}
\begin{definition}
 Let $\boldsymbol D=\boldsymbol D[0,1]$ be the space of functions $x:[0,1]\to\boldsymbol R$ that are right continuous and have left-hand limits.
\end{definition}
\begin{exercise}
 If $x_n\in\boldsymbol D$ and $\|x_n-x\|\to0$, then $x\in\boldsymbol D$.
\end{exercise}
For $x\in \boldsymbol D$ and $T\subset[0,1]$ we will use notation 
\[w_x(T)=w(x,T)=\sup_{s,t\in T}|x(t)-x(s)|,\]
and write $w_x[t,t+\delta]$ instead of  $w_x([t,t+\delta])$.
This should not be confused with the earlier defined modulus of continuity
\[w_x(\delta)=w(x,\delta)=\sup_{0\le t\le 1-\delta}w_x[t,t+\delta],\]
Clearly, if $T_1\subset T_2$, then $w_x(T_1)\le w_x(T_2)$. Hence $w_x(\delta)$ is monotone over $\delta$.
\begin{example}\label{frp}
 Consider $x_n(t)=$ the fractional part of $nt$. It has regular downward jumps of  size $1$. For example, $x_1(t)=t$ for $t\in[0,1)$, and $x_1(1)=0$. Another example: $x_2(t)=2t$ for $t\in[0,1/2)$, $x_2(t)=2t-1$ for $t\in[1/2,1)$, and $x_2(1)=0$. 
Placing an interval $[t,t+\delta]$ around a jump, we find $w_{x_n}(\delta)\equiv1$.
\end{example}
\begin{lemma}\label{L12.1}
Consider an arbitrary $x\in \boldsymbol D$. For each  $\epsilon>0$, there exist points $0=t_0<t_1<\ldots<t_v=1$ such that
 \[w_x[t_{i-1},t_i)<\epsilon,\quad i=1,2,\ldots,v.\]
It follows that $x$ is bounded, and that $x$ can be uniformly approximated by simple functions constant over intervals, so that it is Borel measurable.
It follows also  that $x$ has at most countably many jumps. 
\end{lemma}
Proof. To prove the first statement, let $t^\circ=t^\circ(\epsilon)$ be the supremum of those $t\in[0,1]$ for which $[0,t)$ can be decomposed into finitely many subintervals satisfying $w_x[t_{i-1},t_i)<\epsilon$. We show in three steps that $t^\circ=1$.

Step 1. Since $x(0)=x(0+)$, we have $w_x[0,\eta_0)<\epsilon$ for some small positive $\eta_0$. Thus $t^\circ>0$.

Step 2. Since $x(t^\circ-)$ exists, we have $w_x[t^\circ-\eta_1,t^\circ)<\epsilon$ for some small positive $\eta_1$, which implies that the interval $[0,t^\circ)$ can itself be so decomposed. 

Step 3. Suppose  $t^\circ=\tau$, where $\tau<1$. From $x(\tau)=x(\tau+)$ using the argument of the step 1, we see that according to the definition of  $t^\circ$ we should have $t^\circ>\tau$.

The last statement of the lemma follows from the fact that for any natural $n$, there exist  at most finitely many points $t$ at which $|x(t)-x(t-)|\ge n^{-1}$. 

\begin{exercise}
 Find a bounded function $x\notin \boldsymbol D$ with the following property: for any set  $0=t_0<t_1<\ldots<t_v=1$ there exists an $i$ such that  $w_x[t_{i-1},t_i)\ge1$.
\end{exercise}

\begin{definition}
Let $\delta\in(0,1)$. A set $0=t_0<t_1<\ldots<t_v=1$ is called $\delta$-sparse if $t_i-t_{i-1}>\delta$ for $i=1,\ldots,v$. Define an analog of the modulus of continuity $w_x(\delta)$ by
\begin{align*}
 w'_x(\delta)&=w'(x,\delta)=\inf_{\{t_i\}}\max_{1\le i\le v}w_x[t_{i-1},t_i),
\end{align*}
where the infimum extends over all $\delta$-sparse sets $\{t_i\}$. The function $w_x'(\delta)$ is called a {\it cadlag modulus} of $x$.
\end{definition}
\begin{exercise}
 Using Lemma \ref{L12.1} show that a function $x:[0,1]\to\boldsymbol R$ belongs to $\boldsymbol D$ if and only if $w_x'(\delta)\to0$.
\end{exercise}

\begin{exercise}
Compute  $w'_x(\delta)$ for $x=1_{[0,a)}$.
\end{exercise}
\begin{lemma}\label{p123}
For any $x$, $w'_x(\delta)$ is non-decreasing over $\delta$, and $w_x'(\delta)\le w_x(2\delta)$. Moreover, for any $x\in \boldsymbol D$,
\[ j_x\le w_x(\delta)\le 2 w_x'(\delta)+j_x,\quad j_x=\sup_{0<t\le1}|x(t)-x(t-)|.\]
\end{lemma}
Proof. Taking a $\delta$-sparse set with $t_i-t_{i-1}\le2\delta$ we get $w_x'(\delta)\le w_x(2\delta)$. To see that $w_x(\delta)\le 2 w_x'(\delta)+j_x$ take a $\delta$-sparse set such that $w_x[t_{i-1},t_i)\le w_x'(\delta)+\epsilon$ for all $i$. If $|t-s|\le\delta$, then $s,t\in[t_{i-1},t_{i+1})$ for some $i$ and $|x(t)-x(s)|\le 2(w_x'(\delta)+\epsilon)+j_x$.

\begin{lemma}\label{(12.28)} Considering triples $t_1,t,t_2$ in [0,1] put
\begin{align*}
 w''_x(\delta)&=w''(x,\delta)=\sup_{t_1\le t_2\le t_1+\delta}\sup_{t_1\le t\le t_2}\{|x(t)-x(t_1)|\wedge|x(t_2)-x(t)|%=\sup_{t_1\le t\le t_2\le t_1+\delta}\{|x(t)-x(t_1)|\wedge|x(t_2)-x(t)|
 \}.
\end{align*}
For any $x$, $w''_x(\delta)$ is non-decreasing over $\delta$, and $w''_x(\delta)\le w'_x(\delta)$.
\end{lemma}
Proof. Suppose that $w'_x(\delta)<w$ and $\{\tau_i\}$ be a $\delta$-sparse set such that $w_x[\tau_{i-1},\tau_i)<w$ for all $i$. If $t_1\le t\le t_2\le t_1+\delta$, then either $|x(t)-x(t_1)|<w$ or $|x(t_2)-x(t)|<w$. Thus $w''_x(\delta)<w$ and letting $w\downarrow w'_x(\delta)$ we obtain $w''_x(\delta)\le w'_x(\delta)$.
\begin{example}
For the functions $x_n(t)=1_{\{t\in[0,n^{-1})\}}$ and $y_n=1_{\{t\in[1-n^{-1},1]\}}$ we have $w''_{x_n}(\delta)=w''_{y_n}(\delta)=0$, although $w'_{x_n}(\delta)=w'_{y_n}(\delta)=1$ for $n\ge\delta^{-1}$.
\end{example}
\begin{exercise}
 Consider $x_n(t)$ from Example \ref{frp}. Find $w'_{x_n}(\delta)$ and $w''_{x_n}(\delta)$ for all $(n,\delta)$.
\end{exercise}
 \begin{figure}
\centering
\includegraphics[width=8cm,height=4cm]{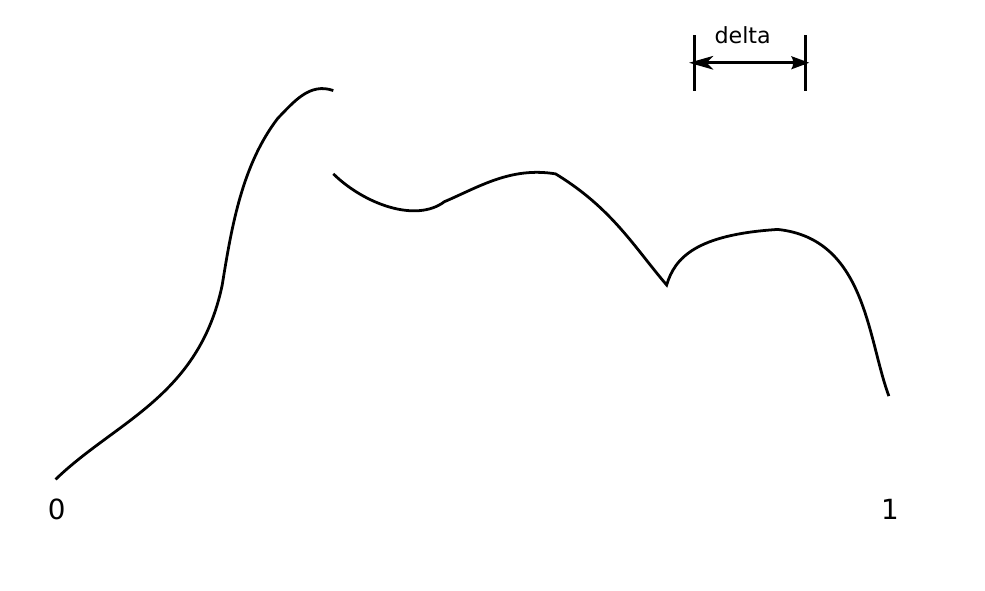}
\caption{Exercise \ref{xwx}.}
\label{wx}
\end{figure}
\begin{exercise}\label{xwx}
 Compare the values of $w'_{x}(\delta)$ and $w''_{x}(\delta)$ for the curve $x$ on the Figure \ref{wx}. 
\end{exercise}

\begin{lemma}\label{(12.32)}
For any $x\in \boldsymbol D$ and $\delta\in(0,1)$,
\begin{align*}
{w'_x(\delta/2)\over24} \le w''_x(\delta)\vee |x(\delta)-x(0)|\vee |x(1-)-x(1-\delta)|\le w'_x(\delta).
  \end{align*}
\end{lemma}
Proof. The second inequality follows from the definition of $w'_x(\delta)$ and Lemma \ref{(12.28)}. For the first inequality it suffices to show that
\begin{align*}
(i)&\quad w_x[t_1,t_2) \le 2(w''_x(\delta)+ |x(t_2)-x(t_1)|),\mbox{ if }t_2\le t_1+\delta,\\
(ii)&\quad w'_x(\delta/2) \le 6\Big(w''_x(\delta)\vee w_x[0,\delta)\vee  w_x[1-\delta,1)\Big),
  \end{align*}
  as these two relations imply
  \begin{align*}
 {w'_x(\delta/2)\over6} &\le w''_x(\delta)\vee w_x[0,\delta)\vee  w_x[1-\delta,1)\\
 &\le (2w''_x(\delta)+ 2|x(\delta)-x(0)|)\vee  (2w''_x(\delta)+ 2|x(1-)-x(1-\delta)|)\\
 &\le (4w''_x(\delta))\vee (4|x(\delta)-x(0)|)\vee(4 |x(1-)-x(1-\delta)|).
\end{align*}
Here we used the trick
\[w_x[1-\delta,1)=\lim_{t\uparrow1}w_x[1-\delta,t)\le2w''_x(\delta)+2\lim_{t\uparrow1} |x(t)-x(1-\delta)|.\]

To see (i), note that, if $t_1\le t< t_2\le t_1+\delta$, then either $|x(t)-x(t_1)|\le w''_x(\delta)$, or $|x(t_2)-x(t)|\le w''_x(\delta)$. In the latter case,  we have
\[|x(t)-x(t_1)|\le |x(t)-x(t_2)|+|x(t_2)-x(t_1)|\le w''_x(\delta)+ |x(t_2)-x(t_1)|.\]
Therefore, for $t_2\le t_1+\delta$,
\[\sup_{t_1\le t< t_2}|x(t)-x(t_1)|\le w''_x(\delta)+ |x(t_2)-x(t_1)|,\]
hence
\[\sup_{t_1\le s,t< t_2}|x(t)-x(s)|\le 2(w''_x(\delta)+ |x(t_2)-x(t_1)|).\]

We prove (ii) in four steps.

Step 1. We will need the following inequality
\[%(iii)\quad 
\qquad\qquad |x(s)-x(t_1)|\wedge|x(t_2)-x(t)|\le 2w''_x(\delta)\mbox{ if }t_1\le s<t\le t_2\le t_1+\delta.\qquad\qquad\qquad(*)\]
To see this observe that, by the definition of $w''_x(\delta)$, either $|x(s)-x(t_1)|\le w''_x(\delta)$ or both $|x(t_2)-x(s)|\le w''_x(\delta)$ and $|x(t)-x(s)|\le w''_x(\delta)$. 
In the second case, using the triangular inequality we get $|x(t_2)-x(t)|\le 2w''_x(\delta)$ .

Step 2. Putting
$$\alpha:=w''_x(\delta)\vee w_x[0,\delta)\vee  w_x[1-\delta,1),\qquad T_{x,\alpha}:=\{t:x(t)-x(t-)>2\alpha\},$$
show that there exist points $0=s_0<s_1<\ldots<s_r=1$ such that $s_i-s_{i-1}\ge\delta$ and 
\[T_{x,\alpha}\subset\{s_0,\ldots,s_r\}.\]

Suppose $u_1,u_2\in T_{x,\alpha}$ and $0<u_1<u_2<u_1+\delta$. Then there are disjoint intervals $(t_1,s)$ and $(t,t_2)$ such that $u_1\in(t_1,s)$, $u_2\in(t,t_2)$, and $t_2-t_1<\delta$. As both these intervals are short enough, we have a contradiction with $(*)$. Thus $(0,1)$ can not contain two points from $T_{x,\alpha}$, within $\delta$ of one another. And neither $[0,\delta)$ nor $[1-\delta,1)$ can contain a point from $T_{x,\alpha}$.

Step 3. Recursively adding middle points for the pairs $(s_{i-1},s_i)$ such that $s_i-s_{i-1}>\delta$ we get and enlarged set $\{s_0,\ldots,s_r\}$ (with possibly a larger $r$) satisfying
\[T_{x,\alpha}\subset\{s_0,\ldots,s_r\},\qquad \delta/2<s_i-s_{i-1}\le\delta,\quad i=1,\ldots,r.\]

Step 4. It remains to show that $w'_x(\delta/2) \le 6\alpha$. Since $\{s_0,\ldots,s_r\}$ from step 3 is a $(\delta/2)$-sparse set, it suffices to verify that 
$$w_x[s_{i-1},s_i)\le 6\alpha,\quad i=1,\ldots,r.$$ 
The proof will be completed after we demonstrate that 
$$|x(t_2)-x(t_1)|\le 6\alpha\quad  \mbox{ for }s_{i-1}\le t_1<t_2<s_i.$$ 
Define $\sigma_1$ and $\sigma_2$ by
\begin{align*}
\sigma_1&=\sup\{\sigma\in[t_1,t_2]: \sup_{t_1\le u\le\sigma}|x(u)-x(t_1)|\le2\alpha\},\\
\sigma_2&=\inf\{\sigma\in[t_1,t_2]: \sup_{\sigma\le u\le t_2}|x(t_2)-x(u)|\le2\alpha\}.
  \end{align*}
  If $\sigma_1<\sigma_2$, then there are $\sigma_1<s<t<\sigma_2$ violating $(*)$ due to the fact that by definition of $\alpha$, we have $w''_x(\delta)\le\alpha$. Therefore, $\sigma_2\le\sigma_1$ and it follows that 
  $|x(\sigma_1-)-x(t_1)|\le2\alpha$ and $|x(t_2)-x(\sigma_1)|\le2\alpha$. Since $\sigma_1\in (s_{i-1},s_i)$, we have $|x(\sigma_1)-x(\sigma_1-)|\le2\alpha$ implying 
  $$|x(t_2)-x(t_1)|\le |x(t_2)-x(\sigma_1)|+|x(\sigma_1)-x(\sigma_1-)|+|x(\sigma_1-)-x(t_1)|\le 6\alpha.$$

%\begin{lemma}\label{(12.32)'}
%For any $x\in \boldsymbol D$ and $\delta\in(0,1)$,
%\begin{align*}
%w'_x(\delta/2) \le 24w''_x(\delta)\vee |x(\delta)-x(0)|\vee |x(1)-x(1-\delta)|.
%  \end{align*}
%\end{lemma}
%Proof. The proof is similar to the proof of the first inequality stated in Lemma \ref{(12.32)}. It is based on the modified estimates
%\begin{align*}
%(i)&\quad w_x[t_1,t_2] \le 2(w''_x(\delta)+ |x(t_2)-x(t_1)|),\mbox{ if }t_2\le t_1+\delta,\\
%(ii)&\quad w'_x(\delta/2) \le 6\Big(w''_x(\delta)\vee w_x[0,\delta]\vee  w_x[1-\delta,1]\Big).
%  \end{align*}

\subsection{Two metrics in $\boldsymbol D$ and the Skorokhod topology}
\begin{example}
 Consider
$x(t)=1_{\{t\in[a,1]\}}$ and $y(t)=1_{\{t\in[b,1]\}}$ for $a,b\in(0,1)$. If $a\ne b$, then $\|x-y\|=1$ even when $a$ is very close to $b$. 
For the space $\boldsymbol D$,  the uniform metric is not good and we need another metric.
 \end{example}
\begin{definition}
Let $\Lambda$ denote the class of strictly increasing continuous mappings $\lambda:[0,1]\to[0,1]$ with $\lambda0=0$, $\lambda1=1$.  Denote by $1\in\Lambda$ the identity map $1 t\equiv t$, and put 
$ \|\lambda\|^\circ=\sup_{s<t}\big|\log{\lambda t-\lambda s\over t-s}\big|$. The smaller $ \|\lambda\|^\circ$ is the closer to 1 are the slopes of $\lambda$:
\[ e^{-\|\lambda\|^\circ}\le{\lambda t-\lambda s\over t-s}\le e^{\|\lambda\|^\circ}.\]
\end{definition}
\begin{exercise}
 Let $\lambda,\mu\in\Lambda$. Show that $$\|\lambda\mu-\lambda\|\le \|\mu-1\|\cdot e^{\|\lambda\|^\circ}.$$
\end{exercise}

\begin{definition}
For $x,y\in \boldsymbol D$ define  
 \begin{align*}
 d(x,y)&=\inf_{\lambda\in\Lambda}\{\|\lambda-1\|\vee\|x-y\lambda\|\},\\
 d^\circ(x,y)&=\inf_{\lambda\in\Lambda}\{\|\lambda\|^\circ\vee\|x-y\lambda\|\},
\end{align*}
\end{definition}
\begin{exercise}
 Show that $d(x,y)\le \|x-y\|$ and $d^\circ(x,y)\le \|x-y\|$.
\end{exercise}

\begin{example}
 Consider
$x(t)=1_{\{t\in[a,1]\}}$ and $y(t)=1_{\{t\in[b,1]\}}$ for $a,b\in(0,1)$. Clearly, if $\lambda(a)=b$, then $\|x-y\lambda\|=0$ and otherwise $\|x-y\lambda\|=1$. Thus
 \begin{align*}
 d(x,y)&=\inf\{\|\lambda-1\|:\lambda\in\Lambda, \lambda(a)=b\}=|a-b|,\\
 d^\circ(x,y)&=\Big(\inf\{\|\lambda\|^\circ:\lambda\in\Lambda, \lambda(a)=b\}\Big)\wedge1=\Big(\big|\log{a\over b}\big|\vee\big|\log{1-a\over 1-b}\big|\Big)\wedge1,
\end{align*}
so that $d(x,y)\to0$ and $d^\circ(x,y)\to0$ as $b\to a$.
 \end{example}
 \begin{exercise}\label{abc}
 Given $0<b<a<c<1$, find $d(x,y)$ for 
$$x(t)=2\cdot1_{\{t\in[a,1]\}},\qquad y(t)=1_{\{t\in[b,1]\}}+1_{\{t\in[c,1]\}}.$$
Does $d(x,y)\to0$ as $b\to a$ and $c\to a$?
\end{exercise}

\begin{lemma}\label{p126}
 Both $d$ and $d^\circ$ are metrics  in $\boldsymbol D$, and  $d\le e^{d^\circ}-1$.
\end{lemma}
Proof. Note that $d(x,y)$ is the infimum of those $\epsilon>0$ for which there exists a $\lambda\in\Lambda$ with
 \begin{align*}
 &\sup_t|\lambda t-t|=\sup_t|t-\lambda^{-1}t|<\epsilon,\\
 &\sup_t|x(t)-y(\lambda t)|=\sup_t|x(\lambda^{-1}t)-y(t)|<\epsilon.
\end{align*}
Of course $d(x,y)\ge0$, $d(x,y)=0$ implies $x=y$, and $d(x,y)=d(y,x)$. To see that $d$ is a metric we have to check the triangle inequality $d(x,z)\le d(x,y)+d(y,z)$. It follows from
 \begin{align*}
 &\|\lambda_1\lambda_2-1\|\le\|\lambda_1-1\|+\|\lambda_2-1\|,\\
 &\|x-z\lambda_1\lambda_2\|\le\|x-y\lambda_2\|+\|y-z\lambda_1\|.
\end{align*}
Symmetry and the triangle inequality for $d^\circ$ follows from $ \|\lambda^{-1}\|^\circ= \|\lambda\|^\circ$ and the inequality $$ \|\lambda_1\lambda_2\|^\circ\le\|\lambda_1\|^\circ+\|\lambda_2\|^\circ.$$ That $d^\circ(x,y)=0$ implies $x=y$ follows from 
$d\le e^{d^\circ}-1$ which is a consequence of $\|x-y\lambda\|\le e^{\|x-y\lambda\|}-1$ and 
\[\|\lambda-1\|=\sup_{0\le t\le 1}t\Big| {\lambda t-\lambda0\over t-0}-1\Big|\le e^{\|\lambda\|^\circ}-1.\]
The last inequality uses $|u-1|\le e^{|\log u|}-1$ for $u>0$.
\begin{example}
 Consider $j_x$, the maximum jump in $x\in \boldsymbol D$. Clearly, $|j_x-j_y|<\epsilon$ if $\|x-y\|<\epsilon/2$, and so $j_x$ is continuous in the uniform topology. It is also continuous in the Skorokhod topology. Indeed, if $d(x,y)<\epsilon/2$, then there is a $\lambda$ such that $\|\lambda-1\|<\epsilon/2$ and $\|x-y\lambda\|<\epsilon/2$. Since  $j_y=j_{y\lambda}$, we conclude using continuity in the uniform topology
 $|j_x-j_{y}|=|j_x-j_{y\lambda}|<\epsilon$.
\end{example}
\begin{lemma}\label{12.2}
If $d(x,y)=\delta^2$ and $\delta\le1/3$,  then $d^\circ(x,y)\le4\delta+w'_x(\delta)$.
\end{lemma}
Proof. We prove that if $d(x,y)<\delta^2$ and $\delta\le1/3$,  then $d^\circ(x,y)<4\delta+w'_x(\delta)$.
Choose  $\mu\in\Lambda$ such that $\|\mu-1\|<\delta^2$ and $\|x\mu^{-1}-y\|<\delta^2$. Take $\{t_i\}$ to be a $\delta$-sparse set satisfying $w_x[t_{i-1},t_i)\le w'_x(\delta)+\delta$ for each $i$. Take $\lambda$ to agree with $\mu$ at the points $\{t_i\}$ and to be linear in between. Since $\mu^{-1}\lambda t_i=t_i$, we have $t\in[t_{i-1},t_i)$ if and only if $\mu^{-1}\lambda t\in[t_{i-1},t_i)$, and therefore
\[|x(t)-y(\lambda t)|\le |x(t)-x(\mu^{-1}\lambda t)|+|x(\mu^{-1}\lambda t)-y(\lambda t)|< w'_x(\delta)+\delta+\delta^2\le4\delta+w'_x(\delta).\]
Now it is enough to verify that $\|\lambda\|^\circ<4\delta$. Draw a picture to see that the slopes of $\lambda$ are always between ${\delta\pm2(\delta^2)\over \delta}=1\pm2\delta$. Since $|\log(1\pm2\delta)|< 4\delta$ for sufficiently small $\delta$, we get $\|\lambda\|^\circ<4\delta$.

\begin{theorem}\label{12.1}
The metrics  $d$ and $d^\circ$ are equivalent and generate the same, so called {\it Skorokhod topology}. 
\end{theorem}
Proof. By definition $d(x_n,x)\to0$ ($d^\circ(x_n,x)\to0$) if and only if there is a sequence $\lambda_n\in\Lambda$ such that 
$\|\lambda_n-1\|\to0$ ($\|\lambda_n\|^\circ\to0$) and  $\|x_n\lambda_n-x\|\to0$. If $d^\circ(x_n,x)\to0$, then $d(x_n,x)\to0$ due to $d\le e^{d^\circ}-1$.
The reverse implication follows from Lemma \ref{12.2}.

\begin{definition}
Denote by  $\mathcal D$ the   Borel  $\sigma$-algebra formed from the open and closed sets  in $({\boldsymbol D},d)$, or equivalently $({\boldsymbol D},d^\circ)$, using the operations of countable intersection, countable union, and set difference.
\end{definition}

\begin{lemma}\label{(12.14)}
 Skorokhod convergence $x_n\to x$ in ${\boldsymbol D}$ implies $x_n(t)\to x(t)$ for continuity points $t$ of $x$. Moreover, if $x$ is continuous on $[0,1]$, then Skorokhod convergence implies uniform convergence.
\end{lemma}
Proof. Let  $\lambda_n\in\Lambda$ be such that 
$\|\lambda_n-1\|\to0$ and  $\|x_n-x\lambda_n\|\to0$. The first assertion follows from
\[|x_n(t)-x(t)|\le |x_n(t)-x(\lambda_nt)|+|x(\lambda_nt)-x(t)|.\]
The second assertion is obtained from 
\[\|x_n-x\|\le \|x_n-x\lambda_n\|+w_x(\|\lambda_n-1\|).\]
\begin{example}
Put $x_n(t )=1_{\{t\in[a-2^{-n},1]\}}+1_{\{t\in[a+2^{-n},1]\}}$ and $x(t )=2\cdot 1_{\{t\in[a,1]\}}$ for some $a\in(0,1)$. We have $x_n(t)\to x(t)$ for continuity points $t$ of $x$, however $x_n$ does not converge to $x$ in the Skorokhod topology.
\end{example}

\begin{exercise}
Fix $\delta\in(0,1)$ and consider $w'_x(\delta)$  as a function of $x\in\boldsymbol D$. 

(i) The function $w'_x(\delta)$ is continuous with respect to the uniform metric
since
$$|w'_x(\delta) - w'_y(\delta)|\le2 ||x-y||.$$
Hint: show that  $w_x[t_{i-1},t_i)\le w_y[t_{i-1},t_i)+2 ||x-y||$.

(ii) However $w'_x(\delta)$ is not continuous  with respect to  $(\boldsymbol D,d)$. Verify this using $x_n=1_{[0,\delta+2^{-n})}$ and $x=1_{[0,\delta)}$.
\end{exercise}

\begin{exercise}
Show that  $h(x)=\sup_tx(t)$ is a continuous mapping from $\boldsymbol D$ to $\boldsymbol R$. Hint: show first that for any $\lambda\in\Lambda$, 
\[|h(x)-h(y)|\le \|x-y\lambda\|.\]
\end{exercise}

\subsection{Separability and completeness of $\boldsymbol D$}
\begin{figure}
\centering
\includegraphics[width=6cm]{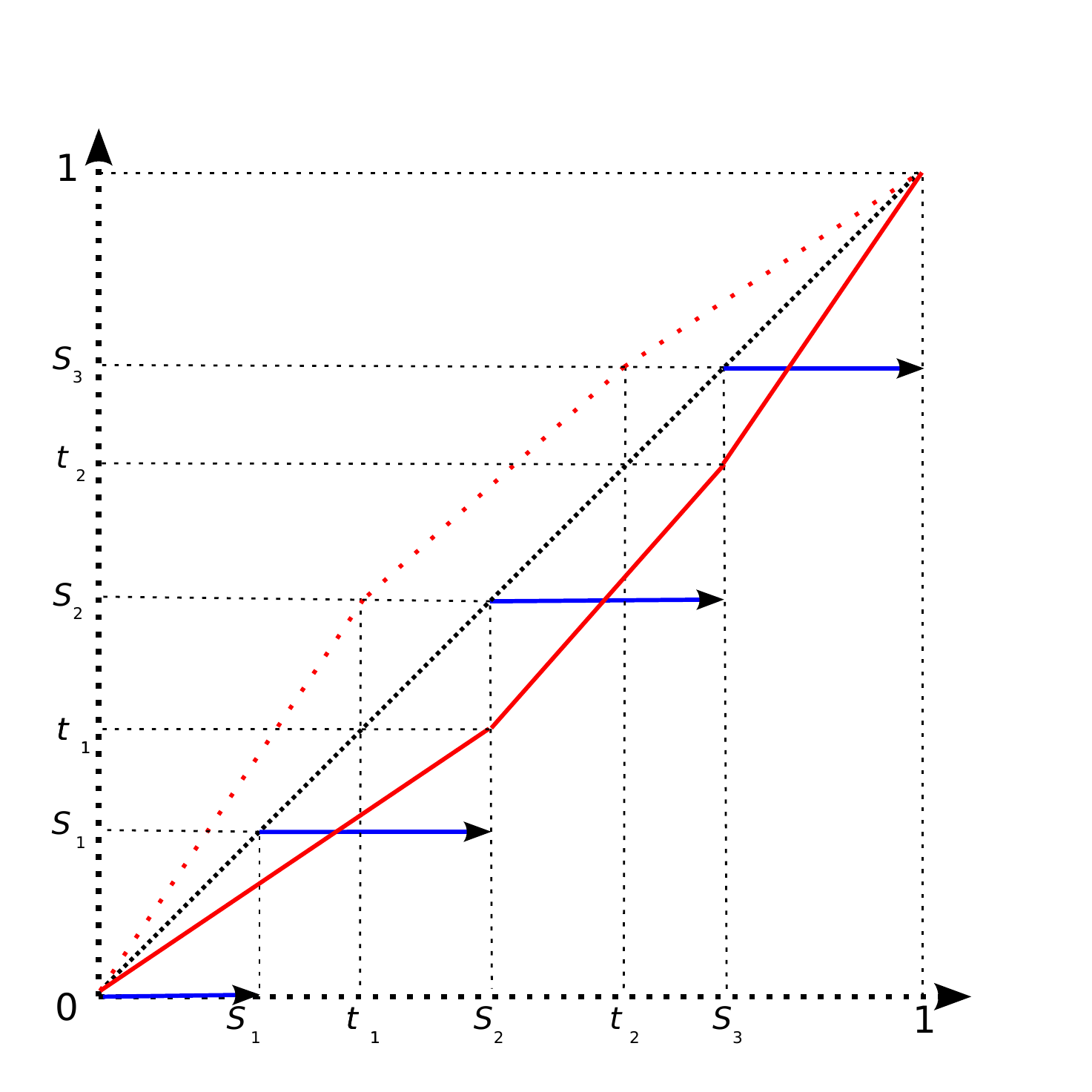}
\caption{Lemma  \ref{L12.3}. Comparing the blue graph of the function $\kappa$ to the red graph of the function $\lambda^{-1}$.}
\label{Figa}
\end{figure}

\begin{lemma}\label{L12.3}
Given $0=s_0<s_1<\ldots<s_k=1$ define a non-decreasing map\\ $\kappa:[0,1]\to [0,1]$ by setting 
\[
\kappa t=\left\{
\begin{array}{ll}
 s_{j-1} &   \mbox{for }  t\in[s_{j-1},s_j),\ j=1,\ldots,k,\\
1  &    \mbox{for }  t=1.
\end{array}
\right.
\]
If $\max_j(s_j-s_{j-1})\le\delta$, then $d( x\kappa,x)\le \delta\vee w'_x(\delta)$ for any $x\in \boldsymbol D$.
\end{lemma}
Proof. Given $\epsilon>0$ find a $\delta$-sparse set $\{t_i\}$  satisfying $w_x[t_{i-1},t_i)\le w'_x(\delta)+\epsilon$ for each $i$. Let $\lambda\in\Lambda$ be linear between $\lambda t_i=s_j$ for $t_i\in(s_{j-1},s_j]$, $j=1,\ldots,k$ and $\lambda 0=0$. Since $\|\lambda-1\|\le\delta$, it suffices to show that $|x(\kappa t)-x(\lambda^{-1}t)|\le w'_x(\delta)+\epsilon$. This holds if $t$ is 0 or 1, and it is enough to show that, for $t\in(0,1)$, both $\kappa t$ and  $\lambda^{-1} t$ lie in the same $[t_{i-1},t_i)$, see Figure \ref{Figa}. We prove this by showing that  $\kappa t<t_i$ is equivalent to $\lambda^{-1} t<t_i$, $i=1,\ldots,k$. Suppose that $t_i\in(s_{j-1},s_j]$. Then 
\[\kappa t<t_i\quad\Rightarrow\quad \kappa t<s_j\quad\Rightarrow\quad \kappa t\le s_{j-1}\quad\Rightarrow\quad \kappa t<t_i.\]
Thus $\kappa t<t_i$ is equivalent to $\kappa t<s_j$ which in turn  is equivalent to $ t<s_j$. On the other hand, 
$\lambda t_i=s_j$, and hence  $ t<s_j$ is equivalent to $t<\lambda t_i$ or $\lambda^{-1} t<t_i$.
\begin{example}
 Let $x_n(t)=1_{\{t\in[t_0,t_0+2^{-n})\}}$. Since $d(x_n,x_{n+1})=2^{-n}-2^{-n-1}=2^{-n-1}$, the sequence $(x_n)$ is fundamental. However, it is not $d$-convergent. Indeed, $x_n(t)\to0$ for all $t\ne t_0$ and Skorokhod convergence $x_n\to x$  in ${\boldsymbol D}$ by Lemma \ref{(12.14)}, should imply $x(t)=0$ for all points of continuity of $x$. Since $x\in{\boldsymbol D}$ has at most countably many points of discontinuities, by right continuity we conclude  that $x\equiv0$. Moreover, since the limit $x\equiv0$ is continuous, we must have $\|x_n-x\|\to0$. But $\|x_n\|\equiv1$.
\end{example}

\begin{theorem}\label{12.2}
The space $\boldsymbol D$ is separable under  $d$ and $d^\circ$, and is complete under  $d^\circ$.
\end{theorem}
Proof. $Separability$ for $d$. Put $s_j=j/k$, $j=1,\ldots,k$. Let $B_k$ be the set of functions having a constant, rational value over each $[s_{j-1},s_j)$ and a rational value at $t=1$. Then $B=\cup B_k$ is countable. Now it is enough to prove that given $x\in{\boldsymbol D}$ and $\epsilon>0$ we can find some $y\in B$ such that $d(x,y)<2\epsilon$. Choosing $k$ such that  $k^{-1}<\epsilon$ and $w'_x(k^{-1})<\epsilon$ we can find  $y\in B_k$ satisfying $d( x\kappa,y)<\epsilon$, for $\kappa$ defined as in  Lemma \ref{L12.3}. It remains to see that $d( x\kappa,x)<\epsilon$ according to Lemma \ref{L12.3}.

$Completeness$. We show that any $d^\circ$-fundamental sequence $x_n\in \boldsymbol D$ contains a subsequence $y_k=x_{n_k}$ that is $d^\circ$-convergent. Choose $n_k$ in such a way that $d^\circ(y_k,y_{k+1})<2^{-k}$. Then $\Lambda$ contains $\mu_k$ such that $\|\mu_k\|^\circ<2^{-k}$ and $\|y_k\mu_k-y_{k+1}\|<2^{-k}$.

We suggest a choice of $\lambda_k\in\Lambda$ such that $\|\lambda_k\|^\circ\to0$ and $\|y_k\lambda_k-y\|\to0$ for some $y\in \boldsymbol D$. To this end put $\mu_{k,m}=\mu_k\mu_{k+1}\ldots \mu_{k+m}$. From
\begin{align*}
 \| \mu_{k,m+1}-\mu_{k,m}\|&\le\|  \mu_{k+m+1}-1\|\cdot e^{\|\mu_k\mu_{k+1}\ldots \mu_{k+m}\|^\circ}\\
 &\le(e^{\| \mu_{k+m+1}\|^\circ}-1)\cdot e^{\|\mu_k\|^\circ+\ldots +\|\mu_{k+m}\|^\circ}\\
 &\le 2^{-k-m}\cdot e^{2^{-k}+\ldots+2^{-k-m}}<2^{-k-m+2}
\end{align*}
we conclude that for a fixed $k$ the sequence of functions $\mu_{k,m}$ is uniformly fundamental. Thus there exists a $\lambda_k$ such that $\|\mu_{k,m}-\lambda_k\|\to 0$ as $m\to\infty$. To prove that $\lambda_k\in\Lambda$ we use 
\begin{align*}
\log \Big| {\mu_{k,m}t- \mu_{k,m}s\over t-s}\Big|\le\|\mu_k\mu_{k+1}\ldots \mu_{k+m}\|^\circ<2^{-k+1}.
\end{align*}
Letting here $m\to\infty$ we get $\|\lambda_k\|^\circ\le2^{-k+1}$. Since $\|\lambda_k\|^\circ$ is finite we conclude that $\lambda_k$ is strictly increasing and therefore $\lambda_k\in\Lambda$. 

Finally, observe that 
\[\|y_k\lambda_k-y_{k+1}\lambda_{k+1}\|=\|y_k\mu_k\lambda_{k+1}-y_{k+1}\lambda_{k+1}\|=\|y_k\mu_k-y_{k+1}\|<2^{-k}.\]
It follows, that the sequence $y_k\lambda_k\in \boldsymbol D$ is uniformly fundamental and hence $\|y_k\lambda_k-y\|\to0$ for some $y$. Observe that $y$ must lie in $\boldsymbol D$. Since $\|\lambda_k\|^\circ\to0$, we obtain $d^\circ(y_k,y)\to0$.

\subsection{Relative compactness in the Skorokhod topology}
First comes an analogue of the Arzela-Ascoli theorem in terms of $w'_x(\delta)$, and then a convenient alternative in terms of $w''_x(\delta)$.
\begin{theorem} \label{12.3}
A set $A\subset \boldsymbol D$ is relatively compact in the Skorokhod topology  iff 
\begin{align*}
(i)& \quad \sup_{x\in A}\|x\|<\infty,  \\
(ii)& \quad \lim_{\delta\to0}\sup_{x\in A} w'_x(\delta)=0. %\label{sCp} 
\end{align*}
 \end{theorem}
 Proof of sufficiency only. Put $\alpha=\sup_{x\in A}\|x\|$.  For a given $\epsilon>0$,

put $H_\epsilon=\{\alpha_i\}$, where $-\alpha=\alpha_0<\alpha_1<\ldots<\alpha_k=\alpha$ and $\alpha_j-\alpha_{j-1}\le\epsilon$,

and choose $\delta<\epsilon$ so that $w'_x(\delta)<\epsilon$ for all $x\in A$. 

\noindent According to Lemma \ref{L12.3} for any $\kappa=\kappa_{\{s_j\}}$ satisfying $\max_j(s_{j-1}-s_j)\le\delta$, we have $d( x\kappa,x)\le \epsilon$ for all $x\in A$. Take  
$B_\epsilon$ be the set of $y\in \boldsymbol D$ that assume on each $[s_{j-1},s_j)$ a constant value from $H_\epsilon$ and $y(1)\in H_\epsilon$. For any $x\in A$ there is a $y\in B_\epsilon$ such that $d( x\kappa,y)\le \epsilon$. Thus $B_\epsilon$ forms a $2\epsilon$-net for $A$ in the sense of $d$ and $A$ is  totally bounded in the sense of $d$. 

But we must show that $A$ is  totally bounded in the sense of $d^\circ$, since this is the metric under which $\boldsymbol D$ is complete. This is true as according Lemma \ref{12.2}, the set $B_{\delta^2}$ is an $\epsilon'$-net for $A$, where $\epsilon'=4\delta+\sup_{x\in A}w'_x(\delta)$ can be chosen arbitrary small.

\begin{theorem} \label{12.4}
A set $A\subset \boldsymbol D$ is relatively compact in the Skorokhod topology  iff 
\begin{align*}
(i)& \quad \sup_{x\in A}\|x\|<\infty,  \\
(ii)& \quad 
\left\{
\begin{array}{l}
 \lim_{\delta\to0}\sup_{x\in A} w''_x(\delta)=0,\\
 \lim_{\delta\to0}\sup_{x\in A} |x(\delta)-x(0)|=0,\\
 \lim_{\delta\to0}\sup_{x\in A} |x(1-)-x(1-\delta)|=0.
\end{array}
\right.
 %\label{sCp} 
\end{align*}
 \end{theorem}
 Proof. It is enough to show that (ii) of Theorem \ref{12.4} is equivalent to (ii) of Theorem \ref{12.3}. This follows from Lemma \ref{(12.32)}.
%\begin{corollary} \label{12.4'}
%A set $A\subset \boldsymbol D$ is relatively compact in the Skorokhod topology  provided
%\begin{align*}
%(i)& \quad \sup_{x\in A}\|x\|<\infty,  \\
%(ii)& \quad 
%\left\{
%\begin{array}{l}
% \lim_{\delta\to0}\sup_{x\in A} w''_x(\delta)=0,\\
% \lim_{\delta\to0}\sup_{x\in A} |x(\delta)-x(0)|=0,\\
% \lim_{\delta\to0}\sup_{x\in A} |x(1)-x(1-\delta)|=0.
%\end{array}
%\right.
% %\label{sCp} 
%\end{align*}
% \end{corollary}
% Proof. Apply Lemma \ref{(12.32)'}.

\section{Probability measures on $\boldsymbol D$ and random elements}

\subsection{Finite-dimensional distributions on $\boldsymbol D$}
Finite-dimensional sets play in $\boldsymbol D$ the same role as they do in $\boldsymbol C$.
\begin{definition}
Consider projection mappings $\pi_{t_1,\ldots,t_k}:\boldsymbol D\to\boldsymbol R^k$. For $T\subset[0,1]$, define in $\mathcal D$ the subclass $\mathcal D_f(T)$ of finite-dimensional sets $\pi_{t_1,\ldots,t_k}^{-1}(H)$, where $k$ is arbitrary, $t_i$ belong to $T$, and $H\in\mathcal R^k$.
\end{definition}
\begin{theorem}\label{12.5}  Consider projection mappings $\pi_{t_1,\ldots,t_k}:\boldsymbol D\to\boldsymbol R^k$. The following three statements hold.

(a) The projections $\pi_0$ and $\pi_1$ are continuous, and for $t\in(0,1)$, $\pi_t$ is continuous at $x$ if and only if $x$ is continuous at $t$.

(b) Each $\pi_{t_1,\ldots,t_k}$ is a measurable map.

(c) If $T$ contains 1 and is dense in $[0,1]$, then $\sigma\{\pi_t:t\in T\}=\sigma\{\mathcal D_f(T)\}=\mathcal D$  and $\mathcal D_f(T)$ is a separating class.
 \end{theorem}
 Proof. (a) Since each $\lambda\in\Lambda$ fixes 0 and 1, $\pi_0$ and $\pi_1$ are continuous: for $i=0,1$,
 \[d(x,y)\ge\inf_{\lambda\in\Lambda}\|x-y\lambda\|\ge |x(i)-y(i)|=|\pi_i(x)-\pi_i(y)|.\]
 Suppose that $0<t<1$. If $x$ is continuous at $t$, then by Lemma \ref{(12.14)}, $\pi_t$ is continuous at $x$. Suppose, on the other hand, that $x$ is discontinuous at $t$. If $\lambda_n\in\Lambda$ carries $t$ to $t-1/n$ and is linear on $[0,t]$ and $[t,1]$, and if $x_n(s)=x(\lambda_ns)$, then $d(x_n,x)\to0$ but $x_n(t)\nrightarrow x(t)$.
 
 (b) A mapping into $\boldsymbol R^k$ is measurable if each component mapping is. Therefore it suffices to show that 
 $\pi_t$ is measurable. Since $\pi_1$ is continuous, we may assume $t<1$. We use the pointwise convergence
 \[h_\epsilon(x):=\epsilon^{-1}\int_t^{t+\epsilon}x_sds\to\pi_t(x),\quad \epsilon\to0,\quad x\in \boldsymbol D.\]
If  $x^n\to x$ in the Skorokhod topology, then $x^n_t\to x_t$ for continuity points $t$ of $x$, and we conclude that  $x^n\to x$ almost surely. By the Bounded Convergence Theorem, the almost sure convergence $x^n\to x$ and the uniform boundedness of $(x^n)$ imply
$h_\epsilon(x^n)\to h_\epsilon(x)$.
Thus for each $\epsilon$, $h_\epsilon$ is continuous and therefore measurable, implying that  its limit $\pi_t$ is also measurable. 

(c) By right-continuity and the assumption that $T$ is dense, it follows that $\pi_0$ is measurable with respect to 
$\sigma\{\mathcal D_f(T)\}$. So we may as well assume that $0\in T$.

Suppose $s_0,\ldots,s_k$ are points in $T$ satisfying $0=s_0<\ldots<s_k=1$. For 
$\alpha=(\alpha_0,\ldots,\alpha_k)\in\boldsymbol R^{k+1}$ define $V\alpha\in \boldsymbol D$ by 
\[
(V\alpha)(t)=\left\{
\begin{array}{ll}
 \alpha_{j-1} &   \mbox{for }  t\in[s_{j-1},s_j),\ j=1,\ldots,k,\\
\alpha_k  &    \mbox{for }  t=1.
\end{array}
\right.
\]
Clearly, $V:\boldsymbol R^{k+1}\to \boldsymbol D$ is continuous implying that $\kappa=V\pi_{s_0,\ldots,s_k}$ is measurable $\sigma\{\mathcal D_f(T)\}/\mathcal D$.

Since $T$ is dense, for any $n$ we can choose $s_0^n,\ldots,s_k^n$ so that $\max_i(s_i^n-s_{i-1}^n)<n^{-1}$. Put $\kappa^n=V\pi_{s_0^n,\ldots,s_k^n}$. With this choice define a map $A_n:\boldsymbol D\to \boldsymbol D$ by $A_nx=x\kappa^n$. By Lemma \ref{L12.3}, $A_nx\to x$ for each $x$. We conclude that the identity map is measurable $\sigma\{\mathcal D_f(T)\}/\mathcal D$ and therefore $\mathcal D\subset\sigma\{\mathcal D_f(T)\}$. Finally, since $\mathcal D_f(T)$ is a $\pi$-system, it is a separating class.

%\begin{corollary}\label{p135} 
%Both $w'(\cdot,\delta)$ and $w''(\cdot,\delta)$ are $\mathcal D$-measurable.
%\end{corollary}
%Proof.

\begin{definition}
  Let $D_c$ be the set of count paths: nondecreasing functions $x\in \boldsymbol D$  with $x(t)\in\mathbb Z$ for each $t$, and $x(t)- x(t-)=1$ at points of discontinuity. 
\end{definition}
\begin{exercise}
Find $d(x,y)$ for $x,y\in D_c$ in terms of the jump points of these two count paths. How does a fundamental sequence $(x_n)$ in $D_c$ look for large $n$? 
Show that $D_c$ is closed in the Skorokhod topology. 
\end{exercise}
\begin{lemma}
Let $T_0=\{t_1,t_2,\ldots\}$  be a countable, dense set in $[0,1]$, and put
 $\pi(x)=(x(t_1),x(t_2),\ldots).$ 
 
 (a) The mapping $\pi:\boldsymbol D\to\boldsymbol R^\infty$ is $\mathcal D/\mathcal R^\infty$-measurable.  
 
 (b) If $x,x_n\in D_c$ are such that $\pi (x_n)\to \pi(x)$, then 
 $x_n\to x$ in the Skorokhod topology. 
 \end{lemma}
 Proof.  (a) In terms of notation of Section \ref{wcR}, 
 $$\pi^{-1}(\pi^{-1}_kH)=\pi^{-1}_{t_1,\ldots,t_k}H\in\mathcal D\qquad\mbox{ for }H\in \mathcal R^k,$$ 
 and the finite-dimensional sets $\pi^{-1}_kH$ generate $\mathcal R^\infty$.

(b) Convergence $\pi (x_n)\to \pi(x)$ implies $x_n(t_i)\to x(t_i)$, which in turn means that $x_n(t_i)=x(t_i)$ for $n>n_i$, for all $i=1,2,\ldots$. A function in $x\in D_c$ has only finitely many discontinuities, say $0<s_1<\ldots< s_k\le1$. For a given $\epsilon$ choose points $u_i$ and $v_i$ in $T_0$ in such a way that $u_i<s_i\le v_i<u_{i}+\epsilon$ and the intervals $[v_{i-1},u_i], i=1,\ldots,k$ are disjoint, with $v_0=0$. Then for $n$ exceeding some $n_0$, $x_n$ agrees with $x$ over each $[v_{i-1},u_i]$ and has a single jump in each $[u_i,v_i]$. If $\lambda_n\in\Lambda$ carries $s_i$ to the point in $[u_i,v_i]$ where $x_n$ has a jump and is defined elsewhere by linearity, then $\|\lambda_n-1\|\le\epsilon$ and $x_n(\lambda_nt)\equiv x(t)$ implying $d(x_n,x)\le\epsilon$ for $n>n_0$.

\begin{theorem}\label{12.6} 
  Let $T_0$ be a countable, dense set in $[0,1]$. If $P_n(D_c)=P(D_c)=1$ and $P_{n}\pi^{-1}_{t_1,\ldots,t_k}\Rightarrow P\pi^{-1}_{t_1,\ldots,t_k}$ for all $k$-tuples in $T_0$, then $P_n\Rightarrow P$.
 \end{theorem}
 Proof. The idea is, in effect, to embed $D_c$ in $\boldsymbol R^\infty$ and apply Theorem \ref{E2.4}. By hypothesis, $P_{n}\pi^{-1}\pi^{-1}_k\Rightarrow P\pi^{-1}\pi^{-1}_k$, but since in  $\boldsymbol R^\infty$ weak convergence is the same thing as weak convergence of finite-dimensional distributions, it follows that $P_{n}\pi^{-1}\Rightarrow P\pi^{-1}$ in  $\boldsymbol R^\infty$. For $A\subset \boldsymbol D$, define $A^*=\pi^{-1}(\pi A)^-$. If $A\in\mathcal D$, then
 \[\limsup_nP_n(A)\le \limsup_nP_n(A^*)=\limsup_nP_n(\pi^{-1}(\pi A)^-)\le P(\pi^{-1}(\pi A)^-)=P(A^*).\]
 Therefore, if $F\in\mathcal D$ is closed, then
 \[\limsup_nP_n(F)=\limsup_nP_n(F\cap D_c)\le P((F\cap D_c)^*)=P((F\cap D_c)^*\cap D_c).\]
 
 It remains to show that if $F\in\mathcal D$ is closed, then $(F\cap D_c)^*\cap D_c\subset F$.
 Take an   $x\in(F\cap D_c)^*\cap D_c$. Since $x\in(F\cap D_c)^*$, we have $\pi(x)\in(\pi(F\cap D_c))^-$ and there is a sequence $x_n\in F\cap D_c$ such that $\pi (x_n)\to \pi(x)$. Because $x\in D_c$, the previous lemma gives $x_n\to x$ in the Skorokhod topology. Since $x_n\in F$ and $F\in\mathcal D$ is closed, we conclude that $x\in F$.

\begin{corollary}\label{E12.3}
 Suppose for each $n$, $\xi_{n1},\ldots,\xi_{nn}$ are iid indicator r.v. with $\mathbb P(\xi_{ni}=1)=\alpha/n$. If $X^n_t=\sum_{i\le nt} \xi_{ni}$, then  $X^n\Rightarrow X$ in $\boldsymbol D$, where $X$ is the Poisson process with parameter $\alpha$.
\end{corollary}
Proof. The random process $X^n_t=\sum_{i\le nt} \xi_{ni}$ has independent increments. Its finite-dimensional distributions weakly converge to that of the Poisson process $X$ with  
$$\mathbb P(X_t-X_s=k)={\alpha^k(t-s)^k\over k!}e^{-\alpha(t-s)}\qquad\mbox{ for }0\le s<t\le1.$$

\begin{exercise}
 Suppose that $\xi$ is uniformly distributed over $[{1\over3},{2\over3}]$, and consider the random functions
 \[X_t=2\cdot1_{\{t\in[\xi,1]\}},\quad X_t^n=1_{\{t\in[\xi-n^{-1},1]\}}+1_{\{t\in[\xi+n^{-1},1]\}}.\]
 Show that $X^n\nRightarrow X$, even though $(X^n_{t_1},\ldots,X^n_{t_k})\Rightarrow (X_{t_1},\ldots,X_{t_k})$ for all $(t_1,\ldots,t_k)$. Why does Theorem \ref{12.6} not apply?
 
\end{exercise}
\begin{lemma}\label{13.0}
Let  $P$ be a probability measure on $(\boldsymbol D,\mathcal D)$. Define $T_P\subset[0,1]$ as the collection of $t$ such that  the projection $\pi_t$ is $P$-almost surely continuous. The set $T_P$ contains 0 and 1, and its complement in $[0,1]$ is at most countable. For $t\in(0,1)$, $t\in T_P$ is equivalent to $P\{x:x(t)\neq x(t-)\}=0$.
 \end{lemma}
 Proof.  Recall Theorem \ref{12.5} (a) and put $J_t=\{x:x(t)\neq x(t-)\}$ for a $t\in(0,1)$. We have to show that $P(J_t)>0$ is possible for at most countably many $t$. Let $J_t(\epsilon)=\{x:|x(t)- x(t-)|>\epsilon\}$. For fixed, positive $\epsilon$ and $\delta$, there can be at most finitely many $t$ for which $P(J_t(\epsilon))\ge \delta$. Indeed, if $P(J_{t_n}(\epsilon))\ge\delta$ for infinitely many distinct $t_n$, then
 \[P(J_{t_n}(\epsilon)\mbox{ i.o.})=P(\limsup_{n}J_{t_n}(\epsilon))\ge \limsup_{n\to\infty}P(J_{t_n}(\epsilon))\ge \delta,\]
 contradicting the fact that for a single $x\in \boldsymbol D$ the jumps can exceed $\epsilon$ at only finitely many points, see Lemma \ref{L12.1}. Thus  
 $P(J_t(\epsilon))>0$ is possible for at most countably many $t$. The desired result follows from 
  \[\{t:P(J_t)>0\}=\bigcup_n\{t:P(J_t(n^{-1}))>0\},\]
  which in turn is a consequence of $P(J_t(\epsilon))\uparrow P(J_t)$ as $\epsilon\downarrow 0$.

 \begin{theorem}\label{13.1}
Let $P_n, P$ be probability measures on $(\boldsymbol D,\mathcal D)$. If the sequence  $(P_n)$ is tight and $P_{n}\pi^{-1}_{t_1,\ldots,t_k}\Rightarrow P\pi^{-1}_{t_1,\ldots,t_k}$ holds whenever $t_1,\ldots,t_k$ lie in $T_P$, then $P_n\Rightarrow P$.
 
\end{theorem}
Proof. We will show that if a subsequence $(P_{n'})\subset (P_n)$ converges weakly to some $Q$, then $Q=P$. Indeed, if $t_1,\ldots,t_k$ lie in $T_Q$, then 
$\pi_{t_1,\ldots,t_k}$ is continuous on a set of $Q$-measure 1, and therefore, $P_{n'}\Rightarrow Q$ implies by the mapping theorem that $P_{n'}\pi^{-1}_{t_1,\ldots,t_k}\Rightarrow Q\pi^{-1}_{t_1,\ldots,t_k}$. On the other hand, if $t_1,\ldots,t_k$ lie in $T_P$, then $P_{n'}\pi^{-1}_{t_1,\ldots,t_k}\Rightarrow P\pi^{-1}_{t_1,\ldots,t_k}$ by the assumption. Therefore, if $t_1,\ldots,t_k$ lie in $T_Q\cap T_P$, then $Q\pi^{-1}_{t_1,\ldots,t_k}=P\pi^{-1}_{t_1,\ldots,t_k}$. It remains to see  that $\mathcal D_f(T_Q\cap T_P)$ is a separating class by applying Lemma \ref{13.0} and Theorem \ref{12.5}.

\subsection{Tightness criteria in $\boldsymbol D$}
 
\begin{theorem}\label{13.2}
Let $P_n$ be probability measures on $(\boldsymbol D,\mathcal D)$. The sequence  $(P_n)$ is tight if and only if the following two conditions hold:
\begin{align*}
(i)& \quad \lim_{a\to\infty}\limsup_{n\to\infty}P_n(x: \|x\|\ge a)=0,  \\
(ii)& \quad \lim_{\delta\to0}\limsup_{n\to\infty}P_n(x: w'_x(\delta)\ge\epsilon)=0,\mbox{ for each positive }\epsilon. %\label{sCp} 
\end{align*}
Condition (ii) is equivalent to
\begin{align*}
(ii')& \quad \forall\epsilon,  \eta>0;  \exists\delta, n_0>0: \quad  %\label{sCp} 
P_n(x: w'_x(\delta)\ge\epsilon)\le\eta,  \mbox{ for }n>n_0.
\end{align*}
\end{theorem}
Proof. This theorem is proven similarly to Theorem \ref{7.3} using Theorem \ref{12.3}. Equivalence of (ii) and (ii$'$) is due to monotonicity of $w'_x(\delta)$.

\begin{theorem}\label{13.2'}
Let $P_n$ be probability measures on $(\boldsymbol D,\mathcal D)$. The sequence  $(P_n)$ is tight if and only if the following two conditions hold:
\begin{align*}
(i)& \quad \lim_{a\to\infty}\limsup_{n\to\infty}P_n(x: \|x\|\ge a)=0,  \\
(ii)& \quad \forall\epsilon,  \eta>0;  \exists\delta, n_0>0: \quad  %\label{sCp} 
\left\{
\begin{array}{l}
P_n(x: w''_x(\delta)\ge\epsilon)\le\eta,\\
P_n(x:  |x(\delta)-x(0)|\ge\epsilon)\le\eta,\\
P_n(x:  |x(1-)-x(1-\delta)|\ge\epsilon)\le\eta,  
\end{array}
\right.\quad \mbox{ for }n>n_0.
\end{align*}
\end{theorem}
Proof. This theorem follows from Theorem \ref{13.2} with $(i)$ and $(ii')$ using Lemma \ref{(12.32)}. (Recall how  Theorem \ref{12.4} was obtained from Theorem \ref{12.3} using Lemma \ref{(12.32)}.)

\begin{lemma}\label{C13}
 Turn to Theorems \ref{13.2} and \ref{13.2'}. Under (ii) condition (i) is equivalent to the following weaker version:

(i$\,'$) for each $t$ in a set $T$ that is dense in $[0,1]$ and contains 1,
\[\lim_{a\to\infty}\limsup_{n\to\infty}P_n(x: |x(t)|\ge a)=0,\]
%(i$\,''$) the last relation holds for $t=0$ and, recalling notation $ j(x)=\sup_{0<t\le1}|x(t)-x(t-)|$,
%\[\lim_{a\to\infty}\limsup_{n\to\infty}P_n(x: |j(x)|\ge a)=0.\]
%
\end{lemma}
Proof. The implication (i) $\Rightarrow$ (i$'$) is trivial. Assume  (ii) of Theorem \ref{13.2} and  (i$'$). For a given $\delta\in(0,1)$ choose from $T$ points $0<s_1<\ldots<s_k=1$ such that $\max\{s_1,s_2-s_{1},\ldots,s_k-s_{k-1}\}<\delta$. By hypothesis  (i$'$), there exists an $a$ such that
\[\qquad \qquad \qquad \qquad P_n(x: \max_{j}|x(s_j)|\ge a)<\eta,\quad n>n_0\qquad \qquad \qquad \qquad (*)\]
For a given $x$, take a $\delta$-sparse set $(t_0,\ldots,t_v)$ such that  all $w_x[t_{i-1},t_i)<w'_x(\delta)+1$. 
Since each $[t_{i-1},t_i)$ contains an $s_j$, we have 
$$\|x\|\le \max_{j}|x(s_j)|+w'_x(\delta)+1.$$
Using (ii$'$) of Theorem \ref{13.2} and $(*)$, we get $P_n(x: \|x\|\ge a+2)<2\eta$ implying (i).

%\begin{corollary}\label{13.2''}
%A sequence $(P_n)$ of probability measures on $(\boldsymbol D,\mathcal D)$ is tight provided
%\begin{align*}
%(i)& \quad \lim_{a\to\infty}\limsup_{n\to\infty}P_n(x: \|x\|\ge a)=0,  \\
%(ii)& \quad \forall\epsilon,  \eta>0;  \exists\delta, n_0>0: \quad  %\label{sCp} 
%\left\{
%\begin{array}{l}
%P_n(x: w''_x(\delta)\ge\epsilon)\le\eta,\\
%P_n(x:  |x(\delta)-x(0)|\ge\epsilon)\le\eta,\\
%P_n(x:  |x(1)-x(1-\delta)|\ge\epsilon)\le\eta,  
%\end{array}
%\right.\quad \mbox{ for }n>n_0.
%\end{align*}
%\end{corollary}
%Proof. Apply Corollary \ref{12.4'}.

\subsection{A key condition on 3-dimensional distributions}
The following condition plays an important role. 
\begin{definition}\label{kc}
%For $\epsilon>0$ define $A_\epsilon\in\mathcal R^3$ by
%\[A_\epsilon=\{(z_1,z_2,z_3): |z_2-z_1|\ge\epsilon,|z_3-z_2|\ge\epsilon\}.\]
%Consider a system of consistent finite dimensional distributions  $(\mu_{t_1,\ldots,t_k})$. We will write 
%$\mu\in H_{\alpha,\beta,a}$ if there exist  $\alpha>1$, $\beta\ge0$, and a nondecreasing continuous $H:[0,1]\to[0,a]$ such that  for all $\epsilon>0$ and all $0\le t_1\le t_2\le t_3\le1$,
For a probability measure $P$ on $\boldsymbol D$, we will write $P\in H_{\alpha,\beta}$, if there exist  $\alpha>1$, $\beta\ge0$, and a nondecreasing continuous $H:[0,1]\to\boldsymbol R$ such that  for all $\epsilon>0$ and all $0\le t_1\le t_2\le t_3\le1$,
%$P\pi^{-1}_{t_1,\ldots,t_k}=\mu_{t_1,\ldots,t_k}$ and $\mu\in H_{\alpha,\beta}$. 
\[P\pi^{-1}_{t_1,t_2,t_3}\{(z_1,z_2,z_3): |z_2-z_1|\ge\epsilon,|z_3-z_2|\ge\epsilon\}\le \epsilon^{-2\beta}(H(t_3)-H(t_1))^{\alpha}.\]
For a random element $X$ on $\boldsymbol D$ with probability distribution $P$, this condition $P\in H_{\alpha,\beta}$ means that  for all $\epsilon>0$  and all $0\le r\le s\le t\le1$
\[\mathbb P\Big(|X_s-X_r|\ge\epsilon, |X_t-X_s|\ge\epsilon\Big)\le \epsilon^{-2\beta}(H(t)-H(r))^{\alpha}.\]

\end{definition}

\begin{lemma}\label{10.3}
Let a random element $X$ on $\boldsymbol D$ have a probability distribution $P\in H_{\alpha,\beta}$. Then there is a constant $K_{\alpha,\beta}$ depending only on $\alpha$ and $\beta$ such that
\[\mathbb P\Big(\sup_{r\le s\le t}\big(|X_s-X_r|\wedge |X_t-X_s|\big)\ge\epsilon\Big)\le {K_{\alpha,\beta}\over\epsilon^{2\beta}}(H(1)-H(0))^{\alpha}.\]
\end{lemma}
 Proof. 
The stated estimate is obtained  in four consecutive steps.

Step 1. Let $T_k=\{i/2^k,0\le i\le 2^k\}$ and
\begin{align*}
 A_k&=\max\big(|X_s-X_r|\wedge |X_t-X_s|\mbox{ over the adjacent triplets }   r\le s\le t \mbox{ in }  T_k\big),\\
 B_k&=\max\big(|X_s-X_r|\wedge |X_t-X_s|\mbox{ over }   r\le s\le t \mbox{ from }  T_k\big).
\end{align*}
We will show that $B_k\le 2(A_1+\ldots+A_k)$. To this end, for each $t\in T_k$ define a $t_n\in T_{k-1}$ by
\[t_n=
\left\{
\begin{array}{ll}
 t &  \mbox{if }  t\in T_{k-1}, \\
 t -2^{-k}&  \mbox{if }  t\notin T_{k-1}  \mbox{ and }  |X_t-X_{t-2^{-k}}|\le |X_t-X_{t+2^{-k}}|, \\
 t +2^{-k}&  \mbox{if }  t\notin T_{k-1}  \mbox{ and }  |X_t-X_{t-2^{-k}}|> |X_t-X_{t+2^{-k}}|,
 \end{array}
\right.
\]
so that $|X_t-X_{t_n}|\le A_k$. Then for any triplet  $r\le s\le t$ from $T_k$,
\begin{align*}
|X_s-X_r|&\le |X_s-X_{s_n}|+|X_{s_n}-X_{r_n}|+|X_{r}-X_{r_n}| \le |X_{s_n}-X_{r_n}|+2A_k,\\
|X_t-X_s|&\le |X_{t_n}-X_{s_n}|+2A_k.
\end{align*}
Since here $r_n\le s_n\le t_n$ lie in $T_{k-1}$, it follows that $|X_s-X_r|\wedge |X_t-X_s|\le B_{k-1}+2A_k$, and therefore, 
$$B_k\le B_{k-1}+2A_k\le 2(A_1+\ldots+A_k),\quad k\ge1.$$

Step 2. Consider a special case when $H(t)\equiv t$. Using the right continuity of the paths we get from step 1 that 
\[\sup_{r\le s\le t}\big(|X_s-X_r|\wedge |X_t-X_s|\big)\le 2\sum_{k=1}^\infty A_k.\]
This implies that for any $\theta\in(0,1)$,
\begin{align*}
\mathbb P\Big(\sup_{r\le s\le t}\big(|X_s-X_r|&\wedge |X_t-X_s|\big)\ge2\epsilon\Big)\le \mathbb P\Big(\sum_{k=1}^\infty A_k\ge\epsilon\Big)\\
&\le \mathbb P\Big(\sum_{k=1}^\infty A_k\ge\epsilon(1-\theta)\sum_{k=1}^\infty\theta^k\Big)\le \sum_{k=1}^\infty \mathbb P\Big(A_k\ge\epsilon(1-\theta)\theta^k\Big)\\
&\le \sum_{k=1}^\infty\sum_{i=1}^{2^k-1}\mathbb P\Big(|X_{i/2^k}-X_{(i-1)/2^k}|\wedge|X_{(i+1)/2^k}-X_{i/2^k}|\ge\epsilon(1-\theta)\theta^k\Big).
\end{align*}
Applying the key condition with $H(t)\equiv t$ we derive from the previous relation choosing a $\theta\in(0,1)$ satisfying $\theta^{2\beta}>2^{1-\alpha}$, that the stated estimate holds in the special case
\begin{align*}
\mathbb P\Big(\sup_{r\le s\le t}\big(|X_s-X_r|\wedge |X_t-X_s|\big)\ge2\epsilon\Big)
&\le \sum_{k=1}^\infty2^k{2^{(1-k)\alpha}\over (\epsilon(1-\theta)\theta^k)^{2\beta}}\\
&= {2^\alpha\over\epsilon^{2\beta}(1-\theta)^{2\beta}} \sum_{k=1}^\infty(\theta^{-2\beta}2^{1-\alpha})^k.
\end{align*}

Step 3. For a  strictly increasing $H(t)$ take $a$ so that $a^{2\beta}H(1)^\alpha=1$, and define a new process $Y_t$ by $Y_t=aX_{b(t)}$, where the time change $b(t)$ is such that $H(b(t))=tH(1)$. Since
\begin{align*}
\mathbb P\Big(|Y_s-Y_r|\ge\epsilon; |Y_t-Y_s|\ge\epsilon\Big)&=\mathbb P\Big(|X_{b(s)}-X_{b(r)}|\ge a^{-1}\epsilon; |X_{b(t)}-X_{b(s)}|\ge a^{-1}\epsilon\Big)\\
&\le \epsilon^{-2\beta}(t-r)^{\alpha},
\end{align*}
we can apply the result of the step 2 to the new process and prove the statement of the lemma under the assumption of step 3. 

Step 4. If  $H(t)$ is not strictly increasing, put $H_v(t)=H(t)+v t$ for an arbitrary small positive $v$. We have 
\[\mathbb P\Big(|X_s-X_r|\ge\epsilon; |X_t-X_s|\ge\epsilon\Big)\le \epsilon^{-2\beta}(H_v(t)-H_v(r))^{\alpha},\]
and according to step 3
\[\mathbb P\Big(\sup_{r\le s\le t}\big(|X_s-X_r|\wedge |X_t-X_s|\big)\ge\epsilon\Big)\le {K_{\alpha,\beta}\over\epsilon^{2\beta}}(H(1)+v-H(0))^{\alpha}.\]
It remains to let $v$ go to 0. Lemma \ref{10.3} is proven.

\begin{lemma}\label{p143}
If $P\in H_{\alpha,\beta}$, then given a positive $\epsilon$, 
\[P(x: w''_x(\delta)\ge\epsilon)\le {2K_{\alpha,\beta}\over\epsilon^{2\beta}}(H(1)-H(0))(w_H(2\delta))^{\alpha-1},\]
so that $P(x: w''_x(\delta)\ge\epsilon)\to0$ as $\delta\to0$.
\end{lemma}
 Proof. Take $t_i=i\delta$ for $0\le i\le \lfloor1/\delta\rfloor$ and $t_{\lceil1/\delta\rceil}=1$. If $|t-r|\le\delta$, then $r$ and $t$ lie in the same $[t_{i-1},t_{i+1}]$ for some $1\le i\le\lceil1/\delta\rceil-1$. According to Lemma \ref{10.3}, for $X$ with distribution $P$,
\[\mathbb P\Big(\sup_{t_{i-1}\le r\le s\le t\le t_{i+1}}\big(|X_s-X_r|\wedge |X_t-X_s|\big)\ge\epsilon\Big)\le {K_{\alpha,\beta}\over\epsilon^{2\beta}}(H(t_{i+1})-H(t_{i-1}))^{\alpha},\]
and due to monotonicity of $H$,
\begin{align*}
\mathbb P\Big(w''(X,\delta)\ge\epsilon\Big)&\le \sum_{i=1}^{\lceil1/\delta\rceil-1}{K_{\alpha,\beta}\over\epsilon^{2\beta}}(H(t_{i+1})-H(t_{i-1}))^{\alpha}\\
&\le {K_{\alpha,\beta}\over\epsilon^{2\beta}}\Big(\sup_{0\le t\le 1-2\delta}(H(t+2\delta)-H(t))^{\alpha-1}\Big)2(H(1)-H(0))\\
&={2K_{\alpha,\beta}\over\epsilon^{2\beta}}(H(1)-H(0))(w_H(2\delta))^{\alpha-1}. 
\end{align*}
It remains to recall that the modulus of continuity $w_H(2\delta)$ of the uniformly continuous function $H$ converges to 0 as $\delta\to0$.

\subsection{A criterion for existence}
\begin{theorem}\label{13.6}There exists in $\boldsymbol D$ a random element with finite dimensional distributions $\mu_{t_1,\ldots,t_k}$ provided the following three conditions:

(i) the finite dimensional distributions $\mu_{t_1,\ldots,t_k}$ are consistent, see Definition \ref{Kcon},

(ii) there exist  $\alpha>1$, $\beta\ge0$, and a nondecreasing continuous $H:[0,1]\to\boldsymbol R$ such that  for all $\epsilon>0$ and all $0\le t_1\le t_2\le t_3\le1$,
%$P\pi^{-1}_{t_1,\ldots,t_k}=\mu_{t_1,\ldots,t_k}$ and $\mu\in H_{\alpha,\beta}$. 
\[\mu_{t_1,t_2,t_3}\{(z_1,z_2,z_3): |z_2-z_1|\ge\epsilon,|z_3-z_2|\ge\epsilon\}\le \epsilon^{-2\beta}(H(t_3)-H(t_1))^{\alpha},\]

(iii) $\mu_{t,t+\delta}\{(z_1,z_2): |z_2-z_1|\ge\epsilon\}\to0$ as $\delta\downarrow0$ for each $t\in[0,1)$.
\end{theorem}
Proof. The main idea, as in the proof of Theorem \ref{7.1} (a), is to construct a sequence $(X^n)$ of random elements in $\boldsymbol D$ such that the corresponding sequence of distributions $(P_n)$ is tight and has  the desired limit finite dimensional distributions $\mu_{t_1,\ldots,t_k}$. 

Let vector $(X_{n,0},\ldots,X_{n,2^n})$ have distribution $\mu_{t_{0},\ldots,t_{2^n}}$, where $t_i\equiv t^n_i=i2^{-n}$, and define
\[
X^n_t=\left\{
\begin{array}{ll}
X_{n,i}  &  \mbox{for } t\in[i2^{-n},(i+1)2^{-n}),\quad i=0,\ldots,2^n-1,  \\
X_{n,2^n}  &   \mbox{for } t=1.
\end{array}
\right.
\]
The rest of the proof uses Theorem \ref{13.2'} and is divided in four steps.

Step 1. For all $\epsilon>0$  and $r,s,t\in T_n=\{t_{0},\ldots,t_{2^n}\}$ we have by (ii),
\[\mathbb P\Big(|X_s^n-X_r^n|\ge\epsilon, |X_t^n-X_s^n|\ge\epsilon\Big)\le \epsilon^{-2\beta}(H(t)-H(r))^{\alpha}.\]
It follows that in general for $0\le r \le s\le t\le 1$
\[\mathbb P\Big(|X_s^n-X_r^n|\ge\epsilon, |X_t^n-X_s^n|\ge\epsilon\Big)\le \epsilon^{-2\beta}(H(t)-H(r-2^{-n}))^{\alpha},\]
where $H(t)=H(0)$ for $t<0$. Slightly modifying Lemma \ref{p143} we obtain that  given a positive $\epsilon$, there is a constant $K_{\alpha,\beta,\epsilon}$
\[P_n(x: w''_x(\delta)\ge\epsilon)\le K_{\alpha,\beta,\epsilon}(H(1)-H(0))(w_H(3\delta))^{\alpha-1},\quad \mbox{for }\delta\le 2^{-n}.\]
This gives the first part of (ii) in Theorem \ref{13.2'}.

Step 2. If $2^{-k}\le\delta$, then 
\[\|X^n\|\le \max_{t\in T_k}|X^n_t|+w''(X^n,\delta).\]
Since the distributions of the first term on the right all coincide for $n\ge k$, it follows by step 1 that condition (i) in Theorem \ref{13.2'} is satisfied.

Step 3. To take care of the second and third parts of (ii) in Theorem \ref{13.2'}, we fix some $\delta_0\in(0,1/2)$, and temporarily assume that
for $\delta\in(0,\delta_0)$, 
\[\qquad \qquad \mu_{0,\delta}\{(z_1,z_2):z_1=z_2\}=1,\qquad \mu_{1-\delta,1}\{(z_1,z_2):z_1=z_2\}=1.\qquad \qquad \qquad (*)\]
In this special case, the second and third parts of (ii) in Theorem \ref{13.2'} hold and we conclude that 
the sequence of distributions $P_n$ of $X_n$ is tight. 

By Prokhorov's theorem, $(X^n)$ has a subsequence  weakly converging in distribution  to a random element $X$ of $\boldsymbol D$ with some distribution $P$. We want to show that $P\pi^{-1}_{t_1,\ldots,t_k}=\mu_{t_1,\ldots,t_k}$. Because of the consistency hypothesis, this holds for dyadic rational $t_i\in \cup_nT_n$. The general case is obtained using the following facts:
\begin{align*}
P_n\pi^{-1}_{t_1,\ldots,t_k}&\Rightarrow P\pi^{-1}_{t_1,\ldots,t_k},\\
P_n\pi^{-1}_{t_1,\ldots,t_k}&=\mu_{t_{n1},\ldots,t_{nk}},\quad \mbox{ for some }t_{ni}\in T_n,\mbox{ provided }k\le 2^n,\\
\mu_{t_{n1},\ldots,t_{nk}}&\Rightarrow \mu_{t_1,\ldots,t_k}.
\end{align*}
The last fact is a consequence of (iii). Indeed, by Kolmogorov's extension theorem, there exists a stochastic process $Z$ with vectors $(Z_{t_{1}},\ldots,Z_{t_{k}})$ having distributions $\mu_{t_1,\ldots,t_k}$. Then by (iii), $Z_{t+\delta}\stackrel{\rm P}{\to}Z_t$ as $\delta\downarrow0$. Using Exercise \ref{inP} we derive 
$(Z_{t_{n1}},\ldots,Z_{t_{nk}})\stackrel{\rm P}{\to}(Z_{t_{1}},\ldots,Z_{t_{k}})$ implying $\mu_{t_{n1},\ldots,t_{nk}}\Rightarrow \mu_{t_1,\ldots,t_k}$.

Step 4. It remains to remove the restriction $(*)$. To this end take
\[
\lambda t=\left\{
\begin{array}{ll}
 0 &   \mbox{for } t\le\delta_0, \\
 {t-\delta_0\over1-2\delta_0} &    \mbox{for } \delta_0<t<1-\delta_0,   \\
 1 &    \mbox{for } t\ge1-\delta_0.
\end{array}
\right.
\]
Define $\nu_{s_1,\ldots,s_k}$ as $\mu_{t_1,\ldots,t_k}$ for $s_i=\lambda t_i$. Then the $\nu_{s_1,\ldots,s_k}$ satisfy the conditions of the theorem with a new $H$, as well as $(*)$, so that there is  a random element $Y$ of $\boldsymbol D$ with these finite-dimensional distributions. Finally,
setting $X_t=Y_{\delta_0+t(1-2\delta_0)}$ we get a process $X$ with the required finite dimensional distributions $P\pi^{-1}_{t_1,\ldots,t_k}=\mu_{t_1,\ldots,t_k}$.

\begin{example}\label{E13.1}
 Construction of a Levy process. Let $\nu_t$ be a measure on the line for which $\nu_t(\boldsymbol R)=H(t)$ is nondecreasing and continuous, $t\in[0,1]$. Suppose for 
 $s\le t$, $\nu_s(A)\le\nu_t(A)$ for all $A\in\mathcal R$ so that $\nu_t-\nu_s$ is a measure with total mass $H(t)-H(s)$. Then there is an infinitely divisible distribution having mean 0, variance $H(t)-H(s)$, and characteristic function 
 \[\phi_{s,t}(u)=\exp \int_{-\infty}^\infty{e^{iuz}-1-iuz\over z^2}(\nu_t-\nu_s)(dz).\]
\end{example}
We can use Theorem \ref{13.6} to construct a random element $X$ of $\boldsymbol D$ with $X_0=0$, for which the increments are independent and
$$\mathbb E(e^{iu_1X_{t_1}}e^{iu_2(X_{t_2}-X_{t_1})}\cdots e^{iu_k(X_{t_k}-X_{t_{k-1}})})=\phi_{0,t_1}(u_1)\phi_{t_1,t_2}(u_2)\cdots \phi_{t_{k-1},t_k}(u_k).$$
Indeed, since $\phi_{r,t}(u)=\phi_{r,s}(u)\phi_{s,t}(u)$ for  $r\le s\le t$, the implied finite-dimensional distributions $\mu_{t_1,\ldots,t_k}$ are consistent. Further, by Chebyshev's inequality and independence, condition (ii) of Theorem \ref{13.6} is valid with $\alpha=\beta=2$:
\begin{align*}
\mu_{t_1,t_2,t_3}\{(z_1,z_3,z_3): |z_2-z_1|\ge\epsilon,|z_3-z_2|\ge\epsilon\}&\le {H(t_2)-H(t_1)\over\epsilon^{2}}\cdot {H(t_3)-H(t_2)\over\epsilon^{2}}\\&\le \epsilon^{-4}(H(t_3)-H(t_1))^{2}.
\end{align*}
Another application of Chebyshev's inequality gives
\begin{align*}
\mu_{t,t+\delta}\{(z_1,z_2):&|z_2-z_1|\ge\epsilon\}\le {H(t+\delta)-H(t)\over\epsilon^{2}}\to0,\quad \delta\downarrow0.
\end{align*}

\section{Weak convergence on $\boldsymbol D$}
Recall that the subset $T_P\subset[0,1]$, introduced in Lemma \ref{13.0}, is the collection of $t$ such that  the projection $\pi_t$ is $P$-almost surely continuous. 
\subsection{Criteria  for weak convergence in $\boldsymbol D$}

\begin{lemma}\label{13.3'}
Let $P$ be a probability measure on $(\boldsymbol D,\mathcal D)$ and $\epsilon>0$. By right continuity of the paths we have $ \lim_{\delta\to0}P(x:  |x(\delta)-x(0)|\ge\epsilon)=0$.
\end{lemma}
Proof. Put $A_\delta=\{x:  |x(\delta)-x(0)|\ge\epsilon\}$. Let $\delta_n\to0$. It suffices to show that $P(A_{\delta_n})\to0$. To see this observe that right continuity of the paths entails 
\[\bigcap_{n\ge1}\bigcup_{k\ge n}A_{\delta_k}=\emptyset,\]
and therefore, $P(A_{\delta_n})\le P(\cup_{k\ge n}A_{\delta_k})\to0$ as $n\to\infty$.

\begin{theorem}\label{13.3}
Let $P_n, P$ be probability measures on $(\boldsymbol D,\mathcal D)$. Suppose $P_{n}\pi^{-1}_{t_1,\ldots,t_k}\Rightarrow P\pi^{-1}_{t_1,\ldots,t_k}$ holds whenever $t_1,\ldots,t_k$ lie in $T_P$. If for every positive $\epsilon$ 
\begin{align*}
(i)& \quad \lim_{\delta\to0}P(x: |x(1)-x(1-\delta)|\ge \epsilon)=0,  \\
(ii)& \quad \lim_{\delta\to0}\limsup_{n\to\infty}P_n(x: w''_x(\delta)\ge\epsilon)=0,%\label{sCp} 
\end{align*}
then $P_n\Rightarrow P$.
\end{theorem}
Proof. This result should be compared with Theorem \ref{7.5} dealing with the space $\boldsymbol C$. 

Recall Theorem \ref{13.1}. We prove tightness by checking conditions (i$'$) in Lemma \ref{C13} and (ii) in Theorem \ref{13.2'}. For each $t\in T_P$, the weakly convergent sequence $P_n\pi^{-1}_t$ is tight which implies (i$'$) with $T_P$ in the role of $T$.

As to  (ii) in Theorem \ref{13.2'} we have to verify only the second and third parts. By hypothesis, $P_{n}\pi^{-1}_{0,\delta}\Rightarrow P\pi^{-1}_{0,\delta}$ so that  for $\delta\in T_P$,
\[\limsup_{n\to\infty}P_n(x:  |x(\delta)-x(0)|\ge\epsilon)\le P(x:  |x(\delta)-x(0)|\ge\epsilon),\]
and the second part follows from Lemma \ref{13.3'}. 

Turning to the third part of (ii) in Theorem \ref{13.2'}, the symmetric to the last argument brings for $1-\delta\in T_P$,
\[\limsup_{n\to\infty}P_n(x:  |x(1)-x(1-\delta)|\ge\epsilon)\le P(x:  |x(1)-x(1-\delta)|\ge\epsilon).\]
Now, suppose that\footnote{Here we use an argument suggested by Timo Hirscher.}
\[|x(1-)-x(1-\delta)|\ge\epsilon.\]
Since
\[|x(1-)-x(1-\delta)|\wedge|x(1)-x(1-)|\le w''_x(\delta),\]
we have either $w''_x(\delta)\ge\epsilon$
or $|x(1)-x(1-)|\le w''_x(\delta)$. Moreover, in the latter case, it is either 
$|x(1)-x(1-\delta)|\ge\epsilon/2$ or $w''_x(\delta)\ge\epsilon/2$ or both.
The last two observations yield
\begin{align*}
 \limsup_{n\to\infty}&P_n(x:  |x(1-)-x(1-\delta)|\ge\epsilon)
 \\
 &\le P(x:  |x(1)-x(1-\delta)|\ge\epsilon/2)
 +\limsup_{n\to\infty}P_n(x: w''_x(\delta)\ge\epsilon/2),
\end{align*}
and the third part readily follows from conditions (i) and (ii).

\begin{theorem}  \label{13.5} For $X^n\Rightarrow X$ on $\boldsymbol D$ it suffices that 

(i) $(X^n_{t_1},\ldots,X^n_{t_k})\Rightarrow (X_{t_1},\ldots,X_{t_k})$ for points $t_i\in T_P$, where $P$ is the probability distribution of $X$, 

(ii) $X_1-X_{1-\delta}\Rightarrow 0$ as $\delta\to0$,

(iii) there exist  $\alpha>1$, $\beta>0 $, and a nondecreasing continuous function $H:[0,1]\to\boldsymbol R$ such that 
\[\mathbb E\Big(|X^n_s-X^n_r|^\beta|X^n_t-X^n_s|^\beta\Big)\le (H(t)-H(r))^{\alpha}\quad \mbox{ for }0\le r\le s\le t\le1.\]
 \end{theorem}
 Proof. By Theorem  \ref{13.3}, it is enough to show that 
 $$\lim_{\delta\to0}\limsup_{n\to\infty}\mathbb P\Big(w''(X^n,\delta)\ge\epsilon\Big)=0.$$
 This follows from Lemma \ref{p143} as (iii) implies that $X_n$ has a distribution $P_n\in H_{\alpha,\beta}$.

\subsection{Functional CLT  on $\boldsymbol D$}
The identity map $c:\boldsymbol C\to \boldsymbol D$ is continuous and therefore measurable $\mathcal C/\mathcal D$. If $\mathbb W$ is a Wiener measure on $(\boldsymbol C,\mathcal C)$, then $\mathbb Wc^{-1}$ is a Wiener measure on $(\boldsymbol D,\mathcal D)$. We denote this new measure by   $\mathbb W$ rather than $\mathbb Wc^{-1}$. Clearly, $\mathbb W(\boldsymbol C)=1$. Let also denote by $W$ a random element of $\boldsymbol D$ with distribution $\mathbb W$.

\begin{theorem}\label{14.1}
 Let $\xi_1,\xi_2,\ldots$ be iid r.v. defined on $(\Omega,\mathcal F,\mathbb P)$. If $\xi_i$ have zero mean and variance $\sigma^2$ and  
 $X^n_t={\xi_1+\ldots+\xi_{\lfloor nt\rfloor}\over\sigma\sqrt n}$, then $X^n\Rightarrow W$.
\end{theorem}
Proof. We apply Theorem \ref{13.5}. Following the proof of Theorem \ref{8.2'} (a) one
gets the convergence of the fdd (i) even for the $X^n$ as they are
defined here. Condition (ii) follows from the fact that the Wiener process $W_t$ has no jumps. We finish the proof by showing that (iii) holds with $\alpha=\beta=2$ and $H(t)=2t$. Indeed, 
$$ \mathbb E\big(|X^n_s-X^n_r|^2|X^n_t-X^n_s|^2\big)=0\quad \mbox{ for }0\le t-r< n^{-1},$$
as either $X^n_s=X^n_r$ or $X^n_t=X^n_s$.
On the other hand, for $t-r\ge n^{-1}$, by independence,
\begin{align*}
 \mathbb E\big(|X^n_s-X^n_r|^2|X^n_t-X^n_s|^2\big)&={\lfloor ns\rfloor -\lfloor nr\rfloor\over n}\cdot{\lfloor nt\rfloor -\lfloor ns\rfloor\over n}\\
& \le \Big({\lfloor nt\rfloor -\lfloor nr\rfloor\over n}\Big)^2 \le (2(t-r))^2.
\end{align*}

\begin{example}\label{p148}
 Define $(\xi_n)$ on $([0,1],\mathcal B_{[0,1]},\lambda)$ using the Rademacher functions $\xi_n(\omega)=2w_n-1$ in terms of the dyadic (binary) representation $\omega=\omega_1\omega_2\ldots$. Then $(\xi_n)$ is a sequence of independent coin tossing outcomes with values $\pm1$. Theorem \ref{14.1} holds with $\sigma=1$:
  $\Big({\xi_1+\ldots+\xi_{\lfloor nt\rfloor}\over\sqrt n}\Big)_{t\in[0,1]}\Rightarrow W.$

\end{example}

\begin{lemma}\label{M21} Consider a probability space $(\Omega,\mathcal F,\mathbb P)$ and let $\mathbb P_0$ be a probability measure absolutely continuous with respect to $\mathbb P$. Let $\mathcal F_0\subset \mathcal F$ be an algebra of events such that for some $A_n\in\sigma(\mathcal F_0)$
\[\mathbb P(A_n|E)\to\alpha,\quad \mbox{ for all }E\in\mathcal F_0 \mbox{ with }\mathbb P(E)>0.\]
 Then $\mathbb P_0(A_n)\to\alpha$.
\end{lemma}
Proof. We have $\mathbb P_0(A)=\int_Ag_0(\omega)\mathbb P(d\omega)$, where $g_0=d\mathbb P_0/d\mathbb P$. It suffices to prove that 
\[\ \qquad\qquad\qquad \qquad\int_{A_n}g(\omega)\mathbb P(d\omega)\to \alpha\int_{\Omega}g(\omega)\mathbb P(d\omega)\qquad\qquad \qquad\qquad (*)\]
if $g$ is $\mathcal F$-measurable and $\mathbb P$-integrable. We prove $(*)$ in three steps.

Step 1. Write $\mathcal F_1=\sigma(\mathcal F_0)$ and denote by $\mathcal F_2$ the class of events $E$ for which 
\[\mathbb P(A_n\cap E)\to\alpha\mathbb P(E).\]
We show that $\mathcal F_1\subset \mathcal F_2$. To be able to apply Theorem \ref{Dyn} we have to show that $\mathcal F_2$ is a $\lambda$-system. Indeed, suppose for a sequence of disjoint sets $E_i$ we have
\[\mathbb P(A_n\cap E_i)\to\alpha\mathbb P(E_i).\]
Let $E=\cup_iE_i$, then by Lemma \ref{Mtest},
\[\mathbb P(A_n\cap E)=\sum_i\mathbb P(A_n\cap E_i)\to\alpha\sum_i\mathbb P(E_i)=\alpha\mathbb P(E).\]

Step 2. Show that  $(*)$ holds for $\mathcal F_1$-measurable functions $g$. Indeed, due to step 1,  relation $(*)$ holds if $g$ is the indicator of an $\mathcal F_1$-set. and hence if it is a simple $\mathcal F_1$-measurable function. If $g$ is $\mathcal F_1$-measurable and $\mathbb P$-integrable function, choose simple $\mathcal F_1$-measurable functions $g_k$ that satisfy $|g_k|\le|g|$ and $g_k\to g$. Now 
\[\Big|\int_{A_n}g(\omega)\mathbb P(d\omega)-\alpha\int_{\Omega}g(\omega)\mathbb P(d\omega)\Big|\le\Big|\int_{A_n}g_k(\omega)\mathbb P(d\omega)-\alpha\int_{\Omega}g_k(\omega)\mathbb P(d\omega)\Big|+(1+\alpha)\mathbb E|g-g_k|.\]
Let first $n\to\infty$ and then $k\to\infty$ and apply the dominated convergence theorem.

Step 3. Finally, take  $g$ to be a $\mathcal F$-measurable and $\mathbb P$-integrable. We use conditional expectation $g_1=\mathbb E(g|\mathcal F_1)$
\[\int_{A_n}g(\omega)\mathbb P(d\omega)=\mathbb E(g1_{\{A_n\}})=\mathbb E(g_11_{\{A_n\}})\to \alpha\mathbb E(g_1)=\alpha\int_{\Omega}g(\omega)\mathbb P(d\omega).\]

\begin{theorem}\label{14.2}
 Let $\xi_1,\xi_2,\ldots$ be iid r.v. defined on $(\Omega,\mathcal F,\mathbb P)$ having zero mean and variance $\sigma^2$. Put 
 $X^n_t={\xi_1+\ldots+\xi_{\lfloor nt\rfloor}\over\sigma\sqrt n}$. If  $\mathbb P_0$ is a probability measure absolutely continuous with respect to $\mathbb P$, then $X^n\Rightarrow W$ with respect to $\mathbb P_0$.
\end{theorem}
Proof. Step 1. Choose $k_n$ such that $k_n\to\infty$ and $k=o(n)$ as $n\to\infty$ and put
\begin{align*}
 \bar X^n_t&={1\over\sigma\sqrt n}\sum_ {i=k_n}^{\lfloor nt\rfloor}\xi_i,\\
 Y_n&={1\over\sigma\sqrt n}\max_{1\le k< k_n}|\xi_1+\ldots+\xi_k|.
\end{align*}
By Kolmogorov's inequality, for any $a>0$, 
$$\mathbb P(Y_n\ge a)\to0, \ \ n\to\infty,$$ 
and therefore
\begin{align*}
 d(X^n,\bar X^n)\le\|X^n-\bar X^n\|=Y_n\Rightarrow0 \quad\mbox{with respect to }\mathbb P.
\end{align*}
Applying Theorem \ref{14.1} and Corollary \ref{3.1} we conclude that $\bar X^n\Rightarrow W$ with respect to $\mathbb P$.

Step 2:  show using Lemma \ref{M21}, that $\bar X^n\Rightarrow W$ with respect to $\mathbb P_0$.
If $A\in\mathcal D$ is a $\mathbb W$-continuity set, then $\mathbb P(A_n)\to\alpha$ for $A_n=\{\bar X^n\in A\}$ and $\alpha=\mathbb W(A)$. Let 
$\mathcal F_0$ be the algebra of the cylinder sets $\{(\xi_1,\ldots,\xi_k)\in H\}$. If $E\in\mathcal F_0$, then $A_n$ are independent of $E$ for large $n$ and by  Lemma \ref{M21}, $\mathbb P_0(\bar X^n\in A)\to\mathbb W(A)$.

Step 3. Since $1_{\{Y_n\ge a\}}\to0$ almost surely with respect to $\mathbb P$, the dominated convergence theorem gives
$$\mathbb P_0(Y_n\ge a)=\int g_0(\omega)1_{\{Y_n\ge a\}} \mathbb P(d\omega)\to0.$$
Arguing as in step 1 we conclude that $d(X^n,\bar X^n)\Rightarrow0$ with respect to $\mathbb P_0$.

Step 4. Applying once again Corollary \ref{3.1} we conclude that $X^n\Rightarrow W$ with respect to $\mathbb P_0$.

\begin{example}\label{p148.} Define $\xi_n$ on $([0,1],\mathcal B_{[0,1]},\lambda_p)$ with
$$\lambda_p(du)=au^{a-1}du,\quad a=-\log_2(1-p),\quad p\in(0,1)$$
again, as in Example \ref{p148}, using the Rademacher functions. If $p={1\over2}$, then $\lambda_p=\lambda$ and we are back to Example \ref{p148}.
With $p\ne{1\over2}$, this corresponds to dependent $p$-coin tossings  with 
$$\int_0^{2^{-n}}\lambda_p(du)=(1-p)^n$$ 
being the probability of having $n$ failures in the first $n$ tossings, and
$$\int_{1-2^{-n}}^1\lambda_p(du)=1-(1-p)^{-\log_2(1-2^{-n})}$$ 
being the probability of having $n$ successes in  the first $n$ tossings. By Theorem \ref{14.2}, even in this case
  $\Big({\xi_1+\ldots+\xi_{\lfloor nt\rfloor}\over\sqrt n}\Big)_{t\in[0,1]}\Rightarrow W.$
\end{example}

\subsection{Empirical distribution functions}
\begin{definition}
 Let $\xi_1(\omega),\ldots,\xi_n(\omega)$ be iid with a distribution function $F$ over $[0,1]$. 
The corresponding empirical process is defined by $Y^n_t=\sqrt n(F_n(t)-F(t))$, where 
$$F_n(t)=n^{-1}(1_{\{\xi_1\le t\}}+\ldots+1_{\{\xi_n\le t\}})$$ 
is the empirical distribution function.
\end{definition}
\begin{figure}
\centering
\includegraphics[width=12cm,height=6cm]{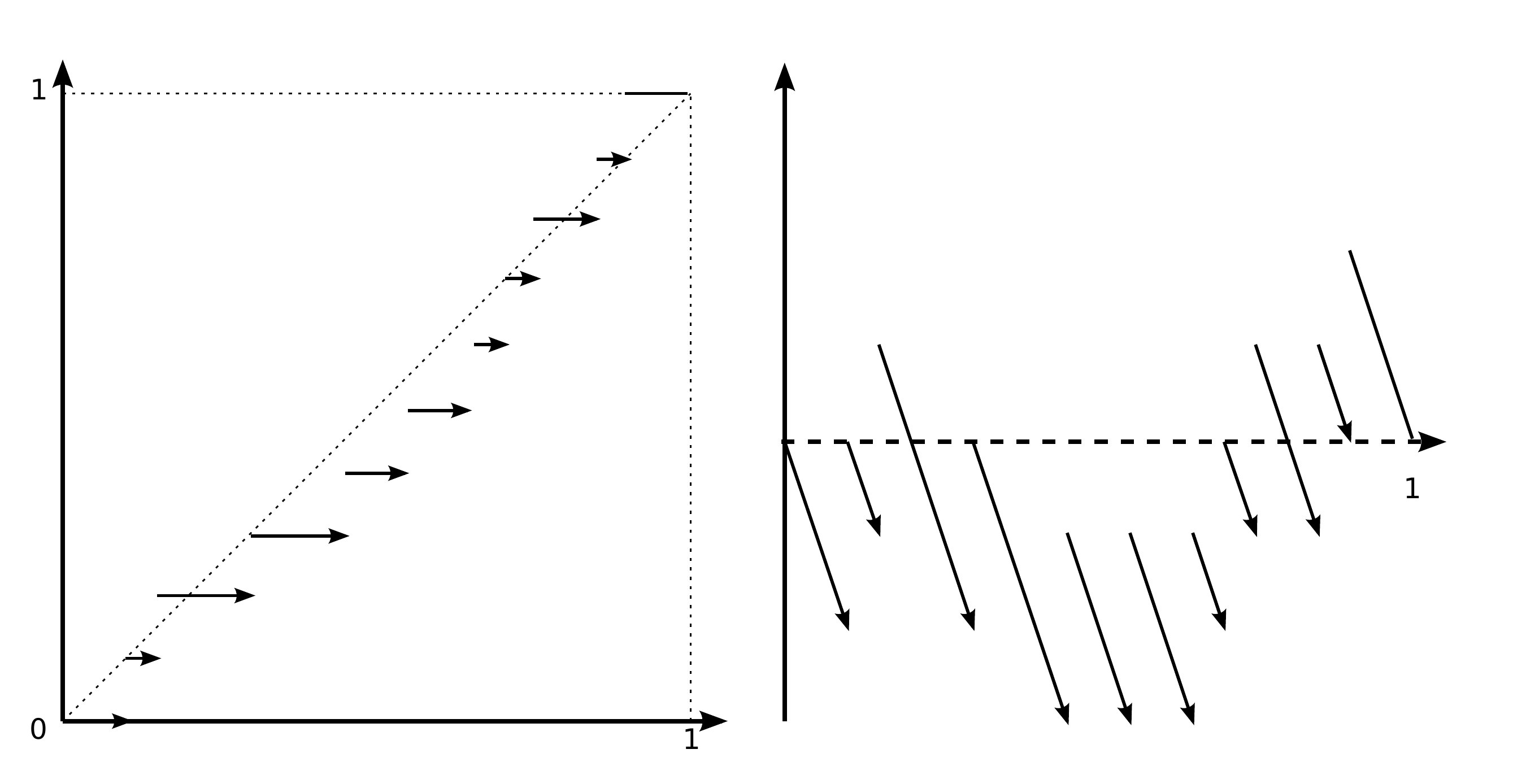}
\caption{An example of the empirical distribution function with $n=10$ for the uniform distribution (left panel) and the corresponding empirical process (right panel).}
\label{empp}
\end{figure}
\begin{lemma}
Let  $(Z_1^n,\ldots,Z_r^n)$ have a multinomial distribution Mn$(n,p_1,\ldots,p_r)$. Then the normalized vector 
$\Big({Z_1^n-np_1\over\sqrt n},\ldots,{Z_r^n-np_r\over\sqrt n}\Big)$ converges in distribution to a multivariate normal distribution with zero means and a covariance matrix
\[
\bf V=\left(
\begin{array}{ccccc}
p_1(1-p_1)  &  -p_1p_2 & -p_1p_3&\ldots&-p_1p_r  \\
 -p_2p_1 &p_2(1-p_2)   & -p_2p_3 &\ldots&-p_2p_r  \\
 -p_3p_1 &-p_3p_2   & p_3(1-p_3) &\ldots&-p_3p_r  \\
\ldots&\ldots&\ldots&\ldots&\ldots \\
-p_rp_1 &-p_rp_2   & -p_rp_3 &\ldots&p_r(1-p_r)
\end{array}
\right).
\]
\end{lemma}
Proof. To apply the continuity property of the multivariate characteristic functions consider
\[\mathbb E\exp\Big(i\theta_1{Z_1^n-np_1\over\sqrt n}+\ldots+i\theta_r{Z_r^n-np_r\over\sqrt n}\Big)=\Big(\sum_{j=1}^r p_j e^{i \tilde\theta _j/\sqrt n}\Big)^n,\]
where $\tilde\theta _j=\theta _j-(\theta_1p_1+\ldots+\theta_rp_r)$. Similarly to the classical case we have
\[\Big(\sum_{j=1}^r p_j e^{i \tilde\theta _j/\sqrt n}\Big)^n=\Big(1-{1\over 2n}\sum_{j=1}^r p_j \tilde\theta _j^2+o(n^{-1})\Big)^n\to 
e^{-{1\over 2}\sum_{j=1}^r p_j \tilde\theta _j^2}=e^{-{1\over 2}(\sum_{j=1}^r p_j \theta _j^2-(\sum_{j=1}^r p_j \theta _j)^2)}.\]
It remains to see that the right hand side equals $e^{-{1\over2}\boldsymbol \theta \mathbf V\boldsymbol \theta ^{\rm t}}$ which follows from the representation
\[
{\bf V}=\left(
\begin{array}{ccc}
p_1  & &0\\
&\ddots& \\
0 &&p_r
\end{array}
\right)-
\left(
\begin{array}{c}
p_1 \\

\vdots \\
p_r
\end{array}
\right)
\Big(p_1,\ldots,p_r\Big).
\]

\begin{theorem}\label{14.3}
 If $\xi_1,\xi_2,\ldots$ are iid  $[0,1]$-valued r.v. with a distribution function $F$, then the empirical process weakly converges $Y^n\Rightarrow Y$ to a random element $(Y_t)_{t\in[0,1]}=(W^\circ_{F(t)})_{t\in[0,1]}$, where $W^\circ$ is the standard Brownian bridge. The limit $Y$ is a Gaussian process specified by $\mathbb E(Y_t)=0$ and $\mathbb E(Y_sY_t)=F(s)(1-F(t))$ for $s\le t$.

\end{theorem}
Proof. We start with the uniform case, $F(t)\equiv t$ for $t\in[0,1]$, by showing $Y^n\Rightarrow W^\circ$, where $W^\circ$ is the Brownian bridge with  $\mathbb E(W^\circ_sW^\circ_t)=s(1-t)$ for $s\le t$. Let 
$$U^n_t=nF_n(t)=1_{\{\xi_1\le t\}}+\ldots+1_{\{\xi_n\le t\}}$$ 
be the number of $\xi_1,\ldots,\xi_n$ falling inside $[0,t]$. Since the increments of $U^n_t$ are described by multinomial joint distributions, by the previous lemma, the fdd of  $Y^n_t={U^n_t-nt\over \sqrt n}$ converge to those of $W^\circ$. Indeed, for $t_1<t_2<\ldots$ and $i<j$,
\[\mathbb E(W^\circ_{t_i}-W^\circ_{t_{i-1}})(W^\circ_{t_j}-W^\circ_{t_{j-1}})=-(t_i-t_{i-1})(t_j-t_{j-1})=-p_ip_j.\]

By Theorem \ref{13.5} it suffices to prove for $t_1\le t\le t_2$ that 
\[\mathbb E\Big((Y^n_t-Y^n_{t_1})^2(Y^n_{t_2}-Y^n_t)^2\Big)\le (t-t_1)(t_2-t)\le (t_2-t_1)^2.\]
In terms of $\alpha_i=1_{\{\xi_i\in(t_1,t]\}}+t_1-t$ and $\beta_i=1_{\{\xi_i\in(t,t_2]\}}+t-t_2$ the first inequality is equivalent to
\[\mathbb E\Big(\big(\sum_{i=1}^n\alpha_i\big)^2\big(\sum_{i=1}^n\beta_i\big)^2\Big)\le n^2(t-t_1)(t_2-t).\]
As we show next, this follows from $\mathbb E(\alpha_i)=\mathbb E(\beta_i)=0$, independence $(\alpha_i,\beta_i)\independent(\alpha_j,\beta_j)$ for $i\ne j$, and the following formulas for the second order moments. Let us write $p_1=t-t_1$, and $p_2=t_2-t$. Since 
\[
\alpha_i=\left\{
\begin{array}{cl}
1-p_1  & \mbox{w.p.  }  p_1, \\
-p_1  & \mbox{w.p.  } 1- p_1,
\end{array}
\right.\ 
\beta_i=\left\{
\begin{array}{cl}
1-p_2  & \mbox{w.p. }  p_2, \\
-p_2  & \mbox{w.p.  } 1- p_2,
\end{array}
\right.\]
and 
\[\alpha_i\beta_i=\left\{
\begin{array}{cl}
-(1-p_1)p_2  & \mbox{w.p.  }  p_1, \\
-p_1(1-p_2)  & \mbox{w.p.  } p_2,\\
p_1p_2  & \mbox{w.p.  } 1-p_1-p_2,
\end{array}
\right.
\]
we have
\begin{align*}
 &\mathbb E(\alpha_i^2)=p_1(1-p_1),\quad \mathbb E(\beta_i^2)=p_2(1-p_2),\quad \mathbb E(\alpha_i\beta_i)=-p_1p_2,\\
 &\mathbb E(\alpha_i^2\beta_i^2)=p_1(1-p_1)^2p_2^2+p_2p_1^2(1-p_2)^2+(1-p_1-p_2)p_1^2p_2^2
 =p_1p_2(p_1+p_2-3p_1p_2)
\end{align*}
and 
\begin{align*}
\mathbb E\Big(\big(\sum_{i=1}^n\alpha_i\big)^2\big(\sum_{i=1}^n\beta_i\big)^2\Big)&=n\mathbb E(\alpha_i^2\beta_i^2)
+n(n-1)\mathbb E(\alpha_i^2)\mathbb E(\beta_i^2)+2n(n-1)(\mathbb E(\alpha_i\beta_i))^2\\
&\le n^2  p_1 p_2 (p_1 +p_2-3p_1 p_2+1-p_1-p_2+3p_1 p_2) \\
&= n^2  p_1 p_2%=n^2p_1p_2(1-p_1-p_2+3p_1p_2)\le n^2p_1p_2
=n^2(t-t_1)(t_2-t).
\end{align*}

This proves the theorem for the uniform case. For a general continuous and strictly increasing $F(t)$ we use the transformation $\eta_i=F(\xi_i)$ into uniformly distributed r.v. If $G_n(t)=G_n(t,\omega)$ is the empirical distribution function of $(\eta_1(\omega),\ldots,\eta_n(\omega))$ and $Z^n_t=\sqrt n(G_n(t)-t)$, then $Z^n\Rightarrow W^\circ$. 

Observe that
\[G_n(F(t))={1_{\{\eta_1\le F(t)\}}+\ldots+1_{\{\eta_n\le F(t)\}}\over n}={1_{\{F(\xi_1)\le F(t)\}}+\ldots+1_{\{F(\xi_n)\le F(t)\}}\over n}=F_n(t).\]
Define $\psi:\boldsymbol D\to \boldsymbol D$ by $(\psi x)(t)=x(F(t))$. If $x_n\to x$ in the Skorokhod topology and $x\in \boldsymbol C$, then the convergence is uniform, so that $\psi x_n\to \psi x$ uniformly and hence in the Skorokhod topology. By the mapping theorem $\psi(Z^n)\Rightarrow \psi(W^\circ)$. Therefore, 
$$Y^n=(Y^n_t)_{t\in[0,1]}=(Z^n_{F(t)})_{t\in[0,1]}=\psi(Z^n)\Rightarrow \psi(W^\circ)=(W^\circ_{F(t)})_{t\in[0,1]}=Y.$$

\begin{figure}
\centering
\includegraphics[width=6cm,height=6cm]{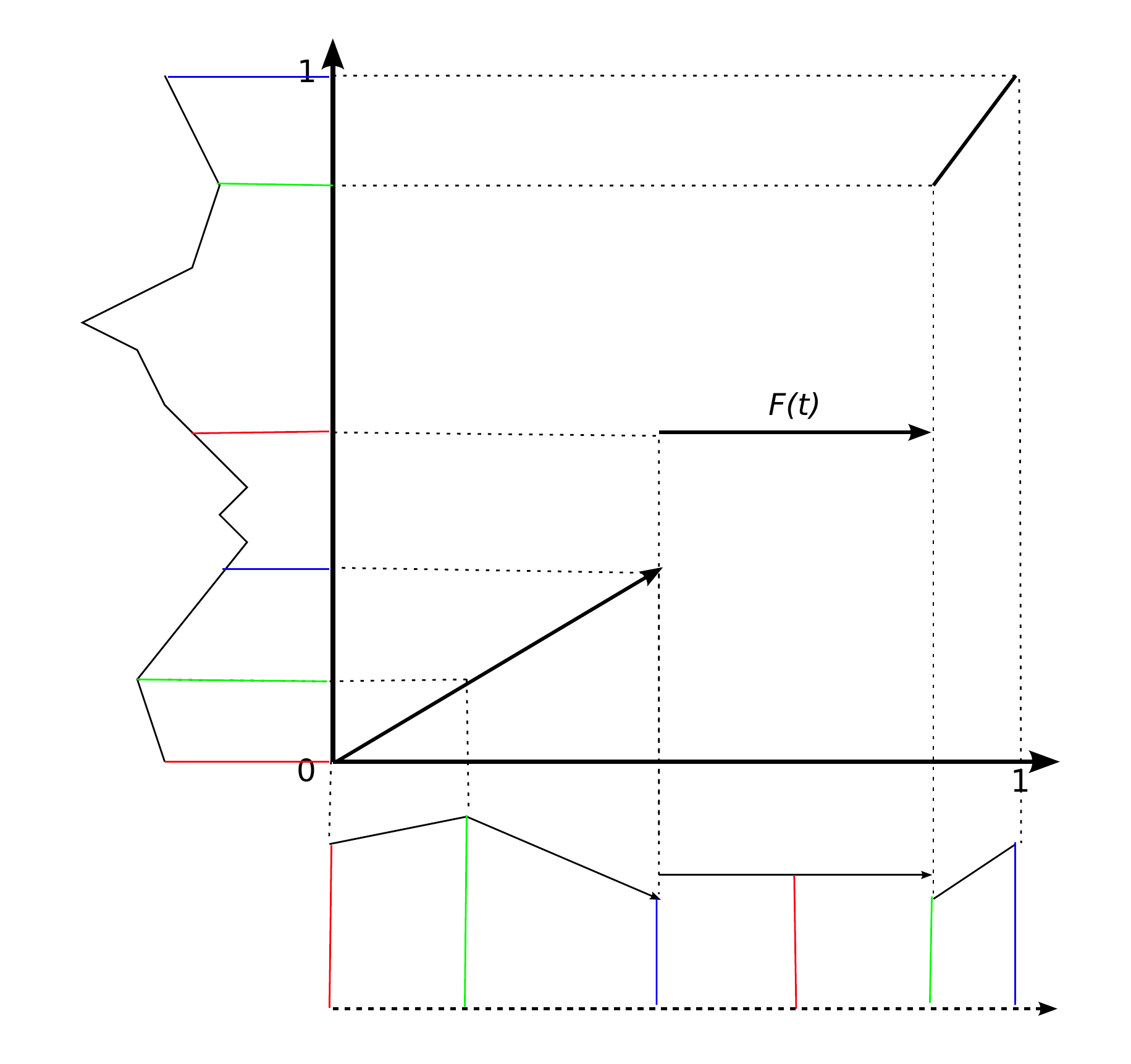}
\caption{The time transformation $\psi:\boldsymbol D\to \boldsymbol D$. Along the $y$-axis an original path $x(t)$ is depicted, along the $x$-axis the time transformed path $(\psi x)(t)=x(F(t))$ is given. The jumps of $F$ translate into path jumps, the constant parts of  $F$ translate into the horizontal pieces of the transformed path.}
\label{titr}
\end{figure}

Finally, for $F(t)$ with jumps and constant parts (see Figure \ref{titr}) the previous argument works provided there exists an iid sequence $\eta_1,\eta_2,\ldots$  of uniformly distributed r.v. as well as iid $\xi'_1,\xi'_2,\ldots$ with distribution function $F$, such that 
$$\{\eta_i\le F(t)\}\equiv\{\xi_i'\le t\},\quad t\in [0,1],i\ge1.$$
This is achieved by starting with uniform $\eta_1,\eta_2,\ldots$ on possibly another probability space and putting $\xi_i'=\phi(\eta_i)$, where $\phi(u)=\inf\{t: u\le F(t)\}$ is the quantile function satisfying
\[\{\phi(u)\le t\}=\{u\le F(t)\}.\]
\begin{example} Kolmogorov-Smirnov test. Let $F$ be continuous. By the mapping theorem we obtain
 \[\sqrt n\sup_t|F_n(t)-F(t)|=\sup_t|Y^n_t|\Rightarrow \sup_t|W^\circ_{F(t)}|=\sup_t|W^\circ_{t}|,\]
 where the limit distribution is given by Theorem \ref{(9.39)}.
\end{example}

\section{The space $\boldsymbol D_\infty=\boldsymbol D[0,\infty)$}\label{secDi}
\subsection{Two metrics on  $\boldsymbol D_\infty$}
To extend the Skorokhod theory to the space $\boldsymbol D_\infty=\boldsymbol D[0,\infty)$ of the cadlag functions on $[0,\infty)$, consider for each $t>0$ the space $\boldsymbol D_t=\boldsymbol D[0,t]$ of the same cadlag functions restricted on $[0,t]$. All definitions for $\boldsymbol D=\boldsymbol D[0,1]$ have obvious analogues for $\boldsymbol D_t$: for example we denote by $d_t(x,y)$ the analogue of $d_1(x,y):=d(x,y)$. Denote $\|x\|_t=\sup_{u\in[0,t]}|x(u)|$.

\begin{example}
 One might try to define Skorokhod convergence $x_n\to x$ on $\boldsymbol D_\infty$ by requiring that $d_t(x_n,x)\to 0$ for each finite $t>0$. This does not work: if $x_n(t)=1_{\{t\in[0,1-n^{-1})\}}$, the natural limit would be $x(t)=1_{\{t\in[0,1)\}}$ but $d_1(x_n,x)=1$ for all $n$. The problem here is that $x$ is discontinuous at $t=1$, and the definition must accommodate discontinuities.
\end{example}
\begin{lemma}\label{L16.1}
 Let $0<u<t<\infty$. If $d_t(x_n,x)\to 0$ and $x$ is continuous at $u$, then $d_u(x_n,x)\to 0$.
\end{lemma}
Proof. By hypothesis, there are time transforms $\lambda_n\in \Lambda_t$ such that $\|\lambda_n-1\|_t\to0$ and $\|x_n-x\lambda_n\|_t\to0$ as $n\to\infty$. Given $\epsilon$, choose $\delta$ so that $|v-u|\le2\delta$ implies $|x(v)-x(u)|\le\epsilon/2$. Now choose $n_0$ so that, if $n\ge n_0$ and $v\le t$, then $|\lambda_n v-v|\le\delta$ and $|x_n(v)-(x\lambda_n)(v)|\le\epsilon/2$. Then, if $n\ge n_0$ and $|v-u|\le\delta$, we have 
\[|\lambda_n v-u|\le|\lambda_n v-v|+|v-u|\le2\delta\]
and hence 
\[|x_n(v)-x(u)|\le|x_n(v)-(x\lambda_n)(v)|+|(x\lambda_n)(v)-x(u)|\le\epsilon.\]
Thus
\[\sup_{|v-u|\le\delta}|x(v)-x(u)|\le\epsilon, \quad \sup_{|v-u|\le\delta}|x_n(v)-x(u)|\le\epsilon \quad \mbox{ for }n\ge n_0.\]

Let
\[u_n=
\left\{
\begin{array}{ll}
u-n^{-1}  &   \mbox{if }  \lambda_n u>u, \\
u  &   \mbox{if }  \lambda_n u=u, \\
\lambda_n^{-1}(u-n^{-1})  &   \mbox{if }  \lambda_n u<u,
\end{array}
\right.
\]
so that $u_n\le u$. Since
\[|u_n -u|\le|\lambda_n^{-1}(u-n^{-1}) -(u-n^{-1})|+n^{-1},\]
 we have $u_n\to u$, and since 
 \[|\lambda_nu_n-u|\le|\lambda_nu_n-u_n|+|u_n -u|,\]
 we also have $\lambda_nu_n\to u$. 
 \begin{figure}
\centering
\includegraphics[width=8cm,height=8cm]{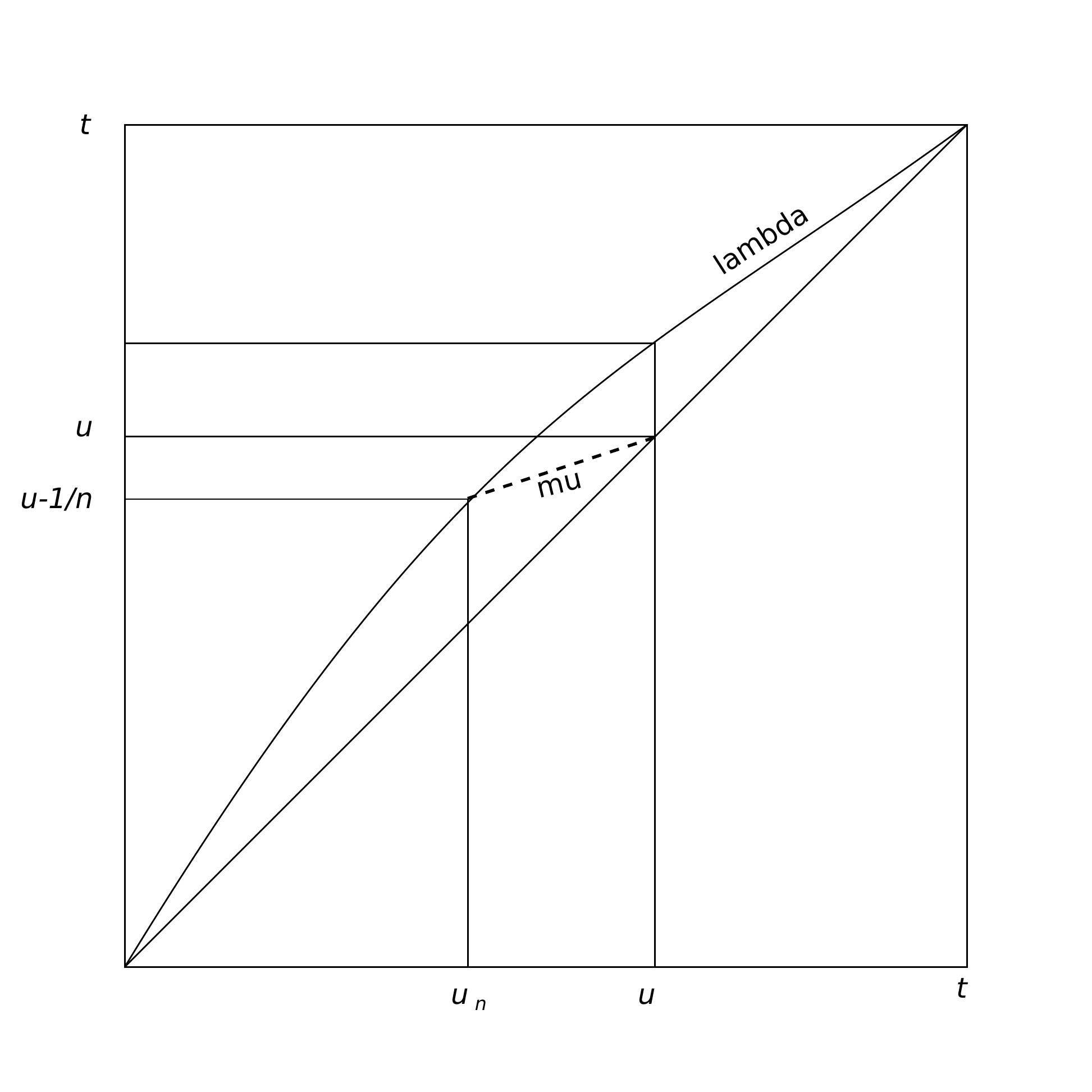}
\caption{A detail of the proof of Lemma \ref{L16.1} and Theorem \ref{16.1}.}
\label{lem}
\end{figure}

 Define $\mu_n\in\Lambda_u$ so that $\mu_nv=\lambda_nv$ for $v\in[0,u_n]$ and interpolate linearly on $(u_n,u]$ aiming at the diagonal point $\mu_nu=u$, see Figure \ref{lem}. By linearity, $|\mu_nv-v|\le |\lambda_nu_n-u_n|$ for $v\in[u_n,u]$ and we have $\|\mu_n-1\|_u\to0$. 
 
 It remains to show that $\|x_n-x\mu_n\|_u\to0$. To do this we choose $n_1$ so that 
 \[u_n>u-\delta\quad \mbox{ and }\quad \lambda_nu_n>u-\delta\quad \mbox{for }n\ge n_1.\]
 If $v\le u_n$, then 
 \[|x_n(v)-(x\mu_n)(v)|=|x_n(v)-(x\lambda_n)(v)|\le \|x_n-x\lambda_n\|_t.\]
 On the other hand, if $v\in[u_n,u]$ and $n\ge n_1$, then  $v\in[u-\delta,u]$ and $\mu_nv\in[u-\delta,u]$ implying for $n\ge n_1\vee n_0$
 \[|x_n(v)-(x\mu_n)(v)|\le|x_n(v)-x(u)|+|x(u)-(x\mu_n)(v)|\le 2\epsilon.\]
 The proof is finished.
\begin{definition}
For any natural $i$, define a map $\psi_i:\boldsymbol D_\infty\to \boldsymbol D_i$ by 
$$(\psi_ix)(t)=x(t)1_{\{t\le i-1\}}+(i-t)x(t)1_{\{i-1<t\le i\}}$$ 
making the transformed function $(\psi_ix)(t)$ continuous at $t=i$. 
\end{definition}

\begin{definition}\label{p168}
 Two topologically equivalent metrics $d_\infty(x,y)$ and $d_\infty^\circ(x,y)$ are defined on $\boldsymbol D_\infty$ in terms of $d(x,y)$ and $d^\circ(x,y)$ by
 \[d_\infty(x,y)=\sum_{i=1}^\infty{1\wedge d_i(\psi_ix,\psi_iy)\over 2^i},\qquad d_\infty^\circ(x,y)=\sum_{i=1}^\infty{1\wedge d_i^\circ(\psi_ix,\psi_iy)\over 2^i}.\]

% \[
%g_ix(t)=g_i(t)x(t),\qquad g_iy(t)=g_i(t)y(t),\qquad g_i(t)=\left\{
%\begin{array}{ll}
%1  &  \mbox{ if }  t\le i-1,\\
%i-t  &  \mbox{ if }  i-1< t<i,\\
%0  &  \mbox{ if }  t\ge i.   
%\end{array}
%\right.
%\]
\end{definition}

The metric properties of $d_\infty(x,y)$ and $d_\infty^\circ(x,y)$ follow from those of $d_i(x,y)$ and $d_i^\circ(x,y)$. In particular, if $d_\infty(x,y)=0$, then $d_i(\psi_ix,\psi_iy)=0$ and $\psi_ix=\psi_iy$ for all $i$, and this implies $x=y$.

\begin{lemma}
The map $\psi_i:\boldsymbol D_\infty\to \boldsymbol D_i$ is continuous.
\end{lemma}
Proof. It follows from the fact that $d_\infty(x_n,x)\to0$ implies $d_i(\psi_ix_n,\psi_ix)\to0$.

\subsection{Characterization of Skorokhod convergence on  $\boldsymbol D_\infty$}
Let $\Lambda_\infty$ be the set of continuous, strictly increasing maps $\lambda:[0,\infty)\to[0,\infty)$ such that $\lambda 0=0$ and $\lambda t\to\infty$ as $t\to\infty$. Denote $\|x\|_\infty=\sup_{u\in[0,\infty)}|x(u)|$. 
\begin{exercise}
 Let $\lambda\in\Lambda_t$ where $t\in(0,\infty]$. Show that the inverse transformation $\lambda^{-1}\in\Lambda_t$ is such that $\|\lambda^{-1}-1\|_t=\|\lambda-1\|_t$.
\end{exercise}
\begin{example}
 Consider the sequence $x_n(t)=1_{\{t\ge n\}}$ of elements of $\boldsymbol D_\infty$. Its natural limit is $x\equiv0$ as $\|x_n-x\|_t=0$ for all $n>t>0$. However, $\|x_n\lambda_n-x\|_\infty=\|x_n\lambda_n\|_\infty=1$ for any choice of $\lambda_n\in\Lambda_\infty$.
\end{example}
\begin{theorem}\label{16.1}
 Convergence  $d_\infty(x_n,x)\to0$ takes place if and only if there is a sequence $\lambda_n\in\Lambda_\infty$ such that
 
  $$\|\lambda_n-1\|_\infty\to0\mbox{ and }\|x_n\lambda_n-x\|_i\to0\mbox{ for each }i.$$
\end{theorem}
Proof. Necessity. Suppose $d_\infty(x_n,x)\to0$. Then $d_i(\psi_ix_n,\psi_ix)\to0$ and there exist $\lambda_n^{(i)}\in\Lambda_i$ such that 
\[\epsilon_n^{(i)}=\|\lambda_n^{(i)}-1\|_i\vee\|(\psi_ix_n)\lambda_n^{(i)}-\psi_ix\|_i\to0,\quad n\to\infty\quad \mbox{for each }i.\]
Choose $n_i>1$ such that $n\ge n_i$ implies $\epsilon_n^{(i)}<i^{-1}$. Arrange that $n_i<n_{i+1}$, and let 
$$i_1=\ldots=i_{n_1-1}=1,\quad i_{n_1}=\ldots=i_{n_2-1}=2,\quad i_{n_2}=\ldots=i_{n_3-1}=3, \quad\ldots$$
so that  $i_n\to\infty$. Define  $\lambda_n\in\Lambda_\infty$ by
\[\lambda_nt=
\left\{
\begin{array}{ll}
 \lambda_n^{(i_n)}t &  \mbox{if }  t\le i_n, \\
 t &    \mbox{if }  t> i_n. 
 \end{array}
\right.
\]
Then  
$$\|\lambda_n-1\|_\infty=\|\lambda_n^{(i_n)}-1\|_{i_n}\le \epsilon_n^{(i_n)}<i^{-1}_n\to0.$$ 
Now fix $i$. If $n$ is large enough, then $i< i_n-1$ and
\begin{align*}
\|x_n\lambda_n-x\|_i&=\|(\psi_{i_n}x_n)\lambda_n-\psi_{i_n}x\|_i\\
&\le \|(\psi_{i_n}x_n)\lambda_n-\psi_{i_n}x\|_{i_n}=\|(\psi_{i_n}x_n)\lambda_n^{(i_n)}-\psi_{i_n}x\|_{i_n}\le \epsilon_n^{(i_n)}<i^{-1}_n\to0. 
\end{align*}

Sufficiency. Suppose that  there is a sequence $\lambda_n\in\Lambda_\infty$ such that, firstly, $\|\lambda_n-1\|_\infty\to0$, and secondly, 
 $\|x_n\lambda_n-x\|_i\to0$   for each $i$. Observe that for some $C_i$,
 $$\|x\|_i\le C_i \mbox{ and }\|x_n\|_i\le C_i\mbox{ for all }(n, i) .$$
 Indeed, by the first assumption, for large $n$ we have $\lambda_n(2i)>i$ implying 
$\|x_n\|_i\le \|x_n\lambda_n\|_{2i}$, where by the second assumption, $\|x_n\lambda_n\|_{2i}\to\|x\|_{2i}$.

 Fix an $i$. It is enough to show that $d_i(\psi_{i}x_n,\psi_{i}x)\to0$. As in the proof of Lemma \ref{L16.1} define
\[u_n=
\left\{
\begin{array}{ll}
i-n^{-1}  &   \mbox{if }  \lambda_n i<i, \\
i  &   \mbox{if }  \lambda_n i=i, \\
\lambda_n^{-1}(i-n^{-1})  &   \mbox{if }  \lambda_n i>i,
\end{array}
\right.
\]
and $\mu_n\in\Lambda_i$ so that $\mu_nv=\lambda_nv$ for $v\in[0,u_n]$ interpolating linearly on $(u_n,i]$ towards $\mu_ni=i$. As before, $\|\mu_n-1\|_i\to0$ and it suffices to check that 
$$\|(\psi_{i}x_n)\mu_n-\psi_{i}x\|_i\to0,\quad n\to\infty.$$

To see that the last relation holds suppose $j:=\lambda_n^{-1}(i-1)\le i-1$ (the other case $j>i-1$ is treated similarly) and observe that
 \begin{align*}
\|(\psi_{i}x_n)\lambda_n-\psi_{i}x\|_i&= \Big(\|x_n\lambda_n-x\|_{j}\Big)\vee\Big(\sup_{j< t\le i-1}|(i-\lambda_nt)x_n(\lambda_nt)-x(t)|\Big)\\
&\vee\Big(\sup_{i-1< t\le i}|(i-\lambda_nt)x_n(\lambda_nt)-(i-t)x(t)|\Big)\\
&\le %\Big(\|x_n\lambda_n-x\|_{i}\Big)\vee
\Big(\|x_n\lambda_n-x\|_{i}+C_i\sup_{j< t\le i-1}|(i-1-\lambda_nt)|\Big)\\
&\vee\Big(\sup_{i-1< t\le i}|i-\lambda_nt|\cdot|x_n(\lambda_nt)-x(t)|+\sup_{i-1< t\le i}|(\lambda_nt-t)x(t)|\Big)\\
&\le \|x_n\lambda_n-x\|_{i}+C_i\|\lambda_n-1\|_i 
\to0.
\end{align*}
It follows that for $t\le u_n$,
 \begin{align*}
|(\psi_{i}x_n)(\mu_nt)-(\psi_{i}x)(t)|\le \|(\psi_{i}x_n)\lambda_n-\psi_{i}x\|_i \to0.
\end{align*}
Turning to the case $u_n<t\le i$, given an $\epsilon\in(0,1)$ choose $n_0$ such that for $n>n_0$, $u_n$ and $\mu_nu_n$ both lie in $[i-\epsilon,i]$. Then
\begin{align*}
|(\psi_{i}x_n)(\mu_nt)-(\psi_{i}x)(t)|\le \sup_{u_n< t\le i}|(i-\mu_nt)x_n(\mu_nt)-(i-t)x(t)|\le 2C_i\epsilon.
\end{align*}

\begin{theorem}\label{16.2}
 Convergence  $d_\infty(x_n,x)\to0$ takes place if and only if $d_t(x_n,x)\to0$ for each continuity point $t$ of $x$.
\end{theorem}
Proof. Necessity. If $d_\infty(x_n,x)\to0$, then $d_i(\psi_{i}x_n,\psi_{i}x)\to0$ for each $i$. Given a continuity point $t$ of $x$, take an integer $i$ for which $t<i-1$. According to Lemma \ref{L16.1}, $d_t(x_n,x)=d_t(\psi_{i}x_n,\psi_{i}x)\to0$.

Sufficiency. Choose continuity points $t_i$ of $x$ in such a way that $t_i\uparrow\infty$ as $i\to\infty$.  By hypothesis,
\[d_{t_i}(x_n,x)\to0,\quad n\to\infty,\quad i\ge1.\]
Choose $\lambda_n^{(i)}\in\Lambda_{t_i}$ so that 
\[\epsilon_n^{(i)}=\|\lambda_n^{(i)}-1\|_{t_i}\vee\|x_n\lambda_n^{(i)}-x\|_{t_i}\to0,\quad n\to\infty\quad \mbox{for each }i.\]
Using the argument from the first part of the proof of Theorem \ref{16.1}, define integers $i_n$ in such a way that $i_n\to\infty$ and $\epsilon_n^{(i_n)}<i^{-1}_n$. Put
\[\lambda_nt=
\left\{
\begin{array}{ll}
 \lambda_n^{(i_n)}t &  \mbox{if }  t\le t_{i_n}, \\
 t &    \mbox{if }  t> t_{i_n},
 \end{array}
\right.
\]
so that $\lambda_n\in\Lambda_\infty$. We have  $\|\lambda_n-1\|_\infty\le i^{-1}_n$, and for any given $i$, if  $n$ is sufficiently large so that $i<t_{ i_n}$, then 
\[\|x_n\lambda_n-x\|_i=\|x_n\lambda_n^{(i_n)}-x\|_{i}\le\|x_n\lambda_n^{(i_n)}-x\|_{t_{i_n}}\le \epsilon_n^{(i_n)}<i^{-1}_n\to0.\]
Applying Theorem \ref{16.1} we get $d_\infty(x_n,x)\to0$. 

\begin{exercise}
 Show that the mapping  $h(x)=\sup_{t\ge0}x(t)$ is not continuous on $\boldsymbol D_\infty$. 
 \end{exercise}

\subsection{Separability and completeness of  $\boldsymbol D_\infty$}
\begin{lemma}\label{M6}
 Suppose $(\boldsymbol S_i,\rho_i)$ are metric spaces and consider $\boldsymbol S=\boldsymbol S_1\times \boldsymbol S_2\times\ldots$ together with the metric of coordinate-wise convergence
 \[\rho(x,y)=\sum_{i=1}^\infty{1\wedge \rho_i(x_i,y_i)\over 2^i}.\]
If each $\boldsymbol S_i$ is separable, then $\boldsymbol S$ is separable. If each $\boldsymbol S_i$ is complete, then $\boldsymbol S$ is complete. 
%(iii) If $A_i$ is compact in $S_i$, then $A=A_1\times A_2\times\ldots$ is compact in $S$.
\end{lemma}
Proof. Separability. For each $i$, let $B_i$ be a countable dense subset in $\boldsymbol S_i$ and $x_i^\circ \in \boldsymbol S_i$ be a fixed point. We will show that the countable set $B=B(x_{1}^\circ,x_{2}^\circ,\ldots)$ defined by
\[B=\{x\in \boldsymbol S: x=(x_1,\ldots,x_k,x_{k+1}^\circ,x_{k+2}^\circ,\ldots), x_1 \in B_1, \ldots x_k \in B_k, k\in\boldsymbol N\}\]
is dense in $\boldsymbol S$.
Given an $\epsilon$ and a point $y\in \boldsymbol S$, choose $k$ so that $\sum_{i>k}2^{-i}<\epsilon$ and then choose points $x_i \in B_i$ so that $\rho_i(x_i,y_i)<\epsilon$. With this choice the corresponding point $x\in B$ satisfies $\rho(x,y)<2\epsilon$.

Completeness. Suppose that $x^n=(x^n_1,x^n_2,\ldots)$ are points of $\boldsymbol S$ forming a fundamental sequence. Then each sequence $(x^n_i)$ is fundamental in $\boldsymbol S_i$ and hence $\rho_i(x^n_i,x_i)\to0$ for some $x_i \in \boldsymbol S_i$.
By the M-test, Lemma \ref{Mtest}, $\rho(x^n,x)\to0$, where $x=(x_1,x_2,\ldots)$.

\begin{definition}\label{deep}
 Consider the product space $\overline {\boldsymbol D}=\boldsymbol D_1\times \boldsymbol D_2\times\ldots$ with the coordinate-wise convergence metric (cf Definition \ref{p168})
  \[\rho(\bar x,\bar y)=\sum_{i=1}^\infty{1\wedge d_i^\circ(\bar x_i,\bar y_i)\over 2^i},\qquad \bar x=(\bar x_1,\bar x_2,\ldots),\ \bar y=(\bar y_1,\bar y_2,\ldots)\in \overline {\boldsymbol D}.\]
 Put  $\psi x=(\psi_1 x,\psi_2 x,\ldots)$ for $ x\in \boldsymbol D_\infty$. Then $\psi x\in\overline {\boldsymbol D}$ and $d^\circ _\infty(x,y)=\rho(\psi x,\psi y)$ so that $\psi$ is an isometry of $(\boldsymbol D_\infty,d^\circ _\infty)$ into $(\overline {\boldsymbol D},\rho)$.
\end{definition}
\begin{lemma}\label{L16.2}
 The image $\overline {\boldsymbol D}_\infty:=\psi \boldsymbol D_\infty$ is closed in $\overline {\boldsymbol D}$.
\end{lemma}
Proof. Suppose that $x_n\in \boldsymbol D_\infty$, $\bar x=(\bar x_1,\bar x_2,\ldots)\in \overline {\boldsymbol D}$, and $\rho(\psi x_n,\bar x)\to0$, then $d_i(\psi_i x_n,\bar x_i)\to0$ for each $i$. We must find an $x\in \boldsymbol D_\infty$ such that $\bar x=\psi x$. 

The sequence of functions $\bar x_i\in \boldsymbol D_i$, $i=1,2,\ldots$ has at most countably many points of discontinuity. Therefore, there is a dense set $T\in[0,\infty)$ such that for every $i\ge t\in T$, the function $\bar x_i(\cdot)$ is continuous at $t$. Since $d_i(\psi_i x_n,\bar x_i)\to0$, we have $\psi_i x_n(t)\to\bar x_i(t)$ for all $t\in T\cap[0,i]$. This means that for every $t\in T$ there exists the limit $x(t)=\lim_nx_n(t)$, since $\psi_i x_n(t)= x_n(t)$ for $i>t+1$.

Now $\psi_i x(t)= \bar x_i(t)$ on $T\cap[0,i]$. Hence $x(t)= \bar x_i(t)$ on $T\cap[0,i-1]$, so that $x$ can be extended to a cadlag function on each $[0,i-1]$ and then  to a cadlag function on $[0,\infty)$. We conclude, using right continuity, that $\psi_i x(t)= \bar x_i(t)$ for all $t\in[0,i]$.

\begin{theorem}\label{16.3} The metric space $(\boldsymbol D_\infty,d_\infty^\circ)$ is separable and complete. 
 \end{theorem}
 Proof. According Lemma \ref{M6} the space  $\overline {\boldsymbol D}$ is separable and complete, so are  the closed subspace $\overline {\boldsymbol D}_\infty$ and its isometric copy $\boldsymbol D_\infty$.
 
%\subsection{Characterization of relative compactness}

\subsection{Weak convergence on  $\boldsymbol D_\infty$}
\begin{definition}
For any natural $i$ and any $s\ge i$, define a map $\psi_{s,i}:\boldsymbol D_s\to \boldsymbol D_i$ by 
$$( \psi_{s,i} x)(t)=x(t)1_{\{t\le i-1\}}+(i-t)x(t)1_{\{i-1<t\le i\}}.
$$ 
\end{definition}
\begin{exercise}
 Show that the mapping $\psi_{s,i}$ is continuous.
\end{exercise}

\begin{lemma}\label{L16.3}
  A necessary and sufficient condition for $P_n\Rightarrow P$ on $\boldsymbol D_\infty$ is that $P_n\psi_k^{-1}\Rightarrow P\psi_k^{-1}$ on $\boldsymbol D_k$ for every $k\in\boldsymbol N$.
\end{lemma}
Proof. Since $\psi_k$ is continuous, the necessity follows from the mapping theorem.

For the sufficiency we need the isometry $\psi$ from Definition \ref{deep} and the inverse isometry $\psi^{-1}$:
\[\boldsymbol D_\infty\stackrel{\psi_k}{\to}\boldsymbol D_k,\qquad \boldsymbol D_\infty\stackrel{\psi}{\to}\overline {\boldsymbol D},\qquad \overline {\boldsymbol D}_\infty\stackrel{\psi^{-1}}{\to}\boldsymbol D_\infty.\]
Define two more mappings 
$$\overline {\boldsymbol D}\stackrel{\zeta_k}{\to}\boldsymbol D_1\times\ldots \times\boldsymbol D_k,\qquad \boldsymbol D_k\stackrel{\chi_k}{\to}\boldsymbol D_1\times\ldots \times\boldsymbol D_k$$
by 
\begin{align*}
 \zeta_k(\bar x)=(\bar x_1,\ldots,\bar x_k),\qquad \chi_k(x)=(\psi_{k,1}x,\ldots,\psi_{k,k}x).
\end{align*}

Consider the Borel $\sigma$-algebra $\overline {\mathcal D}$ for $(\overline {\boldsymbol D},\rho)$ and let $\overline {\mathcal D}_f\subset \overline {\mathcal D}$ be the class of sets of the form $\zeta_k^{-1}H$ where  $k\ge1$ and $H\in\mathcal D_1\times\ldots \times\mathcal D_k$, see Definition \ref{dpr}. The remainder of the proof is split into four steps.

Step 1. Applying Theorem \ref{2.4} show that $\mathcal  A=\overline {\mathcal D}_f$ is a convergence-determining class.
Given a ball $B(\bar x,\epsilon)\subset \overline {\boldsymbol D}$, take $k$ so that $2^{-k}<\epsilon/2$ and consider the cylinder sets
\[A_\eta=\{\bar y\in\overline {\boldsymbol D}: d_i^\circ(\bar x_i,\bar y_i)<\eta,i=1,\ldots,k\}\quad\mbox{for }0<\eta<\epsilon/2.\]
Then $\bar x\in A_\eta^\circ=A_\eta\subset B(\bar x,\epsilon)$ implies $A_\eta\in \mathcal  A_{x,\epsilon}$.
It remains to see that the boundaries of $A_\eta$ for different $\eta$ are disjoint.

Step 2. For probability measures $Q_n$ and $Q$ on $\overline {\boldsymbol D}$ show that   if  $Q_n\zeta_k^{-1}\Rightarrow Q\zeta_k^{-1}$ for every $k$, then $Q_n\Rightarrow Q$.  

This follows from the equality $\partial (\zeta_k^{-1}H)=\zeta_k^{-1}\partial H$ for $H\in \mathcal D_1\times\ldots \times\mathcal D_k$, see the proof of Theorem \ref{E2.4}. 

Step 3. Assume that $P_n\psi_k^{-1}\Rightarrow P\psi_k^{-1}$ on $\boldsymbol D_k$ for every $k$ and show that $P_n\psi^{-1}\Rightarrow P\psi^{-1}$ on $\overline {\boldsymbol D}$.

The map $\chi_k$ is continuous: if $x_n\to x$ in $\boldsymbol D_k$, then $\psi_{k,i}x_n\to \psi_{k,i}x$ in $\boldsymbol D_i$, $i\le k$. By the mapping theorem, 
$P_n\psi_k^{-1}\chi_k^{-1}\Rightarrow P\psi_k^{-1}\chi_k^{-1}$, and since $\chi_k\psi_k=\zeta_k\psi$, we get $P_n\psi^{-1}\zeta_k^{-1}\Rightarrow P\psi^{-1}\zeta_k^{-1}$. Referring to step 2 we conclude $P_n\psi^{-1}\Rightarrow P\psi^{-1}$.

Step 4. Show that $P_n\psi^{-1}\Rightarrow P\psi^{-1}$ on $\overline {\boldsymbol D}$ implies $P_n\Rightarrow P$ on $\boldsymbol D_\infty$.

Extend the isometry $\psi^{-1}$ to a map $\eta:\overline {\boldsymbol D}\to\boldsymbol D_\infty$ by putting $\eta(\bar x)=x_0\in\boldsymbol D_\infty$ for all
$\bar x\notin\overline {\boldsymbol D}_\infty$. Then $\eta$ is continuous when restricted to $\overline {\boldsymbol D}_\infty$, and since $\overline {\boldsymbol D}_\infty$ supports $P\psi^{-1}$ and the $P_n\psi^{-1}$, it follows that
\[P_n=P_n\psi^{-1}\eta^{-1}\Rightarrow P\psi^{-1}\eta^{-1}=P.\]

\begin{definition}
 For a probability measure $P$ on $\boldsymbol D_\infty$ define $T_P\subset[0,\infty)$ as the set of $t$ for which $P(J_t)=0$, where $J_t=\{x:x\mbox{ is discontinuous at }t\}$.  (See Lemma \ref{13.0}.)
\end{definition}
\begin{exercise}
 Let $P$ be the probability measure on $\boldsymbol D_\infty$  generated by the Poisson process with parameter $\lambda$. Show that $T_P=[0,\infty)$.
\end{exercise}
\begin{lemma}\label{p174}
 For $x\in \boldsymbol D_\infty$ let $r_t x$ be the restriction of $x$ on $[0,t]$. The function $r_t : \boldsymbol D_\infty\to  \boldsymbol D_t$ is measurable. The set of points at which $r_t$ is discontinuous belongs to $J_t$.
\end{lemma}
Proof. Denote $\delta_k=t/k$. Define the function $r^k_tx\in \boldsymbol D_t$ as having the value $x({i\delta_k})$ on $[{i\delta_k},{(i+1)\delta_k})$ for $0\le i\le k-1$ and the value $x(t)$ at $t$. Since the $\pi_{i\delta_k}$ are measurable $\mathcal D_\infty/\mathcal R_1$, it follows as in the proof of Theorem \ref{12.5} (b) that $r^k_t$ is measurable $\mathcal D_\infty/\mathcal D_t$. By Lemma \ref{L12.3},
\[d_t(r^k_tx,r_tx)\le \delta_k\vee w'_t(x,\delta_k)\to0\mbox{ as }k\to\infty\mbox{ for each  }x\in \boldsymbol D_\infty.\]
Now, to show that $r_t$ is measurable take a closed $F\in\mathcal D_t$. We have $F=\cap_\epsilon F^{2\epsilon}$, where the intersection is over positive rational $\epsilon$. From 
\[r_t^{-1}F\subset \liminf_k(r^k_t)^{-1}F^{\epsilon}=\bigcup_{j=1}^\infty\bigcap_{k=j}^\infty (r^k_t)^{-1}F^{\epsilon}\subset r_t^{-1}F^{2\epsilon}\]
we deduce that $r_t^{-1}F= \cap_\epsilon \liminf_k(r^k_t)^{-1}F^{\epsilon}$ is measurable. Thus $r_t$ is measurable.

To prove the second assertion take an $x\in \boldsymbol D_\infty$ which is continuous at $t$. If $d_\infty(x_n,x)\to0$, then by Theorem \ref{16.2}, 
$$d_t(r_tx_n,r_tx)=d_t(x_n,x)\to0.$$
In other words, if $x\notin J_t$, then  $r_t$ is continuous at  $x$.

\begin{theorem}\label{16.7}
 A necessary and sufficient condition for $P_n\Rightarrow P$ on $\boldsymbol D_\infty$ is that $P_nr_t^{-1}\Rightarrow Pr_t^{-1}$ for each $t\in T_P$.
\end{theorem}
Proof. If $P_n\Rightarrow P$ on $\boldsymbol D_\infty$, then $P_nr_t^{-1}\Rightarrow Pr_t^{-1}$ for each $t\in T_P$ due to the mapping theorem and Lemma \ref{p174}. 

For the reverse implication, it is enough, by Lemma \ref{L16.3}, to show that $P_n\psi_i^{-1}\Rightarrow P\psi_i^{-1}$ on $\boldsymbol D_i$ for every $i$. 
Given an $i$ choose a $t\in T_P$ so that $t\ge i$. %\[
%(\tau_ix)(t)=\left\{
%\begin{array}{ll}
%(\psi_ix)(t)  &  \mbox{for }   t\in[0, i],\\
% 0 &    \mbox{for }  t>i
%\end{array}
%\right.
%\]
Since $\psi_i= \psi_{t,i} \circ r_t$, the mapping theorem gives
\[P_n\psi_i^{-1}=(P_nr_t^{-1}) \psi_{t,i}^{-1}\Rightarrow (Pr_t^{-1}) \psi_{t,i}^{-1}= P\psi_i^{-1}.\]

\begin{exercise}
 Let $W^\circ$ be the standard Brownian bridge. For $t\in[0,\infty)$ put $W_t=(1+t)W^\circ_{t\over1+t}$. Show that such defined random element $W$ of $\boldsymbol D_\infty$ is a Gaussian process with zero means and covariance function $\mathbb E(W_sW_t)=s$ for $0\le s\le t<\infty$. This is a Wiener process $W=(W_t,0\le t<\infty)$. Clearly, $r_t(W)$ is a Wiener process which is a random element  of $\boldsymbol D_t$.
\end{exercise}

\begin{corollary}
 Let $\xi_1,\xi_2,\ldots$ be iid r.v. defined on $(\Omega,\mathcal F,\mathbb P)$. If $\xi_i$ have zero mean and variance $\sigma^2$ and  $X^n_t={\xi_1+\ldots+\xi_{\lfloor nt\rfloor}\over\sigma\sqrt n}$, then $X^n\Rightarrow W$  on $\boldsymbol D_\infty$.
\end{corollary}
Proof. By Theorem \ref{14.1}, $X^n\Rightarrow W$  on $\boldsymbol D_1$. The same proof gives $X^n\Rightarrow W$  on $\boldsymbol D_t$ for each $t\in [0,\infty)$. In other words, $r_t(X^n)\Rightarrow r_t(W)$ for each $t\in [0,\infty)$, and it remains to apply Theorem \ref{16.7}.

\begin{corollary}
 Suppose for each $n$, $\xi_{n1},\ldots,\xi_{nn}$ are iid indicator r.v. with $\mathbb P(\xi_{ni}=1)=\alpha/n$. If $X^n_t=\sum_{i\le nt} \xi_{ni}$, then  $X^n\Rightarrow X$ on $\boldsymbol D_\infty$, where $X$ is the Poisson process with parameter $\alpha$.
\end{corollary}
Proof. Combine Corollary \ref{E12.3} and Theorem \ref{16.7}.

\end{document}